\documentclass[12pt]{article} 
\usepackage{latexsym, %showkeys,
amscd, 
graphicx, epsfig, amssymb}

\newtheorem{thm}{Theorem}[section]
\newtheorem{defn}[thm]{Definition}
\newtheorem{prop}[thm]{Proposition}
\newtheorem{cor}[thm]{Corollary}
\newtheorem{lemma}[thm]{Lemma}
\newtheorem{rema}[thm]{Remark}

\newcommand{\halmos}{\rule{1ex}{1.4ex}}

\newcommand{\binom}[2]{{{#1}\choose {#2}}}
\newcommand{\text}[1]{\mbox{\rm #1}}

\newcommand{\nn}{\nonumber \\}

 \newcommand{\res}{\mbox{\rm Res}}

\newcommand{\edo}{\mbox{\rm End}\ }
 \newcommand{\pf}{{\it Proof.}\hspace{2ex}}
 \newcommand{\epf}{\hspace*{\fill}\mbox{$\halmos$}}
 \newcommand{\epfv}{\hspace*{\fill}\mbox{$\halmos$}\vspace{1em}}

\newcommand{\wt}{\mbox{\rm wt}\ }

\newcommand{\lwt}{\mbox{\rm wt}^{L}\;}
\newcommand{\rwt}{\mbox{\rm wt}^{R}\;}
\newcommand{\slwt}{\mbox{\rm {\scriptsize wt}}^{L}\,}
\newcommand{\srwt}{\mbox{\rm {\scriptsize wt}}^{R}\,}

\newcommand{\lbar}{\bigg\vert}

\newcommand{\A}{\mathcal{A}}
\newcommand{\Y}{\mathcal{Y}}
\newcommand{\C}{\mathbb{C}}
\newcommand{\Z}{\mathbb{Z}}
\newcommand{\R}{\mathbb{R}}
\newcommand{\Q}{\mathbb{Q}}
\newcommand{\N}{\mathbb{N}}

\newcommand{\V}{\mathcal{V}}
\newcommand{\one}{\mathbf{1}}
\newcommand{\BY}{\mathbb{Y}}
\newcommand{\ds}{\displaystyle}

\title{ {\bf Full field algebras} }
\date{}
\author{Yi-Zhi Huang and Liang Kong}
\begin{document}

\bibliographystyle{alpha}
\maketitle

\begin{abstract} 
We introduce a notion of full field algebra which is essentially
an algebraic formulation of the notion of genus-zero full
conformal field theory. For any vertex operator algebras $V^{L}$
and $V^{R}$, $V^{L}\otimes V^{R}$ is naturally a full field algebra
and we introduce a notion of full field algebra over $V^{L}\otimes V^{R}$. 
We study the structure of full field algebras over $V^{L}\otimes
V^{R}$ using modules and intertwining operators for $V^{L}$ and $V^{R}$.
For a simple vertex operator algebra $V$ 
satisfying certain natural finiteness and
reductive conditions needed 
for  the Verlinde conjecture to hold, we construct 
a bilinear form on the space of intertwining operators for $V$ and 
prove the nondegeneracy and other basic properties of this form. 
The proof of the nondegenracy of the bilinear form depends not only 
on the theory of intertwining operator algebras but also on the 
modular invariance for intertwining operator algebras through the 
use of the results obtained in the proof of the Verlinde conjecture
by the first author. Using 
this nondegenerate bilinear form, we construct a full field algebra
over $V\otimes V$ and an invariant bilinear form on this algebra. 
\end{abstract}

\renewcommand{\theequation}{\thesection.\arabic{equation}}
\renewcommand{\thethm}{\thesection.\arabic{thm}}
\setcounter{equation}{0}
\setcounter{thm}{0}
\setcounter{section}{-1}

\section{Introduction}

In the present paper, we
solve the problem of constructing a genus-zero full conformal field
theory (a conformal field theory on genus-zero Riemann surfaces
containing both chiral and antichiral
parts) from representations of a simple vertex operator algebra $V$
satisfying the following conditions: 
(i) $V_{(n)}=0$ for $n<0$, $V_{(0)}=\mathbb{C}\mathbf{1}$,
and $W_{(0)}=0$ for any irreducible $V$-module $W$ which is not equivalent
to $V$. %$V'$ as a $V$-module is isomorphic to $V$.
(ii) Every $\mathbb{N}$-gradable weak $V$-module is completely
reducible. (iii) $V$ is $C_{2}$-cofinite. %Note that if 
%$V$ satisfies (i)' $V_{(n)}=0$ for $n<0$, $V_{(0)}=\mathbb{C}\mathbf{1}$,
%$W_{(0)}=0$ for any irreducible $V$-module $W$ which is not equivalent
%to $V$, then (i) is satisfied.  
Note %also 
that the last two conditions are equivalent to 
a single condition that every weak $V$-module is completely reducible
(see \cite{L} and
\cite{ABD}).

Conformal field theories in its original form, as formulated 
by Belavin, Polyakov and Zamolodchikov \cite{BPZ} and 
by Kontsevich and Segal \cite{S1} \cite{S2} \cite{S3}, 
have both chiral and antichiral parts. The fundamental work \cite{MS1} \cite{MS2}
of Moore and Seiberg is also based on the existence of 
such full conformal field theories with both chiral and antichiral 
parts. 
In mathematics, however, it is mostly chiral 
conformal field theories that are
constructed and studied.  
To use conformal field theory to solve mathematical problems and to 
understand mathematical results such as
mirror symmetry, we need full conformal field theories, not 
just chiral or antichiral ones. 

In the case of conformal field theories associated to tori, 
Tsukada first constructed and studied these theories mathematically 
in his Ph. D. thesis under the direction of I. Frenkel (see \cite{Ts}). 
Assuming the existence of the structure of a modular tensor category on
the category of modules for a vertex operator algebra,  the existence
of conformal blocks with monodromies compatible with the modular tensor
category and all the necessary convergence properties, 
Felder, Fr\"{o}hlich, Fuchs and Schweigert \cite{FFFS} and 
Fuchs, Runkel, Schweigert and Fjelstad
\cite{FRS1} \cite{FRS2} \cite{FRS3}
\cite{FFRS}
studied open-closed
conformal field theories (in particular full (closed) conformal 
field theories) using the theory of tensor categories and
three-dimensional topological field theories. They constructed 
correlation functions as states in some three-dimensional topological 
field theories and they showed
the existence of consistent operator product expansion coefficients for
bulk operators. However, since these works were based on the fundamental
assumptions mentioned above, an explicit
construction of the corresponding full conformal field theories,
even in the genus-zero case,  is still
needed.

In \cite{KO}, Kapustin and Orlov studied full conformal field theories
associated to tori. They introduced a notion of vertex algebra which is
more general than the original notion of vertex algebra \cite{B}
or vertex
operator algebra \cite{FLM2} 
by allowing both chiral and antichiral parts. In
\cite{R1} and \cite{R2}, Rosellen studied these algebras in details.
However, the construction of the full conformal field theories
associated to affine Lie algebras (the WZNZ models) and to the Virasoro
algebra (the minimal models) was still an open problem.  More generally,
we would like to construct full conformal field theories from the
representations of a vertex operator algebra satisfying reasonable
conditions.  Also, since the braid group representations obtained from
the representations of these vertex operator algebras are not one
dimensional in general, it seems that the corresponding full conformal
field theories in general do not satisfy the axioms for the algebras
introduced and studied in \cite{KO}, \cite{R1} and \cite{R2}. 

In \cite{H5}, \cite{H6} and \cite{H8}, the first author
constructed genus-zero chiral theories, genus-one chiral theories and
modular tensor categories from the representations of simple
vertex operator algebras satisfying the three conditions above. 
Since modular tensor categories give modular functors (see \cite{T} and 
\cite{BK}), these results also give modular functors.
One of the remaining problems is to
construct a full conformal field theory from a 
chiral theory and an antichiral theory obtained from the 
chiral theory. In the present paper,
we solve this problem in the genus-zero case by constructing 
a full conformal field theory corresponding to what physicists
call a diagonal theory (see, for example, \cite{MS3})
from the representations of a simple vertex 
operator algebra satisfying the conditions above. The 
same full conformal field theory is also constructed by the 
second author using the theory of tensor categories in \cite{K}. 
The genus-one case and the higher-genus case can be obtained 
from the construction
of the genus-zero theories in this paper and the properties of
genus-one and higher-genus chiral theories. These will be discussed in 
future publications. 

Technically, we construct a genus-zero full conformal 
field theory as follows: We first 
introduce a notion of full field algebra and several variants, which 
are essentially algebraic formulations of 
the notion of genus-zero full conformal field theory (for a 
precise discussion of the equivalence of this notion of 
full field algebra and its variants with geometric formulations
of genus-zero conformal field theories in terms of operads, see \cite{K}).
For a simple vertex operator algebra satisfying the three conditions above,
by the results in \cite{H5}, 
we have an intertwining operator algebra,  which is 
equivalent to a genus-zero chiral conformal field theory
(see \cite{H3} and \cite{H3.5}). The genus-zero chiral conformal field theory
also gives a genus-zero antichiral conformal field theory. 
We construct a nondegenerate bilinear form on the space of intertwining 
operators and use this bilinear form to 
put the genus-zero chiral and antichiral conformal field
theories together. We show that the resulting mathematical object is 
a full field algebra satisfying additional properties 
and thus gives a genus-zero full conformal field theory. 
One interesting aspect of our construction is that our construction
(actually the proof of the nondegeneracy of the bilinear form on the space of 
the intertwining operators)
needs the theorem proved in \cite{H7} (see also \cite{H9})
stating that the Verlinde conjecture holds for 
such a vertex operator algebra. This theorem in \cite{H7},
and thus also our construction of genus-zero full conformal field theories, 
depend not only on 
genus-zero chiral theories constructed in \cite{H5},
but also on genus-one chiral theories constructed  in \cite{H6}.

This paper is organized as follows: In Section 1, we introduce the notion 
of full field algebra and several variants
and discuss their basic properties. In Section 2,
we discuss basic relations between intertwining operator algebras
and full field algebras. This is a section preparing for 
our construction in Section 3. Our construction of full field algebras
is given in Section 3. We also construct invariant bilinear  forms
on these full field algebras in the same section.

\paragraph{Acknowledgment} 
The first author is partially supported 
by NSF grant DMS-0401302.

\renewcommand{\theequation}{\thesection.\arabic{equation}}
\renewcommand{\thethm}{\thesection.\arabic{thm}}
\setcounter{equation}{0}
\setcounter{thm}{0}

\section{Definitions and basic properties}

Let $\mathbb{F}_{n}(\C) =\{ (z_{1}, \ldots, z_{n})\in \C^n \;|\; z_i\neq z_j
\mbox{ if } i\neq j \}$. 
For an $\R$-graded vector space $F=\coprod_{r\in \R}F_{(r)}$,
we let $\overline{F}=\prod_{r\in \R}F_{(r)}$ be the algebraic completion
of $F$. For $r\in \R$, let 
$P_{r}$ be the projection from $F$ or $\overline{F}$ to $F_{(r)}$.
A series $\sum f_{n}$ in $\overline{F}$ is said to be {\it
absolutely convergent} if for any $f'\in F'$, $\sum|\langle f', f_{n}\rangle|$
is convergent. The sums $\sum|\langle f', f_{n}\rangle|$ for $f'\in 
F'$ define a linear functional on $F'$. We call this linear functional 
the {\it sum} of the series and denote it by the same notation 
$\sum f_{n}$. If the homogeneous subspaces of $F$ are all finite-dimensional, 
then $\overline{F}=(F')^{*}$ and, in this case, the sum 
of an absolutely convergent series is always in $\overline{F}$. 
When the sum is in $\overline{F}$,
we say that the series is {\it absolutely convergent in $\overline{F}$}. 

\begin{defn} \label{ffa}
{\rm A {\it full field algebra} is an $\R$-graded vector space 
$F=\coprod_{r\in\R}F_{(r)}$ (graded by {\it total conformal weight}
or simply {\it total weight}),  
equipped with  {\it correlation function maps}
$$\begin{array}{rrcl}
m_{n}: & F^{\otimes n}\times
\mathbb{F}_{n}(\C)&\to& \overline{F}\\
&(u_{1}\otimes \cdots \otimes u_{n}, (z_{1}, \dots, z_{n}))&\mapsto&
m_{n}(u_{1}, \dots, u_{n}; z_{1}, \bar{z}_{1}, \dots, z_{n}, \bar{z}_n),
\end{array}$$
for $n\in \Z_+$ and
a distinguished element $\mathbf{1}$ called {\it vacuum}
satisfying the following axioms:

\begin{enumerate}

\item For $n\in \Z_+$, 
$m_n(u_{1}, \dots, u_{n}; z_{1}, \bar{z}_{1}, \dots, z_{n}, \bar{z}_n)$ is 
linear in $u_{1}, \dots, u_{n}$ and smooth in the real and imaginary 
parts of $z_{1}, \dots, z_{n}$.

\item For $u\in F$, $m_{1}(u; 0, 0)=u$. 

\item For $n\in \Z_{+}$, $u_{1}, \dots, u_{n}\in F$, 
\begin{eqnarray*}
&m_{n+1}(u_{1}, \dots, u_{n}, \one; z_{1}, \bar{z}_{1}, \dots, z_{n}, \bar{z}_{n}, 
z_{n+1}, \bar{z}_{n+1})&\nn
&=m_n(u_{1}, \dots, u_{n}; z_{1}, \bar{z}_{1}, \dots, z_{n}, \bar{z}_{n}).&
\end{eqnarray*}

\item The {\it  convergence property}:   
For $k, l_{1}, \ldots, l_{k} \in \Z_+$ and 
$u_{1}^{(1)},\ldots, u_{l_{1}}^{(1)},\ldots, u_{1}^{(k)}$, $\dots,
u_{l_{k}}^{(k)} \in F$,  the series 
\begin{eqnarray}\label{ffa-conv-axiom}
\lefteqn{\sum_{r_{1}, \ldots, r_k} m_{k}(P_{r_{1}} m_{l_{1}}(u_{1}^{(1)}, 
\ldots, u_{l_{1}}^{(1)}; z_{1}^{(1)}, \bar{z}_{1}^{(1)}, \dots, 
z_{l_{1}}^{(1)}, \bar{z}_{l_{1}}^{(1)} ), \dots,}  \nn
&& 
P_{r_k} m_{l_{k}} (u_{1}^{(k)}, \ldots, u_{l_{k}}^{(k)}, 
z_{1}^{(k)}, \bar{z}_{1}^{(k)}, \ldots, z_{l_{k}}^{(k)}, 
\bar{z}_{l_{k}}^{(k)}); z^{(0)}_{1}, \bar{z}_{1}^{(0)}, \ldots, 
z^{(0)}_{k}, \bar{z}^{(0)}_k)\nn
&&
\end{eqnarray}
converges absolutely to  
\begin{eqnarray}\label{sum}
\lefteqn{m_{l_{1}+ \cdots + l_k}(u_{1}^{(1)}, \dots, u_{l_k}^{(k)};
z_{1}^{(1)}+z^{(0)}_{1}, \bar{z}_{1}^{(1)}+ \bar{z}^{(0)}_{1}, \dots,
z_{l_{1}}^{(1)}+z^{(0)}_{1}, }  \nn
&&\quad
\bar{z}_{l_{1}}^{(1)}+\bar{z}^{(0)}_{1},  
\dots, 
z_{1}^{(k)}+z^{(0)}_k,  
\bar{z}_{1}^{(k)}+ \bar{z}^{(0)}_k, \dots,
z_{l_k}^{k}+z^{(0)}_k, \bar{z}_{l_k}^{k}+\bar{z}^{(0)}_k).  \nn
&&
\end{eqnarray}
when $|z_p^{(i)}| + |z_q^{(j)}|< |z^{(0)}_i
-z^{(0)}_j|$ for $i,j=1, \ldots, k$, $i\ne j$ and for $p=1, 
\dots,  l_i$ and $q=1, \dots, l_j$.

\item The {\it permutation property}: For any $n\in \Z_{+}$ and 
any $\sigma\in S_{n}$,
we have 
\begin{eqnarray}
\lefteqn{m_{n}(u_{1}, \dots, u_{n}; z_{1}, \bar{z}_{1}, 
\dots, z_{n}, \bar{z}_{n})} \nn
&& = m_{n}(u_{\sigma(1)}, \dots, u_{\sigma(n)}; 
z_{\sigma(1)}, \bar{z}_{\sigma(1)}, \dots, 
z_{\sigma(n)}, \bar{z}_{\sigma(n)})   \label{ffa-perm-axiom}
\end{eqnarray}
for $u_{1}, \dots, u_{n}\in F$ and $(z_{1}, \dots, z_{n})\in 
\mathbb{F}_{n}(\C)$.

\item Let $\mathbf{d}$ be the 
grading operator, that is, the operator defined by $\mathbf{d}f=rf$ 
for $f\in F_{(r)}$. Then for $n\in \Z_{+}$, $a\in \R$, $u_{1}, 
\dots, u_{n}\in F$,
\begin{eqnarray*}
\lefteqn{e^{a\mathbf{d}}m_{n}(u_{1}, \dots, u_{n}; z_{1}, \bar{z}_{1},
\dots, z_{n}, \bar{z}_{n})}\nn
&&=m_{n}(e^{a\mathbf{d}}u_{1}, \dots, e^{a\mathbf{d}}u_{n}; 
e^{a}z_{1}, e^{a}\bar{z}_{1},
\dots, e^{a}z_{n}, e^{a}\bar{z}_{n}).
\end{eqnarray*}

\end{enumerate}
}

\end{defn}

We denote the full field algebra defined above by 
$(F, m, \mathbf{1})$ or simply by $F$.
In the definition above, we use the notations 
$$m_{n}(u_{1}, \dots, u_{n}; z_{1},
\bar{z}_{1} \dots, z_{n}, \bar{z}_n)$$ instead of
$$m_{n}(u_{1}, \dots,
u_{n}; z_{1},\dots, z_{n})$$
to emphasis that 
these are in general not holomorphic in $z_{1},\dots, z_{n}$.
For $u'\in F'$, $u_{1}, \dots, u_{n}\in F$, 
$$\langle u', m_{n}(u_{1}, \dots, u_{n}; z_{1},
\bar{z}_{1} \dots, z_{n}, \bar{z}_n)\rangle$$
as a function of $z_{1}, \dots, z_{n}$ 
is called a {\it correlation function}.
{\it Homomorphisms} and {\it isomorphisms} for full field algebras are
defined in the obvious way. 

\begin{rema}\label{weak-c-p}
{\rm Note that in the convergence property, we require that 
the multisum is absolutely convergent. This is stronger than 
the following convergence property:
For $k, l \in \Z_+$ and 
$u_{1},\ldots,  u_{k-1}, v_{1}, \dots, v_{l}\in F$,  the series 
\begin{eqnarray*}
\lefteqn{\sum_{r} m_{k}(u_{1}, \dots, u_{k-1}, 
P_{r} m_{l} (v_{1}, \dots, v_{l};
z_{1}^{(k)}, \bar{z}_{1}^{(k)}, \dots, z_{l}^{(k)}, 
\bar{z}_{l}^{(k)}); }  \nn
&& \quad\quad\quad\quad\quad\quad\quad\quad\quad\quad\quad\quad
\quad\quad\quad\quad\quad\quad z_{1}^{(0)}, \bar{z}_{1}^{(0)}, \ldots, 
z_{k}^{(0)}, \bar{z}_k^{(0)})\nn
&&
\end{eqnarray*}
converges absolutely to  
\begin{eqnarray*}
\lefteqn{m_{k+ l}(u_{1}, \dots, u_{k-1}, v_{1}, \dots, v_{l};
z^{(0)}_{1}, \bar{z}^{(0)}_{1}, \dots,
z^{(0)}_{k-1}, }  \nn
&&\quad
\bar{z}^{(0)}_{k-1},  
\dots, 
z_{1}^{(k)}+z^{(0)}_k,  
\bar{z}_{1}^{(k)}+ \bar{z}^{(0)}_k, \dots,
z_{l}^{k}+z^{(0)}_k, \bar{z}_{l}^{k}+\bar{z}^{(0)}_k).  \nn
&&
\end{eqnarray*}
when $z_{i}^{(0)}\ne z_{j}^{(0)}$ for $i,j=1, \ldots, k$ and 
$|z_p^{(k)}|< |z^{(0)}_k
-z^{(0)}_i|$ for $i=1, \ldots, k-1$,  and for $p=1, \dots, l$.
However, for the purpose of constructing genus-zero 
conformal field theory satisfying geometric axioms, 
this version of the convergence property is actually enough.}
\end{rema}

Let $(F, m, \mathbf{1})$ be a full field algebra and 
let 
$$\begin{array}{rrcl}
\mathbb{Y}:& F^{\otimes 2} \times \C^{\times} &\to&
  \overline{F}\\
&(u\otimes v, z, \bar{z})&\mapsto& \mathbb{Y}(u; z, \bar{z})v
\end{array}$$ 
be given by
$$\mathbb{Y}(u; z, \bar{z})v = m_{2}(u \otimes v; z,\bar{z}, 0, 0)$$
for $u, v\in F$. The map $\mathbb{Y}$ is called the {\it full vertex operator map}
and for $u\in F$, $\mathbb{Y}(u; z, \bar{z})$ is called the {\it 
full vertex operator} associated to $u$. 
We have the following immediate consequences of the definition:

\begin{prop}
\begin{enumerate}

\item The {\it identity property}: $\mathbb{Y}(\mathbf{1}; z, \bar{z})=I_{F}$.

\item The {\it creation property}: ${\ds \lim_{z\to 0}\mathbb{Y}(u; z, \bar{z})
\mathbf{1}=u}$.

\item For $f\in F$,
\begin{equation}\label{d-bracket}
\left[ \mathbf{d}, \mathbb{Y}(u; z, \bar{z})\right] 
= \left( z\frac{\partial}{\partial z} + \bar{z}\frac{\partial}{\partial
  \bar{z}} \right)   \mathbb{Y}(u; z, \bar{z})
+\mathbb{Y}(\mathbf{d} u; z,\bar{z}) 
\end{equation}

\item The total weight of the vacuum $\one$ is $0$, that is, 
$\mathbf{d}\one=0$.
\end{enumerate}
\end{prop}
\pf
For $u\in F$, 
\begin{eqnarray*}
\mathbb{Y}(\mathbf{1}; z, \bar{z})u&=&m_{2}(\one, u; z, \bar{z}, 
0, 0)\nn
&=&m_{2}(u, \one; 0, 0, z, \bar{z})\nn
&=&m_{1}(u; 0, 0)\nn
&=&u.
\end{eqnarray*}

For $u\in F$, 
\begin{eqnarray*}
\lim_{z\to 0}\mathbb{Y}(u; z, \bar{z})
\mathbf{1}&=&\lim_{z\to 0}m_{2}(u, \one; z, \bar{z}, 
0, 0)\nn
&=&\lim_{z\to 0}m_{1}(u; z, \bar{z})\nn
&=&m_{1}(u; 0, 0)\nn
&=&u.
\end{eqnarray*}

For $u, v\in F$ and $a\in \R$,
\begin{eqnarray}\label{d-bracket-0}
e^{a\mathbf{d}}\mathbb{Y}(u; z, \bar{z})v
&=&\mathbb{Y}(e^{a\mathbf{d}}u; e^{a}z, e^{a}\bar{z})e^{a\mathbf{d}}v.
\end{eqnarray}
Taking derivatives of both sides of (\ref{d-bracket-0}) with respect to 
$a$, letting $a=0$ and noticing that $v$ is arbitrary, we obtain (\ref{d-bracket}).

From the identity property, 
$$\left(z\frac{\partial}{\partial z}
+\bar{z}\frac{\partial}{\partial \bar{z}}\right)
\Y(\one; z, \bar{z})=0.$$
Then by (\ref{d-bracket}), we obtain
\begin{equation}\label{d-bracket-1}
\left[\mathbf{d}, \mathbb{Y}(\one; z, \bar{z})\right] 
=\mathbb{Y}(\mathbf{d} \one; z,\bar{z}).
\end{equation}
Since $\mathbb{Y}(\one; z, \bar{z})=I_{F}$, 
(\ref{d-bracket-1}) gives
\begin{equation}\label{d-bracket-2}
\mathbb{Y}(\mathbf{d} \one; z,\bar{z})=0.
\end{equation}
Applying (\ref{d-bracket-2}) to $\one$, taking the 
limit $z\to 0$ on both sides of the resulting formula and using 
the creation property, 
we obtain $\mathbf{d} \one=0$. So the total weight of $\one$ is $0$. 
\epfv

Now we discuss two important properties of full field algebras which follow also
immediately from the definition. 

\begin{prop}[Associativity]
For $u_{1}, u_{2}, u_{3}\in F$, 
\begin{equation}   \label{asso-1}
\mathbb{Y}(u_{1}; z_{1}, \bar{z}_{1})\mathbb{Y}(u_{2}; z_{2}, \bar{z}_{2})u_{3} = 
\mathbb{Y}(\mathbb{Y}(u_{1}; z_{1}-z_{2}, \bar{z}_{1}-\bar{z}_{2})u_{2}; z_{2},\bar{z}_{2})u_{3}
\end{equation}
when $|z_{1}|>|z_{2}|>|z_{1}-z_{2}|>0$. 
\end{prop}
\pf
The convergence property
says, in particular, that the series 
\begin{equation} \label{prod-series}
\mathbb{Y}(u_{1}; z_{1}, \bar{z}_{1})\mathbb{Y}(u_{2}; z_{2}, 
\bar{z}_{2})u_{3}=
\sum_{n\in \R} \mathbb{Y}(u_{1}; z_{1}, \bar{z}_{1}) P_n 
\mathbb{Y}(u_{2}; z_{2}, \bar{z}_{2})u_{3}, 
\end{equation} 
(a {\it product} of full vertex operators)
converges absolutely in $\overline{F}$ for $u_{1}, u_{2}, u_{3} \in F$ when 
$|z_{1}|>|z_{2}|>0$. The convergence
property also says, in particular, that 
the series  
\begin{eqnarray}  \label{it-series}
&&\mathbb{Y}(\mathbb{Y}(u_{1}; z_{1}-z_{2},\bar{z}_{1}-\bar{z}_{2})u_{2}, 
z_{2},\bar{z}_{2})u_{3}  \nn
&&\hspace{2cm}= \sum_{n\in \R}\mathbb{Y}(P_n 
\mathbb{Y}(u_{1}; z_{1}-z_{2},\bar{z}_{1}-\bar{z}_{2})u_{2}, z_{2},\bar{z}_{2})u_{3}
\end{eqnarray}
(an {\it iterate} of full vertex operators)
converges absolutely for $u_{1}, u_{2}, u_{3}\in F$ 
when $|z_{2}|>|z_{1}-z_{2}|>0$. 
Moreover, the convergence property also says that both (\ref{prod-series}) 
and (\ref{it-series}) converge absolutely to 
$$
m_{3}(u_{1}, u_{2}, u_{3}; z_{1}, \bar{z}_{1}, z_{2}, \bar{z}_{2}, 0, 0).
$$
This proves the associativity.
\epfv

\begin{prop}[Commutativity]\label{commu-general}
For $u_{1}, u_{2}, u_{3}\in F$, 
\begin{eqnarray}   
\mathbb{Y}(u_{1}; z_{1}, \bar{z}_{1})\mathbb{Y}(u_{2}; z_{2}, 
\bar{z}_{2})u_{3}, \label{prod-series-1}\\
\mathbb{Y}(u_{2}; z_{2}, \bar{z}_{2})\mathbb{Y}(u_{1}; z_{1}, 
\bar{z}_{1})u_{3}, \label{prod-series-2}
\end{eqnarray}
are the expansions of 
$$
m_{3}(u_{1}, u_{2}, u_{3}; z_{1}, \bar{z}_{1}, z_{2}, \bar{z}_{2}, 0, 0).
$$
in the sets given by 
$|z_{1}|>|z_{2}|>0$ and $|z_{2}|>|z_{1}|>0$, respectively. 
\end{prop}
\pf
By the convergence property, we know that  (\ref{prod-series-1}) 
and (\ref{prod-series-2}) converge absolutely to 
$$m_{3}(u_{1}, u_{2}, u_{3}; z_{1}, \bar{z}_{1}, z_{2}, \bar{z}_{2}, 0, 0)
$$
and 
$$
m_{3}(u_{2}, u_{1}, u_{3}; z_{2}, \bar{z}_{2}, z_{1}, \bar{z}_{1}, 0, 0),
$$
respectively, when $|z_{1}|>|z_{2}|>0$ and $|z_{2}|>|z_{1}|>0$, 
respectively. By the permutation property, 
$$m_{3}(u_{1}, u_{2}, u_{3}; z_{1}, \bar{z}_{1}, z_{2}, \bar{z}_{2}, 0, 0)
=m_{3}(u_{2}, u_{1}, u_{3}; z_{2}, \bar{z}_{2}, z_{1}, \bar{z}_{1}, 0, 0).$$
Thus (\ref{prod-series-1}) 
and (\ref{prod-series-2}) converge absolutely to 
$$m_{3}(u_{1}, u_{2}, u_{3}; z_{1}, \bar{z}_{1}, z_{2}, \bar{z}_{2}, 0, 0)
$$
when $|z_{1}|>|z_{2}|>0$ and $|z_{2}|>|z_{1}|>0$, 
respectively. So they are the expansions of 
$$
m_{3}(u_{1}, u_{2}, u_{3}; z_{1}, \bar{z}_{1}, z_{2}, \bar{z}_{2}, 0, 0).
$$
in the sets given by 
$|z_{1}|>|z_{2}|>0$ and $|z_{2}|>|z_{1}|>0$, respectively. 
\epfv

Before proving more properties, we would like to 
discuss first the problem of constructing full field algebras. 
It is clear that vertex operator algebras have structures of 
full field algebras. 
Let $(V^L, Y^L, \one^L, \omega^L)$ and 
$(V^R, Y^R, \one^R, \omega^R)$ be two vertex operator algebras. 
Consider the graded vector space $V^L\otimes V^R$ equipped with the 
correlation function maps, the vacuum and the operator $\mathbf{d}$ 
given as follows:
For $n\in \Z_{+}$, $u^{L}_{1}, \dots, u^{L}_{n}\in V^{L}$ and 
$u^{R}_{1}, \dots, u^{R}_{n}\in V^{R}$,
$m_{n}(u^{L}_{1}\otimes u^{R}_{1},
\dots, u^{L}_{n}\otimes u^{R}_{n}; z_{1}, \bar{z}_{1}, 
\dots, z_{n}, \bar{z}_{n})$ 
are given by the analytic extensions of
$$(Y^{L}(u^{L}_{1}, z_{1})\otimes Y^{R}(u^{R}_{1}, \bar{z}_{1}))
\cdots (Y^{L}(u^{L}_{n}, z_{n})\otimes Y^{R}(u^{R}_{n}, \bar{z}_{n}))
\mathbf{1}.$$
Then we take the vacuum $\mathbf{1}=\mathbf{1}^{L}\otimes \mathbf{1}^{R}$
and the operators $\mathbf{d}=
L^L(0)\otimes I_{V^R} + I_{V^L}\otimes L^R(0)$%,
%$D^{L}=L^{L}(-1)\otimes I_{V^{R}}$ and $D^{R}=I_{V^{L}}\otimes 
%L^{R}(-1)$
.
In particular, the full vertex operators are given by
$$
\mathbb{Y}(u^L\otimes u^R, z) v^L\otimes v^R 
= Y^L(u^L, z)v^L \otimes Y^R(u^R, \bar{z})v^R.
$$
for $u^L, v^L\in V^L, u^R, v^R\in V^R$ and $z\in \C^{\times}$. 

We have:

\begin{prop}
The vector space $V^L\otimes V^R$ equipped with the 
correlation function maps and the vacuum 
$\mathbf{1}$ given above is 
a full field algebra.
\end{prop}
\pf
The proof is a straightforward and easy verification.
\epfv

Note that there is also a vertex operator algebra structure
on $V^{L}\otimes V^{R}$. For simplicity, we shall 
use the notation $V^{L}\otimes V^{R}$ to denote both the 
vertex operator algebra and the full field algebra structure. 
It should be easy to see which structure we will be using 
in the remaining part of this paper.

The full field algebra $V^L\otimes V^R$ 
in general does not give a genus-one theory, that is, 
suitable $q$-traces, even in the case that they are convergent,
of the full vertex operators in general are not  modular 
invariant. For chiral theories, we know from
\cite{H6} 
that if we consider the intertwining operator
algebras constructed from irreducible 
modules for suitable vertex operator algebras, we do have
modular invariance. 
So it is then natural to look for full field algebras 
from suitable extensions of $V^L\otimes V^R$ by
$V^L\otimes V^R$-modules. 

Note that  $V^L\otimes V^R$ has an  
$\Z \times \Z$-grading with grading operators being
$L^L(0)\otimes I_{V^R}$ and $I_{V^L}\otimes L^R(0)$.  If a
full field algebra is an extension of $V^L\otimes V^R$
by $V^L\otimes V^R$-modules, it has an 
$\R \times \R$-grading. 

For any $\R \times \R$-graded vector space 
$F=\coprod_{(m,n)\in \R\times \R} F_{(m,n)}$, 
we have a {\it left  grading operator} $\mathbf{d}^L$ and a
{\it right grading operator}
$\mathbf{d}^R$ defined by $\mathbf{d}^L u=m u, \mathbf{d}^R u =nu$ 
for $u\in E_{(m,n)}$, where $m$ ($n$) is called the
left (right) weight of $u$ and is
denoted by $\text{\rm wt}^L u$ ($\text{\rm wt}^R u$). 
For $m, n\in \R$, let $P_{m,n}$ be the 
projection from $F \to F_{(m,n)}$. 
We still use $F'$ and $\overline{F}$ to denote the 
graded dual and the algebraic completion of $F$, but note that 
they are with respect to the $\R \times \R$-grading, not any $\R$-grading
induced from the $\R \times \R$-grading. 

\begin{defn}  \label{RR-grading-ffa}
{\rm An $\R \times \R$-graded full field algebra
is a full field algebra $(F, m, \mathbf{1})$ equipped with an 
$\R \times \R$-grading on $F$ (graded by {\it left conformal weight}
or {\it left weight} and {\it right conformal weight}
or {\it right weight} and thus equipped with left and right 
grading operators $\mathbf{d}^{L}$
and $\mathbf{d}^{R}$) and operators $D^{L}$ and $D^{R}$ 
satisfying the following conditions:
\begin{enumerate}

\item The {\it grading compatibility}: 
$\mathbf{d}=\mathbf{d}^{L}+\mathbf{d}^{R}$.

\item The {\it single-valuedness property}: 
$e^{2\pi i (\mathbf{d}^L-\mathbf{d}^R)} = I_F$. 

\item The {\it convergence property}: For $k, l_{1}, \ldots, 
l_{k} \in \Z_+$ and 
$u_{1}^{(1)},\ldots, u_{l_{1}}^{(1)},\ldots, u_{1}^{(k)}$, $\dots,
u_{l_{k}}^{(k)} \in F$,  the series 
\begin{eqnarray}   \label{series-case-2-0}
\lefteqn{\sum_{p_1,q_1, \dots, p_k, q_k} 
m_k(P_{p_1,q_1}m_{l_1}(u_1^{(1)},\dots, u_{l_1}^{(1)}; 
z_1^{(1)}, \bar{z}_1^{(1)}, 
\dots, z_{l_1}^{(1)}, \bar{z}_{l_1}^{(1)}), \dots,} \nn
&&\quad 
P_{p_k,q_k}m_{l_k}(u_1^{(k)},\dots, u_{l_k}^{(k)}; z_1^{(k)}, 
\bar{z}_1^{(k)}, 
\dots, z_{l_k}^{(k)}, \bar{z}_{l_k}^{(k)}); 
z_1^{(0)}, \bar{z}_1^{(0)}, \dots, z_k^{(0)}, \bar{z}_k^{(0)} )\nn
&&
\end{eqnarray}
converges absolutely to  (\ref{sum})
when $|z_p^{(i)}| + |z_q^{(j)}|< |z^{(0)}_i
-z^{(0)}_j|$ for $i,j=1, \ldots, k$, $i\ne j$ and for $p=1, 
\dots,  l_i$ and $q=1, \dots, l_j$.

\item {\it The $\mathbf{d}^L$- and $\mathbf{d}^R$-bracket properties}: 
For $u\in F$, 
\begin{eqnarray}
\left[ \mathbf{d}^L, \mathbb{Y}(u; z,\bar{z})\right] 
&=& z\frac{\partial}{\partial z} \mathbb{Y}(u; z, \bar{z})
+\mathbb{Y}(\mathbf{d}^L u; z,\bar{z})     \label{d-l} \\
\left[ \mathbf{d}^R, \mathbb{Y}(u; z, \bar{z})\right]
&=& \bar{z}\frac{\partial}{\partial \bar{z}} \mathbb{Y}(u; z,\bar{z})
+ \mathbb{Y}(\mathbf{d}^R u; z, \bar{z}). \label{d-r}
\end{eqnarray}

\item The {\it $D^{L}$- and $D^{R}$-derivative property}: For $u\in F$,
\begin{eqnarray}
&&\left[ D^L, \mathbb{Y}(u; z,\bar{z})\right] = 
\mathbb{Y}(D^{L}u; z,\bar{z}) =
\frac{\partial}{\partial z}\mathbb{Y}(u; z, \bar{z}),\\
&&\left[ D^R, \mathbb{Y}(u; z, \bar{z})\right] = 
\mathbb{Y}(D^{R}u; z, \bar{z}) =
\frac{\partial}{\partial \bar{z}}\mathbb{Y}(u; z, \bar{z}). \label{D-L-R-b-d}
\end{eqnarray}
\end{enumerate}}
\end{defn}

We denote the $\R \times \R$-graded 
full field algebra defined above by 
$(F, m, \mathbf{1}, D^{L}, D^{R})$ or simply by $F$. But note that 
there is a refined grading on $F$ now.

\begin{rema}  {\rm
Note that for $\R \times \R$-graded 
full field algebra, there is also a weaker convergence property 
similar to the one in Remark \ref{weak-c-p}.}
\end{rema}

\begin{rema}  {\rm
The single-valuedness property actually says that $F$ is  graded by 
a subgroup $\{(m,n)\in \R \times \R\; |\; m-n \in \Z\}$ of $\R \times \R$. 
This single-valuedness indeed corresponds to a certain single-valuedness
condition
in the geometric axioms for full conformal field theories. }
\end{rema}

We have the following immediate consequences of the definition above.

\begin{prop}
\begin{enumerate}

\item The pair (left weight, right weight) for $\one$ is $(0, 0)$,
that is, $\mathbf{d}^{L}\one=\mathbf{d}^{R}\one=0$.

\item The pairs (left weight, right weight) for $D^{L}$ and $D^{R}$ are 
$(1, 0)$ and $(0, 1)$, respectively, that is, 
\begin{eqnarray*}
[\mathbf{d}^{L}, D^{L}]&=&D^{L},\\
{[\mathbf{d}^{R}, D^{L}]}&=&0,\\
{[\mathbf{d}^{L}, D^{R}]}&=&0,\\
{[\mathbf{d}^{R}, D^{R}]}&=&D^{R}.
\end{eqnarray*}

\item $D^{L}\one=D^{R}\one =0$. 
\end{enumerate}
\end{prop}
\pf
From the identity property, 
$$z\frac{\partial}{\partial z}
\mathbb{Y}(\one; z, \bar{z})=0.$$
Then by (\ref{d-l}), we obtain
\begin{equation}\label{d-l-1}
\left[\mathbf{d}^{L}, \mathbb{Y}(\one; z, \bar{z})\right] 
=\mathbb{Y}(\mathbf{d}^{L} \one; z,\bar{z}).
\end{equation}
Since $\mathbb{Y}(\one; z, \bar{z})=I_{F}$, we obtain 
\begin{equation}\label{d-l-2}
\mathbb{Y}(\mathbf{d}^{L} \one; z,\bar{z})=0
\end{equation}
from (\ref{d-l-1}).
Applying (\ref{d-l-2}) to $\one$, 
taking the limit $z\to 0$ on both sides of (\ref{d-l-2}) and using 
the creation property, 
we obtain $\mathbf{d}^{L} \one=0$. So the left weight of $\one$ is $0$.
Similarly, we can prove that the right weight of $\one$ is $0$. 

Applying both sides of (\ref{d-l}) to $\one$, taking the 
limit $z\to 0$ and then using the creation property and the fact 
$\mathbf{d}^{L} \one =0$
we have just proved, we obtain
\begin{equation}\label{limit-z-der}
\lim_{z\to 0}z\frac{\partial}{\partial z}\mathbb{Y}(u; z, \bar{z})\one=0
\end{equation}
for $u\in F$. 
Applying $\frac{\partial}{\partial z}$ to both sides of 
(\ref{d-l}) and using the 
$D^{L}$-derivative property,
we obtain
\begin{eqnarray}\label{d-l-3}
\lefteqn{\left[\mathbf{d}^{L}, \mathbb{Y}(D^{L}u; z, \bar{z})\right]}\nn
&&= \frac{\partial}{\partial z}z\frac{\partial}{\partial z}
\mathbb{Y}(u; z, \bar{z})
+\mathbb{Y}(D^{L}\mathbf{d}^{L} u; z,\bar{z})\nn
&&=z\frac{\partial}{\partial z}
\frac{\partial}{\partial z}\mathbb{Y}(u; z, \bar{z})
+\frac{\partial}{\partial z}\mathbb{Y}(u; z, \bar{z})
+\mathbb{Y}(D^{L}\mathbf{d}^{L} u; z,\bar{z})\nn
&&=z\frac{\partial}{\partial z}
\mathbb{Y}(D^{L}u; z, \bar{z})
+\mathbb{Y}(D^{L}u; z, \bar{z})
+\mathbb{Y}(D^{L}\mathbf{d}^{L} u; z,\bar{z}).
\end{eqnarray}
Applying (\ref{d-l-3}) to $\one$, 
taking the limit $z\to 0$ on both sides of (\ref{d-l-3}) and using the
creation property, 
(\ref{limit-z-der}) and $\mathbf{d}^{L}\one=0$, we obtain
$$\mathbf{d}^{L}D^{L}u=D^{L}u+D^{L}\mathbf{d}^{L} u,$$
proving that the left weight of $D^{L}$ is $1$. 
Similarly, we can prove that the right weight of $D^{L}$ is $0$,
the left weight of $D^{R}$ is $0$ and the right weight of $D^{R}$
is $1$.

Using the $D^{L}$- and $D^{R}$-derivative properties and the creation 
property, we see immediately that $D^{L}\one=D^{R}\one =0$. 
\epfv

For an $\R \times \R$-graded full field algebra
$(F, m, \mathbf{1}, D^{L}, D^{R})$, 
we now introduce a formal vertex operator map. 
We shall use the convention that for any $z\in \C^{\times}$,
$\log z=\log |z|+\sqrt{-1}\arg z$ where $0\le \arg z <2\pi$.
For $u\in F$, we use $\lwt u$ and $\rwt u$ to denote the 
left and right weights, respectively, of $u$. 
Let $u, v\in F$ and $w'\in F'$ be homogeneous elements. We have
\begin{eqnarray} \label{L-L-0-equ-1}
\langle w', [\mathbf{d}^L, \mathbb{Y}(u; z,\bar{z})]v\rangle &=& \langle 
(\mathbf{d}^L)' w', \mathbb{Y}(u; z, \bar{z})v \rangle-
\langle w', \mathbb{Y}(u; z, \bar{z}) \mathbf{d}^Lv\rangle  \nn
&=& (\lwt w'- \lwt v)\langle w', \mathbb{Y}(u; z, \bar{z})v\rangle 
\end{eqnarray}
where $(\mathbf{d}^L)'$ is the adjoint of $\mathbf{d}^L$. On the other hand, 
\begin{eqnarray} \label{L-L-0-equ-2}
\lefteqn{\langle w', \mathbb{Y}(\mathbf{d}^L u; z, \bar{z})v + 
z\frac{\partial}{\partial z} \mathbb{Y}(u; z, \bar{z})v \rangle}  \nn
&&=\left( \lwt u +  z\frac{\partial}{\partial z}\right) 
\langle w', \mathbb{Y}(u; z, \bar{z})v\rangle .
\end{eqnarray}
Let $f(z, \bar{z})=\langle w', \mathbb{Y}(u; z, \bar{z})v\rangle$.
Then
by (\ref{d-l}), (\ref{L-L-0-equ-1}) and (\ref{L-L-0-equ-2}), we have
\begin{equation}  \label{L-L-0-equ}
z\frac{\partial}{\partial z} f(z, \bar{z}) = (\lwt w'- \lwt u 
-\lwt v) f(z, \bar{z}).
\end{equation}
Similarly, using (\ref{d-r}), we have
\begin{equation}  \label{L-R-0-equ}
\bar{z}\frac{\partial}{\partial \bar{z}} f(z, \bar{z}) = 
(\rwt w'- \rwt u -\rwt v) f(z, \bar{z}).
\end{equation}
The general solution of the system (\ref{L-L-0-equ}) and (\ref{L-R-0-equ}) 
is
\begin{equation} \label{solution-1}
C z^{\slwt w'- \slwt u -\slwt v}\bar{z}^{\srwt w'- \srwt u -\srwt v}
\end{equation}
where $C\in \C$. Note that $f(z, \bar{z})$ is a single-valued
function and that by the single-valuedness of the full field algebra $F$,
$$(\lwt w'- \lwt u -\lwt v)- (\rwt w'- \rwt u -\rwt v)  \in \Z.$$
This means that if we choose any branches of 
$z^{\slwt w'- \slwt u -\slwt v}$ and
$\bar{z}^{\srwt w'- \srwt u -\srwt v}$, then there 
must be a unique constant $C$ such that $f(z, \bar{z})$ is 
equal to (\ref{solution-1}). We choose the branches of 
$z^{\slwt w'- \slwt u -\slwt v}$ and
$\bar{z}^{\srwt w'- \srwt u -\srwt v}$ to be 
$e^{(\slwt w'- \slwt u -\slwt v)\log z}$ and
$e^{(\srwt w'- \srwt u -\srwt v)\;\overline{\log z}}$, respectively. 
So there is a unique $C\in \C$ such that 
\begin{equation} \label{solution-1.5}
f(z, \bar{z})=C e^{(\slwt w'- \slwt u -\slwt v)\log z}
e^{(\srwt w'- \srwt u -\srwt v)\;\overline{\log z}}.
\end{equation}
Hence
$P_{p,q} \mathbb{Y}(u; z,\bar{z}) P_{m, n}$, $m, n, p, q \in \R$
can be written as 
$$
u_{m, n}^{p, q} e^{(p-\slwt u -m)\log z} e^{(q-\srwt u-n)\;\overline{\log
z}}, 
$$ 
where $u_{m ,n}^{p, q}$ are linear maps from 
$F_{(m, n)}\rightarrow F_{(p, q)}$ for $m, n, p, q \in \R$. 
Thus we have the following expansion: 
\begin{equation}  \label{expansion}
\mathbb{Y}(u; z,\bar{z}) =\sum_{r, s\in \R} \mathbb{Y}_{l, r}(u)
e^{(-l-1)\log z} e^{(-r-1)\;\overline{\log
z}}
\end{equation}
where $\mathbb{Y}_{l, r}(u)\in \edo F$ with
$\text{\rm wt}^L \mathbb{Y}_{l, r}(u) = \text{\rm wt}^L u -l-1$ and 
$\text{\rm wt}^R \mathbb{Y}_{l, r}(u) = \text{\rm wt}^R u -r-1$.
Moreover, the expansion above is unique. 
Let $x$ and $\bar{x}$ be independent and commuting formal 
variables. We define the {\it formal full vertex operator}
 $\mathbb{Y}_f$ associated 
to $u\in F$ by
\begin{equation} \label{formal-op-def}
\mathbb{Y}_f(u; x,\bar{x}) = 
\sum_{l, r\in \R}\mathbb{Y}_{l, r}(u) x^{-l-1}\bar{x}^{-r-1}.
\end{equation}
These formal full vertex operators give a {\it formal full vertex 
operator map}
$$\mathbb{Y}_{f}: F\otimes F\mapsto F\{x, \bar{x}\}.$$
For nonzero complex numbers $z$ and $\zeta$, 
we can substitute $e^{r\log z}$ and $e^{s\; \overline{\log \zeta}}$
for $x^{r}$ and $\bar{x}^{s}$, respectively, 
in $\mathbb{Y}_f(u; x,\bar{x})$ to obtain a map 
$\mathbb{Y}_{\rm an}(u; z, \zeta): F\to \overline{F}$ called the 
{\it analytic full vertex operator map}.

The following propositions are clear:

\begin{prop}
For $u\in F$ and $z, \zeta \in \C^{\times}$, we have
\begin{equation}\label{d-conj1}
\mathbb{Y}_{\rm an}(u; z, \zeta) 
= z^{\mathbf{d}^L} \zeta^{\mathbf{d}^R} \mathbb{Y}(
z^{-\mathbf{d}^L} \zeta^{-\mathbf{d}^R}u; 1, 1)
 z^{-\mathbf{d}^L} \zeta^{-\mathbf{d}^R}.
\end{equation}
For formal full vertex operators, we have 
\begin{equation}\label{d-conj2}
\mathbb{Y}_{f}(u; x, \bar{x}) 
= x^{\mathbf{d}^L} \bar{x}^{\mathbf{d}^R} 
\mathbb{Y}(x^{-\mathbf{d}^L} \bar{x}^{-\mathbf{d}^R}u; 1, 1)
 x^{-\mathbf{d}^L} \bar{x}^{-\mathbf{d}^R}.
\end{equation}
\end{prop}

\begin{prop}
For $u\in F$, 
\begin{eqnarray}
\mathbb{Y}_f(\one; x, \bar{x}) u &=& u, \label{form-ident} \\
\lim_{x\rightarrow 0, \bar{x}\rightarrow 0} \mathbb{Y}_f(u; x, \bar{x})
\one
&=& u \label{form-creat},
\end{eqnarray}
where $\lim_{x\rightarrow 0, \bar{x}\rightarrow 0}$ means taking 
the constant term of a power series in $x$ and $\bar{x}$.
In particular, $\mathbb{Y}_{l, r}(u)\one =0$ for all $l, r\in 
\mathbb{R}$ and 
$\mathbb{Y}_{-1,-1}(u)\one = u$. 
\end{prop}

\begin{prop}
For $u\in F$, we have
\begin{eqnarray}   
\left[ D^L, \mathbb{Y}_f(w; x, \bar{x})\right] &=& 
\mathbb{Y}_f(D^L w; x, \bar{x})
= \frac{\partial}{\partial x} \mathbb{Y}_f (w; x, \bar{x}),
\label{D-l-formal}\\
\left[ D^R, \mathbb{Y}_f (w; x, \bar{x})\right] &=& 
\mathbb{Y}_f (D^R w; x, \bar{x})
= \frac{\partial}{\partial \bar{x}} \mathbb{Y}_f(w, x, \bar{x}). 
\label{D-r-formal}
\end{eqnarray}
In particular, we have $D^L\one = D^R\one =0$ and 
for $l, r\in \R$,
\begin{eqnarray}  
\left[ D^L, \mathbb{Y}_{l, r}(u)\right] &=& \mathbb{Y}_{l, r}
(D^{L}u) = -l \mathbb{Y}_{l-1, r}(u),
\label{D-l-comp}\\
\left[ D^R, \mathbb{Y}_{l, r}(u)\right] &=& \mathbb{Y}_{l, r}
(D^{R}u) = -r \mathbb{Y}_{l, r-1}(u),
\label{D-r-comp}
\end{eqnarray}
\end{prop}

We need the following strong version of the creation property: 

\begin{lemma}
For $u\in F$, 
\begin{equation}   \label{strong-creat}
\mathbb{Y}_f(u; x, \bar{x}) \one = e^{xD^L + \bar{x}D^R} u.
\end{equation}
\end{lemma}
\pf
Using (\ref{D-l-formal}) and (\ref{D-r-formal}), we have
\begin{eqnarray}  \label{conj-D-l-r}
\mathbb{Y}_f(e^{x_0D^L + \bar{x}_0D^R}u; x, \bar{x}) 
&=&\mathbb{Y}_f(u; x+x_0, \bar{x}+\bar{x}_0). 
\end{eqnarray}
Now let both sides of (\ref{conj-D-l-r}) act on the vacuum $\one$. Since 
$\mathbb{Y}_f(u; x+x_0, \bar{x}+\bar{x}_0)\one$ involves
only nonnegative integer powers of $x+x_0$ and $\bar{x}+\bar{x}_0$, we can 
take the limit $x\to 0, \bar{x}\to 0$. 
Then replacing $x_0$ and $\bar{x}_0$ by $x$ and
$\bar{x}$, we obtain (\ref{strong-creat}).
\epfv

\begin{prop}[Skew symmetry]
For any $u,v\in F$ and $z\in \C^{\times}$, we have
\begin{equation}  \label{skew-symm}
\mathbb{Y}(u; z, \bar{z})v = e^{zD^L + \bar{z}D^R} 
\mathbb{Y}(v; -z, \overline{-z})u
\end{equation}
and 
\begin{equation}  \label{skew-formal}
\mathbb{Y}_f (u; x, \bar{x})v = e^{xD^L + \bar{x}D^R} 
\mathbb{Y}_f(v; e^{\pi i}x, e^{-\pi i}\bar{x})u.
\end{equation}
\end{prop}
\pf 
From the convergence property, it is clear that, for any $u,v\in F$, 
\begin{equation}
\mathbb{Y}(\mathbb{Y}(u; z_{1}-z_{2}, \bar{z}_{1}-\bar{z}_{2})v; z_{2},\bar{z}_{2}) \one
\end{equation}
converges absolutely to $m_{3}(u,v,\one;z_{1},\bar{z}_{1}, z_{2}, \bar{z}_{2},0,0)$ 
when $|z_{2}|>|z_{1}-z_{2}|>0$, 
and 
\begin{equation}
\mathbb{Y}(\mathbb{Y}(v; z_{2}-z_{1}, \bar{z}_{2}-\bar{z}_{1})u; z_{1},\bar{z}_{1}) \one
\end{equation}
converges absolutely when $|z_{1}|>|z_{1}-z_{2}|>0$
to $m_{3}(v,u,\one;z_{2},\bar{z}_{2}, z_{1}, \bar{z}_{1}, 0,0)$ which is 
equal to
$m_{3}(u,v,\one;z_{1},\bar{z}_{1}, z_{2}, \bar{z}_{2}, 0,0)$ by 
the permutation property.
Hence, using (\ref{strong-creat}), we obtain
\begin{equation}  \label{skew-proof-1}
e^{z_{2}D^L + \bar{z}_{2}D^R} \mathbb{Y}(u; z_{1}-z_{2}, 
\bar{z}_{1}-\bar{z}_{2})v
= e^{z_{1}D^L + \bar{z}_{1}D^R} \mathbb{Y}(v; z_{2}-z_{1}, 
\bar{z}_{2}-\bar{z}_{1})u
\end{equation}
when $|z_{2}|>|z_{1}-z_{2}|>0$ and $|z_{1}|>|z_{1}-z_{2}|>0$. 
We change the variables from $z_{1}, z_{2}$ to $z=z_{1}-z_{2}$ and 
$z_{2}$. Then (\ref{skew-proof-1}) gives
\begin{equation}  \label{skew-proof-1.5}
e^{z_{2}D^L + \bar{z}_{2}D^R} \mathbb{Y}(u; z, \bar{z})v
= e^{(z_{2}+z)D^L + (\bar{z}_{2}+\bar{z})D^R} 
\mathbb{Y}(v; -z, \overline{-z})u
\end{equation}
when $|z_{2}|>|z|>0$ and $|z_{2}+z|>|z|>0$.

Notice that for fixed $z\ne 0$ and $w'\in F'$, 
\begin{equation}  \label{skew-proof-2}
\langle w', e^{z_{2}D^L + \bar{z}_{2}D^R} \mathbb{Y}(u; z,
\bar{z})v\rangle 
\end{equation}
involves only positive integral powers of $z_{2}, \bar{z}_{2}$ and 
thus is a  power series in $z_{2}$ and $\bar{z}_{2}$
absolutely convergent when $|z_{2}|>|z|>0$. 
From complex analysis, we know that a power series in two variables
$z_{2}$ and $\zeta_{2}$
convergent at $z_{2}=z_{2}^{0}$ and $\zeta_{2}=\zeta_{2}^{0}$ 
must be convergent absolutely when $|z_{2}|<|z_{2}^{0}|$ and 
$|\zeta_{2}|<|\zeta_{2}^{0}|$. In particular, when $\zeta_{2}^{0}=
\bar{z}_{2}^{0}$, such a power series
must be absolutely convergent when $|z_{2}|<|z_{2}^{0}|$ and 
$\zeta_{2}=\bar{z}_{2}$. In our case, since for any fixed 
$z$, (\ref{skew-proof-2}) is absolutely convergent 
when $|z_{2}|>|z|>0$, we conclude that (\ref{skew-proof-2})
converges absolutely for all $z_{2}$. Since $w'$ is arbitrary, we see
that the left-hand side of (\ref{skew-proof-1.5}) is absolutely 
convergent in $\overline{F}$ for all $z_{2}$. Since $z$ is also arbitrary, 
by the convergence property again, we see that 
the left-hand side of (\ref{skew-proof-1.5}) is absolutely 
convergent in $\overline{F}$ for all $z$ and $z_{2}$ such that 
$z\ne 0$. 
Similarly, 
the right hand side of (\ref{skew-proof-1.5}) also converges absolutely 
in $\overline{F}$ for all $z$ and $z_{2}$ such that 
$z\ne 0$.

If $e^{-z_{2}D^L - \bar{z}_{2}D^R}$ gives a linear operator on $\overline{F}$,
then we can just multiply both sides of (\ref{skew-proof-1.5}) by 
$e^{-z_{2}D^L - \bar{z}_{2}D^R}$ to obtain (\ref{skew-symm}). In the case that 
the total weights of $F$ is lower-truncated, $e^{-z_{2}D^L - \bar{z}_{2}D^R}$
is indeed a linear operator on $\overline{F}$. In the most general 
case, this might not be true. But we can still obtain (\ref{skew-symm}) 
as follows: Consider the formal series 
\begin{eqnarray}\label{skew-proof-2.5}
\lefteqn{e^{x_{1}D^L + \bar{x}_{1}D^R} 
e^{x_{2}D^L + \bar{x}_{2}D^R} \mathbb{Y}(u; x, \bar{x})v}\nn
&&=e^{(x_{1}+x_{2})D^L + (\bar{x}_{1}+\bar{x}_{2})D^R} 
\mathbb{Y}(u; x, \bar{x})v
\end{eqnarray}
where $x$, $\bar{x}$, $x_{1}$, $\bar{x}_{1}$, $x_{2}$ and $\bar{x}_{2}$
are commuting formal variables. Since $\mathbb{Y}(u; z, \bar{z})v$
is absolutely convergent in $\overline{F}$ when $z\ne 0$, 
we can substitute $z$, $\bar{z}$, $-z_{2}$, $-\bar{z}_{2}$,
$z_{2}$ and $\bar{z}_{2}$ for $x$, $\bar{x}$, $x_{1}$, $\bar{x}_{1}$,
$x_{2}$ and $\bar{x}_{2}$ on the right-hand side 
of (\ref{skew-proof-2.5}), respectively, 
and the resulting series is absolutely convergent
in $\overline{F}$. So we can do the same substitution 
on the left-hand side of (\ref{skew-proof-2.5}) and the resulting 
series is absolutely convergent in $\overline{F}$. 
Similarly, consider the formal
series 
\begin{eqnarray}\label{skew-proof-3}
\lefteqn{e^{x_{1}D^L + \bar{x}_{1}D^R} 
e^{(x_{2}+x)D^L + (\bar{x}_{2}+\bar{x})D^R} 
\mathbb{Y}(v; e^{\pi i}x, e^{-\pi i}\bar{x})u}\nn
&&=e^{(x_{1}+x_{2}+x)D^L + (\bar{x}_{1}+\bar{x}_{2}+\bar{x})D^R} 
\mathbb{Y}(v; e^{\pi i}x, e^{-\pi i}\bar{x})u.
\end{eqnarray}
Since $e^{zD^L + \bar{z}D^R} 
\mathbb{Y}(v; -z, \overline{-z})u$
is absolutely convergent in $\overline{F}$ when $z\ne 0$, 
we can substitute $z$, $\bar{z}$, $-z_{2}$, $-\bar{z}_{2}$,
$z_{2}$ and $\bar{z}_{2}$ for $x$, $\bar{x}$, $x_{1}$, $\bar{x}_{1}$,
$x_{2}$ and $\bar{x}_{2}$ on
the right-hand side and thus also on the left-hand side of 
(\ref{skew-proof-3}) and the resulting series is 
absolutely convergent in $\overline{F}$. The convergence of these 
series and (\ref{skew-proof-1.5}) with suitably chosen $z_{2}$
gives (\ref{skew-symm})

Now  (\ref{skew-formal}) follows immediately:
On the one hand, by (\ref{d-conj2}), we have
\begin{eqnarray}\label{skew-formal-1}
x^{\mathbf{d}^{L}} \bar{x}^{\mathbf{d}^{L}} 
\mathbb{Y}(x^{-\mathbf{d}^{L}} \bar{x}^{-\mathbf{d}^{L}}u; 1, 1) 
x^{-\mathbf{d}^{L}} \bar{x}^{-\mathbf{d}^{L}}v = 
\mathbb{Y}_f(u; x, \bar{x})  v.
\end{eqnarray}
On the other hand, we have 
\begin{eqnarray}\label{skew-formal-2}
&x^{\mathbf{d}^{L}} \bar{x}^{\mathbf{d}^{L}}e^{D^L+D^R} 
\mathbb{Y}(x^{-\mathbf{d}^{L}} \bar{x}^{-\mathbf{d}^{L}}v;
-1, -1) x^{-\mathbf{d}^{L}} \bar{x}^{-\mathbf{d}^{L}}u &\nn
&=
e^{xD^L+\bar{x} D^R}\mathbb{Y}_f(v; e^{\pi i}x, e^{-\pi i}\bar{x}) 
u.&
\end{eqnarray}
Using skew symmetry (\ref{skew-symm}), (\ref{skew-formal-1})
and (\ref{skew-formal-2}),
we obtain (\ref{skew-formal}).
\epf

\begin{defn} \label{g-r-RR-g-ffa}
{\rm
An $\R \times \R$-graded full field algebra $(F, m, \one, 
D^{L}, D^{R})$ is called {\it grading
restricted} if it satisfies the following
grading-restriction conditions: 
\begin{enumerate}

\item There exists $M\in \R$ such that 
$F_{(m,n)}=0$ if $n<M$ or $m<M$.

\item $\dim F_{(m,n)}<\infty$ for $m, n\in \R$. \

\end{enumerate}
We  say that $F$ is {\it lower truncated} if $F$ satisfies
the first grading restriction condition. }
\end{defn} 

In this case, for $u\in F$ and $k\in \R$, we have  
$\sum_{l+r=k} \mathbb{Y}_{l, r}(u) \in \edo F$
with total weight $\wt u -k-2$. 
We denote $\sum_{l+r=k}  \mathbb{Y}_{l, r}(u) $ 
by $\mathbb{Y}_{k-1}(u)$. Then we have the  expansion
\begin{equation}  \label{Y-x-x}
\mathbb{Y}_f(u; x, x) = 
\sum_{k\in \R} \mathbb{Y}_k(u) x^{-k-1},
\end{equation}
where $\wt \mathbb{Y}_k(u) = \wt u - k -1$. For given $u, v\in F$, we have 
$\mathbb{Y}_k(u) w=0$ for sufficiently large $k$.

Let $(V^{L}, Y^{L}, \one^{L}, \omega^{L})$
and $(V^{R}, Y^{R}, \one^{R}, \omega^{R})$ be vertex operator
algebras. Let $\rho$ be an injective homomorphism from
the full field algebra $V^{L}\otimes V^{R}$ to $F$. 
Then we have $\one=\rho(\one^{L}\otimes \one^{R})$,
$\mathbf{d}^{L}\circ \rho=\rho\circ (L^{L}(0)\otimes I_{V^{R}})$,
$\mathbf{d}^{R}\circ \rho=\rho\circ (I_{V^{L}}\otimes  L^{R}(0)))$,
$D^{L}\circ \rho=\rho\circ (L^{L}(-1)\otimes I_{V^{R}})$ and
$D^{R}\circ \rho=\rho\circ (I_{V^{L}}\otimes  L^{R}(-1))$.
Moreover, $F$ has a {\it left conformal element} 
$\rho(\omega^{L}\otimes \one^{R})$ and an {\it right 
conformal element} $\rho(\one^{L}\otimes \omega^{R})$. 
We have the following operators on $F$: 
\begin{eqnarray*}
L^{L}(0)&=&\res_{x}\res_{\bar{x}}\bar{x}^{-1}\mathbb{Y}_{f}(
\rho(\omega^{L}\otimes
\mathbf{1}^{R}); x, \bar{x}),\\
L^{R}(0)&=&\res_{x}\res_{\bar{x}}x^{-1}\mathbb{Y}_{f}(
\rho(\mathbf{1}^{L}\otimes 
\omega^{L}); x, \bar{x}),\\
L^{L}(-1)&=&\res_{x}\res_{\bar{x}}x\bar{x}^{-1}\mathbb{Y}_{f}(
\rho(\omega^{L}\otimes
\mathbf{1}^{R}); x, \bar{x}),\\
L^{R}(-1)&=&\res_{x}\res_{\bar{x}}x^{-1}\bar{x}\mathbb{Y}_{f}(
\rho(\mathbf{1}^{L}\otimes
\omega^{L}); x, \bar{x}).
\end{eqnarray*}
Since these
operators are operators on $F$, it should be easy to distinguish 
them from those operators with the same noptation but 
acting on $V^{L}$ 
or $V^{R}$.

\begin{defn} 
{\rm Let $(V^{L}, Y^{L}, \one^{L}, \omega^{L})$ 
and $(V^{R}, Y^{R}, \one^{R}, \omega^{R})$ be vertex operator algebras. 
A {\it full field algebra over $V^{L}\otimes V^{R}$} 
is a grading-restricted $\R\times \R$-graded full field algebra $(F, m, \one,
D^{L}, D^{R})$ equipped with an injective homomorphism $\rho$ from 
the full field algebra $V^{L}\otimes V^{R}$ to $F$ such that 
$\mathbf{d}^{L}=L^{L}(0)$, $\mathbf{d}^{R}=L^{R}(0)$, 
$D^{L}=L^{L}(-1)$ and $D^{R}=L^{R}(-1)$.}
\end{defn}

We shall denote the full field algebra over $V^{L}\otimes V^{R}$
defined above by $(F, m, \rho)$ or simply by $F$.

The following result allows us to construct full field algebras
using the representation theory of vertex operator algebras:

\begin{thm}  \label{mod-int-op}
Let $(V^{L}, Y^{L}, \one^{L}, \omega^{L})$ 
and $(V^{R}, Y^{R}, \one^{R}, \omega^{R})$ be vertex operator algebras. 
Let  $(F, m, \rho)$ 
be a full field algebra over 
$V^{L}\otimes V^{R}$. Then 
$F$ is a module for 
$V^L\otimes V^R$ viewed as a vertex operator algebra. Moreover,  
$\mathbb{Y}_f(\cdot, x,x)$, 
is an intertwining operator of type $\binom{F}{FF}$. 
\end{thm}
\pf
Let $\mathbb{Y}^{L, R}$ be the vertex operator map for the 
full field  algebra $\rho(V^{L}\otimes V^{R})$. Then we have
\begin{equation}  \label{Y-split}
\mathbb{Y}^{L, R}(\rho(u^L\otimes u^R); z,\bar{z})\rho(v^L\otimes v^R)
=\rho(Y^L(u^L, z)v^L \otimes Y^R(u^R, \bar{z})v^R)
\end{equation}
for $u^L, v^L \in V^L$, $u^R, v^R\in V^R$ and $z\in \C^{\times}$.

Now we show that a splitting formula similar to (\ref{Y-split}) 
holds for vertex operators of the form 
$\mathbb{Y}(\rho(u^{L}\otimes u^{R}); z, \bar{z}): F\to \overline{F}$. 
By the associativity of $\mathbb{Y}$, we have 
\begin{eqnarray} \label{asso-K-K-alg-1}
\lefteqn{\langle w', \mathbb{Y}(\rho(u^L\otimes u^R); z_{1},\bar{z}_{1})
\mathbb{Y}(\rho(v^L\otimes v^R); z_{2},\bar{z}_{2})w\rangle} \nn
&& = \langle w', \mathbb{Y}(\mathbb{Y}^{L, R}
(\rho(u^L\otimes u^R); z_{1}-z_{2}, \bar{z}_{1}-\bar{z}_{2}) 
\rho(v^L\otimes v^R); z_{2},\bar{z}_{2})w\rangle\nn
&&
\end{eqnarray}
when $|z_{1}|>|z_{2}|>|z_{1}-z_{2}|>0$ for $u^L, v^{L}\in V^L$, 
$u^R, v^{R}\in V^R$,
$w\in F$ and $w'\in F'$.
Take $v^{L}=\one^{L}$ and $u^R, v^{R}=\one^R$. Then we have 
\begin{eqnarray} \label{asso-K-K-alg-2}
\lefteqn{\langle w', \mathbb{Y}(\rho(u^L\otimes \one^R); z_{1},\bar{z}_{1})
w\rangle} \nn
&&= \langle w', \mathbb{Y}(\rho(u^L\otimes \one^R); z_{1},\bar{z}_{1})
\mathbb{Y}(\one; z_{2},\bar{z}_{2})w\rangle\nn
&&= \langle w', \mathbb{Y}(\mathbb{Y}^{L, R}
(\rho(u^L\otimes \one^R); z_{1}-z_{2}, \bar{z}_{1}-\bar{z}_{2}) 
\rho(\one^L\otimes \one^R); z_{2},\bar{z}_{2})w\rangle\nn
&&=\langle w', \mathbb{Y}(\rho((Y^L
(u^L, z_{1}-z_{2})
\one^L)\otimes \one^R); z_{2},\bar{z}_{2})w\rangle
\end{eqnarray}
Since the right-hand side of (\ref{asso-K-K-alg-2}) 
is independent of $\bar{z}_{1}$, so is the left-hand side.
Thus we see that 
$\mathbb{Y}(\rho(u^L\otimes \one^R), z, \bar{z})$ depends only on $z$ for all
$u^L\in V^L$ and we shall also denote it by $Y^L(u^L, z)$. (Since 
it acts on $F$, there should be no confusion with the 
vertex operator $Y^L(u^L, z)$ acting on $V^{L}$.) So
$Y^L(u^L, z)$ is a series in powers of $z$. But $Y^L(u^L, z)$
is also single valued. So by (\ref{expansion}), there exists
$u_{n}^{L}\in \edo F$ for $n\in \Z$ such that
$\lwt u_n^L = \wt u^L-n-1$, $\rwt u_n^L =0$ and
$$Y^{L}(u^L, z)=\sum_{n\in \Z}u_{n}^{L}z^{-n-1}.$$
Similarly,  $\mathbb{Y}(\rho(\one^L\otimes u^R); z, \bar{z})$ depends only on 
$\bar{z}$ and  will also be denoted by $Y^R(u^R, \bar{z})$. 
(There should also be no confusion with the 
vertex operator $Y^R(u^R, z)$ acting on $V^{R}$.)
For $u^{R}\in V^{R}$, there exists $u_{n}^{R}\in \edo F$
for $n\in \Z$ such that $\rwt u_n^R = \wt u^R -n-1$, $\lwt u_n^R =0$ 
and
$$Y^{R}(u^R, z)=\sum_{n\in \Z}u_{n}^{R}z^{-n-1}.$$
We also have the formal vertex
operator maps, denoted using the same notations $Y^L$ and $Y^R$,
associated to $Y^L$ and $Y^R$ given by
\begin{eqnarray*}  
Y^L(u^L, x) &=& \sum_{n\in \Z} u_n^L x^{-n-1},  \\
Y^R(u^R,\bar{x})  &=& \sum_{n\in \Z} u_n^R \bar{x}^{-n-1} 
\end{eqnarray*}
for $u^L\in V^L$ and  $u^R\in V^R$.

For $u^L\in V^L, u^R\in V^R$ and $w\in F, w'\in F'$, 
\begin{equation} \label{L-R}
\langle w', Y^L(u^L,z_{1})Y^R(u^R, \bar{z}_{2})w\rangle
=\langle w', \mathbb{Y}(\rho(u^L\otimes \one^{R});
z_{1}, \bar{z}_{1})\mathbb{Y}(\rho(\one^{L}\otimes u^R); 
z_{2}, \bar{z}_{2})w\rangle
\end{equation} 
is absolutely convergent when $|z_{1}|>|z_{2}|>0$, and 
\begin{equation}  \label{R-L}
\langle w', Y^R(u^R, \bar{z}_{2})Y^L(u^L,z_{1})w\rangle
=\langle w', \mathbb{Y}(\rho(\one^{L}\otimes u^R); z_{2}, \bar{z}_{2})
\mathbb{Y}(\rho(u^L\otimes \one^{R});
z_{1}, \bar{z}_{1})w\rangle
\end{equation}
is absolutely convergent when $|z_{2}|>|z_{1}|>0$. They are 
both analytic in
$z_{1}$ and $\bar{z}_{2}$. 
By the convergence property for  full field algebras, 
both side of (\ref{L-R}) and 
(\ref{R-L}) can be extended to a same smooth function on 
$\{ (z_{1}, \bar{z}_{2}) \in (\C^{\times})^2 | z_{1}\neq z_{2} \}$. 
Since the complement of the union of the sets of convergence of
(\ref{L-R}) and (\ref{R-L}) in 
$\{ (z_{1}, \bar{z}_{2}) \in (\C^{\times})^2 | z_{1}\neq z_{2} \}$
 is of lower dimension, by the 
properties of analytic functions, 
it is clear that the extended smooth function is actually 
analytic on $\{ (z_{1}, \bar{z}_{2}) \in (\C^{\times})^2 | z_{1}\neq z_{2} \}$. 

By associativity, we have 
\begin{equation}   \label{L-R-assoc}
\langle w', Y^L(u^L,z_{1})Y^R(u^R, \bar{z}_{2})w\rangle
=\langle w', \mathbb{Y}(\rho((Y^L(u^L, z_{1}-z_{2})\one^L) \otimes u^R); z_{2}, 
\bar{z}_{2})w\rangle
\end{equation}
when $|z_{1}|>|z_{2}|>|z_{1}-z_{2}|>0$. The right-hand side of 
(\ref{L-R-assoc}) has
a well-defined limit as $z_{1}$ goes to $z_{2}$.
Therefore (\ref{L-R}) and (\ref{R-L}) can be further extended to a 
single analytic function on 
$\{ (z_{1},\bar{z}_{2}) \in (\C^{\times})^2 \}$. 
This absence of singularity further implies that 
the left-hand sides of (\ref{L-R}) and (\ref{R-L})
are absolutely convergent and are 
equal for all $z_{1}, \bar{z}_{2} \in \C^{\times}$. 
Let $z_{1}=z_{2}=z$ in (\ref{L-R}), (\ref{R-L}) and (\ref{L-R-assoc}). 
Use the discussion above and the creation property for the vertex 
operator map $Y^{L}$, we obtain 
\begin{equation}
\mathbb{Y}(\rho(u^L\otimes u^R); z, \bar{z}) = Y^L(u^L, z)Y^R(u^R, \bar{z}) =
Y^R(u^R, \bar{z})Y^L(u^L, z),
\end{equation}
or equivalently, in terms of formal vertex operator,
\begin{equation} \label{Y-cl-x-x}
\mathbb{Y}_f(\rho(u^L\otimes u^R); x, \bar{x})= Y^{L}(u^L, x)Y^R(u^R, \bar{x})
=Y^R(u^R, \bar{x})Y^L(u^L, x)
\end{equation}
for all $u^L\in V^L$ and $u^R\in V^R$.
In particular, we have 
$[u_m^L, u_n^R]=0$ for all $u^L \in V^L$ and $u^R\in V^R$. 

Since $F$ is lower truncated, we have
\begin{equation} \label{formal-Y-cl}
\mathbb{Y}_f(\rho(u^L\otimes u^R); x, x)v \in (\text{End}\; F)((x)).
\end{equation} 
for $u^L\in V^L$, $u^R\in V^R$ and $v\in F$. 

The associativity (\ref{asso-K-K-alg-1}) 
together with (\ref{Y-cl-x-x}) and (\ref{formal-Y-cl}) implies the 
associativity for the vertex operator map 
$\mathbb{Y}_f(\rho(\cdot); x, x)\cdot$. Together with the identity property 
this associativity implies that $F$ is a module for the vertex operator 
algebra $V^L\otimes V^R$.

Next we show that $\mathbb{Y}_f(\cdot; x, x)$ is an intertwining operator
of type $\binom{F}{FF}$. 

Since For given $u, v\in F$, we have 
$\mathbb{Y}_k(u) w=0$ for sufficient large $k$, 
the lower-truncation property of $\mathbb{Y}_f(\cdot, x, x)$
holds. For $u\in F$, 
We also have
\begin{eqnarray*}
\mathbb{Y}_f( (D^L + D^R) u; x, x)
&=&\mathbb{Y}_f( (D^L + D^R) u; x, \bar{x})|_{\bar{x}=x}\nn
& =& \left.\left(\left(\frac{\partial}{\partial x}+
\frac{\partial}{\partial\bar{x}}\right)
\mathbb{Y}_f(  u; x, \bar{x})\right)\right|_{\bar{x}=x}\nn
&=& \frac{d}{dx} \mathbb{Y}_f(u;x,x),
\end{eqnarray*}
proving the $D$-derivative property of  $\mathbb{Y}_f(\cdot; x, x)$.

Now, we prove the Jacobi identity
for $\mathbb{Y}_f(\cdot; x, x)$. 
For any fixed $r\in \R$, using the associativity
for the full vertex operator map 
$\mathbb{Y}$ twice, we obtain
\begin{eqnarray} \label{asso-equ}
\lefteqn{\langle w', Y^{L}(u^L, z_{1})Y^{R}(u^{R}, \bar{z}_{2}) 
\mathbb{Y}(u; r, r)w\rangle} \nn
&& = \langle w', \mathbb{Y}(\rho(u^L\otimes \one^{L});
z_{1}, \bar{z}_{1})\mathbb{Y}(\rho(\one^{L}\otimes u^{R}); z_{2}, \bar{z}_{2}) 
\mathbb{Y}(u; r, r)w\rangle\nn
&&=\langle w', 
\mathbb{Y}(\mathbb{Y}(\rho(u^L\otimes \one^{L});
z_{1}-r, \bar{z}_{1}-r)\mathbb{Y}(\rho(\one^{L}\otimes u^{R}); 
z_{2}-r, \bar{z}_{2}-r) 
u; r, r)w\rangle\nn
&&=\langle w', 
\mathbb{Y}(Y^{L}(u^L,
z_{1}-r)Y^{R}(u^{R}, \bar{z}_{2}-r) 
u; r, r)w\rangle\nn
\end{eqnarray}
when $|z_{1}|, |z_{2}|>r>|z_{1}-r|, |z_{2}-r|>0$
for all $u^L\in V^L$, $u^R\in V^R$, $u, w\in F$ and $w'\in F'$.
By the commutativity for the full vertex operator map
$\mathbb{Y}$, 
\begin{equation}  \label{L-R-Y-cl}
\langle w', Y^{L}(u^L, z_{1})Y^{R}(u^{R}, \bar{z}_{2}) 
\mathbb{Y}(u; r, r)w\rangle
\end{equation}
and 
\begin{equation}  \label{Y-cl-L-R}
\langle w', 
\mathbb{Y}(u; r, r)Y^{L}(u^L, z_{1})Y^{R}(u^{R}, \bar{z}_{2}) 
w\rangle
\end{equation}
are absolutely convergent in the regions 
$|z_{1}|, |z_{2}|>r>0$ and $r>|z_{1}|, |z_{2}|>0$, respectively,
to the correlation function
\begin{equation}\label{4-pt-fn}
\langle w', m_{4}(u^L\otimes \one^{L}, 
\one^{R}\otimes u^R, u, w; z_{1}, \bar{z}_{1}, z_{2}, \bar{z}_{2},
r, r, 0, 0).
\end{equation}

By our discussion above, we know that the right-hand side of (\ref{asso-equ}), 
(\ref{L-R-Y-cl}) and  (\ref{Y-cl-L-R}) are all analytic 
in $z_{1}$ and $\bar{z}_{2}$ and that we can take 
$z_{1}=\bar{z}_{2}$ in the right-hand side of (\ref{asso-equ}), 
(\ref{L-R-Y-cl}) and  (\ref{Y-cl-L-R}). Thus after taking 
$z_{1}=\bar{z}_{2}$, the right-hand side of (\ref{asso-equ}), 
(\ref{L-R-Y-cl}) and  (\ref{Y-cl-L-R}) are analytic in $z=z_{1}=\bar{z}_{2}$.
Since the right-hand side of (\ref{asso-equ}), 
(\ref{L-R-Y-cl}) and  (\ref{Y-cl-L-R}) are the expansions of 
(\ref{4-pt-fn}) in the regions $r>|z_{1}-r|, |z_{2}-r|>0$, 
$|z_{1}|, |z_{2}|>r>0$ and $r>|z_{1}|, |z_{2}|>0$, respectively,
we see that we can also let $z_{1}=\bar{z}_{2}$ in (\ref{4-pt-fn})
and the result is also analytic in $z=z_{1}=\bar{z}_{2}$. 
Thus we have proved that 
\begin{eqnarray*}
&\langle w', 
\mathbb{Y}(\mathbb{Y}_{f}(\rho(u^L\otimes u^{R}); z-r, z-r) 
u; r, r)w\rangle,&\\
&\langle w',\mathbb{Y}_{f}(\rho(u^L\otimes u^{R}); z, z) 
\mathbb{Y}(u; r, r)w\rangle,&\\
&\langle w', 
\mathbb{Y}(u; r, r)\mathbb{Y}_{f}(\rho(u^L\otimes u^{R}); z, z) 
w\rangle
\end{eqnarray*}
are absolutely convergent to 
$$\langle w', m_{4}(\rho(u^L\otimes \one^{L}), 
\rho(\one^{R}\otimes u^R), u, w; z, \bar{z}, \bar{z}, z,
r, r, 0, 0)$$
which is in fact analytic in $z$. Using the Cauchy formula for 
contour integrals, we obtain the Cauchy-Jacobi identity 
\begin{eqnarray}\label{cauchy-jacobi}
\lefteqn{\res_{z=\infty}f(z)\langle w',\mathbb{Y}_{f}(\rho(u^L\otimes u^{R});
z, z) 
\mathbb{Y}(u; r, r)w\rangle}\nn
&&-\res_{z=0}f(z)\langle w', 
\mathbb{Y}(u; r, r)\mathbb{Y}_{f}(\rho(u^L\otimes u^{R}); z, z) 
w\rangle\nn
&&=\res_{z=r}f(z)\langle w', 
\mathbb{Y}(\mathbb{Y}_{f}(\rho(u^L\otimes u^{R}); z-r, z-r) 
u; r, r)w\rangle,
\end{eqnarray}
where $f(z)$ is a rational function of $z$ with the only 
possible poles at $z=0, r, \infty$.
Since $w$ and $w'$ are arbitrary, 
this Cauchy-Jacobi identity gives us identities 
for the components of the vertex operator $\mathbb{Y}_{f}(u; x, x)$.
These identities are the component form of the Jacobi identity 
for $\mathbb{Y}_{f}(u; x, x)$. 
\epf

\begin{defn}
{\rm Let $c^{L}, c^{R}\in \C$. 
A {\it conformal full field algebra of central charges $(c^{L}, c^{R})$}
is a grading-restricted 
$\R\times \R$-graded full field algebra $(F, m, \one, D^{L}, D^{R})$
equipped with elements $\omega^{L}$ and $\omega^{R}$ called 
{\it left conformal element} and {\it right conformal element},
respectively, satisfying the following conditions:
\begin{enumerate}

\item The formal full vertex operators
$\mathbb{Y}_{f}(\omega^{L}; x, \bar{x})$ and 
$\mathbb{Y}_{f}(\omega^{R}; x, \bar{x})$ are Laurent series
in $x$ and $\bar{x}$, respectively, that is,
\begin{eqnarray*}
\mathbb{Y}_{f}(\omega^{L}; x, \bar{x})&=&\sum_{n\in \Z}L^{L}(n)x^{-n-2},\\
\mathbb{Y}_{f}(\omega^{R}; x, \bar{x})&=&\sum_{n\in \Z}L^{R}(n)\bar{x}^{-n-2}.\\
\end{eqnarray*}

\item The {\it Virasoro relations}: For $m, n\in \Z$, 
\begin{eqnarray*}
[L^{L}(m), L^{L}(n)]
&=&(m-n)L^{L}(m+n)+\frac{c^{L}}{12}(m^{3}-m)\delta_{m+n, 0},\\
{[L^{R}(m), L^{R}(n)]}
&=&(m-n)L^{R}(m+n)+\frac{c^{R}}{12}(m^{3}-m)\delta_{m+n, 0},\\
{[L^{L}(m), L^{R}(n)]}&=&0.
\end{eqnarray*}

\item $d^{L}=L^{L}(0)$, $d^{R}=L^{R}(0)$, $D^{L}=L^{L}(-1)$ and 
$D^{R}=L^{R}(-1)$.

\end{enumerate}}
\end{defn}

We shall denote the conformal full field algebra by 
$(F, m, \one, \omega^{L}, \omega^{R})$ or simply by $F$.

We have:

\begin{prop} 
Let $(F, m, \one, \omega^{L}, \omega^{R})$ be a conformal full field algebra. 
Then the following {\it commutator formula for Virasoro operators and 
formal full vertex operators} hold: For $u\in F$,
\begin{eqnarray}
\lefteqn{[\mathbb{Y}_{f}(\omega^{L}; x_{1}, \bar{x}_{1}), 
\mathbb{Y}_{f}(u; x_{2}, \bar{x}_{2})]}\nn
&&=\res_{x_{0}}x_{2}^{-1}\delta\left(\frac{x_{1}-x_{0}}{x_{2}}\right)
\mathbb{Y}_{f}(\mathbb{Y}_{f}(\omega^{L}; x_{0}, \bar{x}_{0})u; x_{2}, 
\bar{x}_{2}),\label{omega-commu-l}\\
\lefteqn{[\mathbb{Y}_{f}(\omega^{R}; x_{1}, \bar{x}_{1}), 
\mathbb{Y}_{f}(u; x_{2}, \bar{x}_{2})]}\nn
&&=\res_{\bar{x}_{0}}\bar{x}_{2}^{-1}
\delta\left(\frac{\bar{x}_{1}-\bar{x}_{0}}{\bar{x}_{2}}\right)
\mathbb{Y}_{f}(\mathbb{Y}_{f}(\omega^{R}; x_{0}, \bar{x}_{0})u; x_{2}, 
\bar{x}_{2}).\label{omega-commu-r}
\end{eqnarray}
\end{prop}
\pf
For any $v'\in F'$, $u, v\in F$, we consider 
\begin{equation}\label{m3-omega}
\langle v', m_{3}(\omega^{L}, u, v; z_{1}, \bar{z}_{1}, z_{2}, \bar{z}_{2},
0, 0)\rangle.
\end{equation}
Using the convergence property and the permutation 
property for conformal full field algebras,
we know that it is equal to 
\begin{eqnarray}
&\langle v', \mathbb{Y}(\omega^{L};  z_{1}, \bar{z}_{1})
\mathbb{Y}(u, z_{2}, \bar{z}_{2})v\rangle,&\label{omega-prod1}\\
&\langle v', \mathbb{Y}(u;  z_{2}, \bar{z}_{2})
\mathbb{Y}(\omega^{L}, z_{1}, \bar{z}_{1})v\rangle,&\label{omega-prod2}
\end{eqnarray}
in the regions $|z_{1}|>|z_{2}|>0$, $|z_{2}|>|z_{1}|>0$, respectively. 
By the definition of conformal full field algebra, 
we know that for any fixed $z_{2}\ne 0$,
(\ref{omega-prod1}) and 
(\ref{omega-prod2}) are analytic as functions of $z_{1}$ in the regions 
$|z_{1}|>|z_{2}|>0$ and $|z_{2}|>|z_{1}|>0$, respectively. 
So (\ref{m3-omega}) is analytic as a function 
of  $z_{1}$ in the regions $|z_{1}|>|z_{2}|>0$ and $|z_{2}|>|z_{1}|>0$.
But we know that (\ref{m3-omega}) is smooth as a function 
of  $z_{1}$ in $\C\setminus \{z_{2},  0\}$. 
Thus (\ref{m3-omega}) must be analytic in $\C\setminus \{z_{2},  0\}$. 

We know that (\ref{m3-omega})
is equal to (\ref{omega-prod1}),
(\ref{omega-prod2}) and 
$$\langle v', \mathbb{Y}(\mathbb{Y}(\omega^{L};  z_{1}-z_{2}, 
\bar{z}_{1}-\bar{z}_{2})u; z_{2}, \bar{z}_{2})v\rangle$$
in the regions $|z_{1}|>|z_{2}|>0$,  $|z_{2}|>|z_{1}|>0$
and $z_{2}|>|z_{1}-z_{2}|>0$, respectively. 
Since $F$ is lower truncated, using the Virasoro relation, we see  that 
$\mathbb{Y}_{f}(\omega^{L};  x, \bar{x})u$ and 
$\mathbb{Y}_{f}(\omega^{L};  x, \bar{x})v$ have only finitely many 
terms in negative powers of $x$. Also using the lower-truncation property
of $F$ and the Virasoro relation, we see that for any $w\in F$,
$\langle v', \mathbb{Y}_{f}(\omega^{L};  x, \bar{x})w\rangle$ 
has only finitely many terms in positive powers of $x$. 
Using these facts, we see that the singularities $z_{1}=z_{2}, 0, \infty$ of
(\ref{m3-omega}) are all poles. Using 
the Cauchy formula, we obtain the component form (\ref{omega-commu-l}). 

Similarly, we can prove (\ref{omega-commu-r}). 
\epfv

The following is clear from the definition and Theorem \ref{mod-int-op}:

\begin{prop}
Let $(V^{L}, Y^{L}, \one^{L}, \omega^{L})$ and $(V^{L}, Y^{L}, \one^{L}, 
\omega^{L})$ 
be vertex operator algebras of central 
charges $c^{L}$ and $c^{R}$, respectively.
A full field algebra $(F, m, 
\rho)$ over $V^{L}\otimes V^{R}$
equipped with the left and right conformal elements
$\rho(\omega^{L}\otimes \one^{R})$ and $\rho(\one^{L}\otimes \omega^{R})$
is a conformal full field algebra.
\end{prop}

In view of this proposition, we shall call 
the conformal full field algebra in the proposition above, that is, 
a full field algebra $(F, m, 
\rho)$ over $V^{L}\otimes V^{R}$ equipped with the left and right 
conformal elements $\rho(\omega^{L}\otimes \one^{R})$ and 
$\rho(\one^{L}\otimes \omega^{R})$,  a
{\it conformal full field algebra over $V^{L}\otimes V^{R}$}
and denote it by $(F, m, \rho)$
or simply by $F$.

\renewcommand{\theequation}{\thesection.\arabic{equation}}
\renewcommand{\thethm}{\thesection.\arabic{thm}}
\setcounter{equation}{0}
\setcounter{thm}{0}

\section{Intertwining operator algebras and full field algebras}

Let $V^{L}$ and $V^{R}$ be vertex operator algebras. 
In the preceding section, we have shown that a conformal
full field algebra $(F, m, \rho)$
over $V^{L}\otimes V^{R}$ is a module for the vertex operator 
algebra $V^{L}\otimes V^{R}$ and the $\mathbb{Y}_{f}(\cdot; x, x)$ 
is an intertwining operator
of type $\binom{F}{FF}$. This result suggests a method to construct
conformal 
full field algebras from intertwining operator algebras, 
which are algebras of
intertwining operators for vertex operator algebras and 
were introduced and studied in \cite{H1}, \cite{H2},
\cite{H3}, \cite{H3.5}, \cite{H4} and \cite{H5} by the first author. 

Let $V$ be a vertex operator algebra and for a  $V$-module $W$,
let $C_{1}(W)$ be the subspace of 
$V$ spanned by $u_{-1}w$ for $u\in V_{+}=\coprod_{n\in \Z_{+}}V_{(n)}$ 
and $w\in W$. 

We consider the following conditions for a vertex operator algebra $V$:

\begin{enumerate}

\item  Every $\C$-graded generalized $V$-module is a direct sum of
$\C$-graded irreducible
$V$-modules.

\item There are only finitely many inequivalent $\C$-graded
irreducible $V$-modules and they are all $\R$-graded.

\item Every
$\R$-graded irreducible $V$-module $W$ satisfies
the $C_{1}$-cofiniteness condition, that is, $\dim W/C_{1}(W)<\infty$.

\end{enumerate}

In this section, we fix vertex operator algebras 
$(V^L, Y^{L}, \one^{L}, \omega^{L})$ 
and $(V^{R}, Y^{R}, \one^{R}, \omega^{R})$  satisfying 
these conditions. 
Let $\mathcal{A}^L$ and $\mathcal{A}^R$ be the sets of 
equivalent classes of irreducible modules for $V^L$ and for $V^R$, 
respectively. Let $\{W^{L; a}\;|\;a\in \A^{L}\}$ be a complete set of 
representatives of the equivalence classes in 
$\A^{L}$ and $\{W^{R; b}\;|\;b\in \A^{R}\}$
a complete set of 
representatives of the equivalence classes in $\A^{R}$.

\begin{prop} \label{rational}
The vertex operator algebra $V^L\otimes V^R$ also satisfy the 
conditions above.
\end{prop}
\pf
Let $W$ be a generalized $V^L\otimes V^R$-module. 
Then $W$ is a generalized $V^{L}$-module. So there exist
vector spaces $M^{a}$ for $a\in \A^{L}$ such that 
$W$ is equivalent to the generalized $V^{L}$-module 
$\coprod_{a\in \A^{L}}(W^{L;a}\otimes M^{a})$. Since $W$ 
is also a generalized $V^{R}$-module, $M^{a}$ must be $V^{R}$-modules. 
So they can be written as direct sums of irreducible $V^{R}$-modules 
$W^{R;b}$, $b\in \A^{R}$. So $W$ is equivalent to 
$\coprod_{a\in \A^{L}, b\in \A^{R}}N_{ab}(W^{L;b}\otimes W^{R;a})$
where $N_{ab}\in \N$ for $a\in \A^{L}$, $b\in \A^{R}$.
By Proposition 4.7.2 of \cite{FHL}, $W^{L;a}\otimes W^{R;b}$ are irreducible 
$V^L\otimes V^R$-modules. So $V^L\otimes V^R$ satisfies Condition 1. 
The second condition follows from Theorem 4.7.4 of \cite{FHL}.
The $C_{1}$-cofiniteness follows immediately from 
the fact that 
$$C_{1}(W^{L; a})\otimes W^{R; b}\oplus W^{L; a}\otimes 
C_{1}(W^{R; b})\subset
C_{1}(W^{L; a}\otimes W^{R; b}).$$
\epfv

This result immediately gives:

\begin{cor}  \label{mod-decomp}
Let $F$ be a module for the vertex operator algebra $V^{L}\otimes V^{R}$.
Then as a module for the vertex operator algebra $V^{L}\otimes V^{R}$, 
$F$ is isomorphic to 
\begin{equation}  \label{F-exp}
\coprod_{a\in \mathcal{A}^{L}}\coprod_{b\in \mathcal{A}^{R}} \coprod_{m_{ab}=1}
^{h_{ab}}(W^{L; a})^{(m_{ab})} \otimes  (W^{R; b})^{(m_{ab})}
\end{equation}
\end{cor}

Let $F$ be a module for the vertex operator algebra $V^{L}\otimes V^{R}$
and let $\gamma$ be an isomorphism from (\ref{F-exp}) to $F$. 
Then there exist operators 
$L^{L}(0)$ and $L^{R}(0)$ on $F$ given by 
\begin{eqnarray*}
L^{L}(0)\rho(w^{L}\otimes w^{R})&=&\rho((L^{L}(0)w^{L})\otimes w^{R}),\\
L^{R}(0)\rho(w^{L}\otimes w^{R})&=&\rho(w^{L}\otimes (L^{R}(0)w^{R}))
\end{eqnarray*}
for $w^{L}\in (W^{L; a})^{(m_{ab})}$ and $w^{R}\in (W^{R; b})^{(m_{ab})}$.
Clearly $L^{L}(0)$ and $L^{R}(0)$ commute with each other. 

Let $\Y$ be an intertwining operator of type $\binom{F}{FF}$ and 
$\gamma$ an isomorphism from (\ref{F-exp}) to $F$. 
Let 
\begin{eqnarray*}
\mathbb{Y}^{\mathcal{Y}}: (F\otimes F)\times \C^{\times}&\to&
\overline{F}\nn
(u\otimes v, z)&\mapsto& \mathbb{Y}^{\mathcal{Y}}(u; z, \bar{z})v
\end{eqnarray*}
and 
\begin{eqnarray*}
\mathbb{Y}_{f}^{\mathcal{Y}}: F\otimes F&\to&
F\{x, \bar{x}\}\nn
u\otimes v&\mapsto& \mathbb{Y}_{f}^{\mathcal{Y}}(u; x, \bar{x})v
\end{eqnarray*}
be linear maps given by 
$$\mathbb{Y}^{\mathcal{Y}}(u; z, \bar{z})v=z^{L^L(0)}
\bar{z}^{L^R(0))} \mathcal{Y}(u,  1)
z^{-L^L(0)}  \bar{z}^{-L^R(0)}$$
and 
$$\mathbb{Y}_{f}^{\mathcal{Y}}(u; x, \bar{x})v=x^{L^L(0)}
\bar{x}^{L^R(0)} \mathcal{Y}(u, 1)
x^{-L^L(0)}  \bar{x}^{-L^R(0)},$$
respectively, for $u\in F$. We call $\mathbb{Y}^{\mathcal{Y}}$ 
and $\mathbb{Y}_{f}^{\mathcal{Y}}$
the {\it splitting} and {\it formal splitting} of $\Y$, respectively. 

\begin{prop}\label{int-decomp}
Let $\Y$ be an intertwining operator of type $\binom{F}{FF}$,
$\mathbb{Y}^{\mathcal{Y}}$ and $\mathbb{Y}_{f}^{\mathcal{Y}}$,
the splitting and formal splitting of $\Y$, respectively,
and $\gamma$ an isomorphism from (\ref{F-exp}) to $F$. 
Then for any $a_{1}, a_{2}\in \A^{L}$, $b_{1}, b_{2}\in 
\A^{R}$, $1\le m_{a_{1}b_{1}}\le h_{a_{1}b_{1}}$ and $1\le m_{a_{2}b_{2}}\le 
h_{a_{2}b_{2}}$, there exist intertwining operators 
$\Y_{a_{1}a_{2}}^{L;m_{a_{3}b_{3}};a_{3}}$ and 
$\Y_{b_{1}b_{2}}^{R; m_{a_{3}b_{3}}; b_{3}}$
for $a_{3}\in \A^{L}$, $b_{3}\in \A^{R}$ and $m_{a_{3}b_{3}}=1, \dots, 
h_{a_{3}b_{3}}$ of types 
$\binom{(W^{L;a_{3}})^{(m_{a_{3}b_{3}})}}
{(W^{L;a_{1}})^{(m_{a_{1}b_{1}})}\;(W^{L;a_{2}})^{(m_{a_{2}b_{2}})}}$ and 
$\binom{(W^{R;b_{3}})^{(m_{a_{3}b_{3}})}}
{(W^{R;b_{1}})^{(m_{a_{1}b_{1}})}\;(W^{R;b_{2}})^{(m_{a_{2}b_{2}})}}$, 
respectively, such that for
$u^L\otimes u^R \in (W^{L; a_{1}})^{(m_{a_{1}b_{1}})}\otimes 
(W^{R; b_{1}})^{(m_{a_{1}b_{1}})}$ and
$v^L\otimes v^R \in (W^{L; a_{2}})^{(m_{a_{2}b_{2}})}\otimes 
(W^{R; b_{2}})^{(m_{a_{2}b_{2}})}$, we have 
\begin{eqnarray}  \label{ch-ach}
\lefteqn{\mathbb{Y}^{\mathcal{Y}}(\gamma(u^L\otimes u^R); z,\bar{z}) 
\gamma(v^L \otimes v^R)}\nn
&& =
\sum_{a_{3}\in \A^{L}}\sum_{b_{3}\in \A^{R}}
\sum_{m_{a_{3}b_{3}}=1}^{h_{a_{3}b_{3}}}\gamma(
\Y_{a_{1}a_{2}}^{L;m_{a_{3}b_{3}};a_{3}}(u^L, z)v^L \otimes 
\Y_{b_{1}b_{2}}^{R; m_{a_{3}b_{3}}; b_{3}}(u^R,\bar{z})v^R).\nn
\end{eqnarray}
Similarly, for the formal full vertex operator, we have 
\begin{eqnarray}  \label{ch-ach-f}
\lefteqn{\mathbb{Y}_{f}^{\mathcal{Y}}(\gamma(u^L\otimes u^R); x,\bar{x}) 
\gamma(v^L \otimes v^R)}\nn
&& =
\sum_{a_{3}\in \A^{L}}\sum_{b_{3}\in \A^{R}}
\sum_{m_{a_{3}b_{3}}=1}^{h_{a_{3}b_{3}}}\gamma(
\Y_{a_{1}a_{2}}^{L;m_{a_{3}b_{3}};a_{3}}(u^L, x)v^L \otimes 
\Y_{b_{1}b_{2}}^{R; m_{a_{3}b_{3}}; b_{3}}(u^R,\bar{x})v^R).\nn
\end{eqnarray}
\end{prop}
\pf 
Since $\mathbb{Y}^{\Y}$ restricted to 
$$\gamma((W^{L; a_{1}})^{(m_{a_{1}b_{1}})}\otimes 
(W^{R; b_{1}})^{(m_{a_{1}b_{1}})})\otimes 
\gamma((W^{L; a_{2}})^{(m_{a_{2}b_{2}})}\otimes 
(W^{R; b_{2}})^{(m_{a_{2}b_{2}})})$$
is an intertwining operator of type 
$$\binom{F}{\gamma((W^{L; a_{1}})^{(m_{a_{1}b_{1}})}\otimes 
(W^{R; b_{1}})^{(m_{a_{1}b_{1}})})\;\;
\gamma((W^{L; a_{2}})^{(m_{a_{2}b_{2}})}\otimes 
(W^{R; b_{2}})^{(m_{a_{2}b_{2}})})},$$
it was proved in \cite{DMZ} that 
(\ref{ch-ach}) is true when $z=\bar{z}=r>0$. Then we have
\begin{eqnarray*}  \label{Y-cl-z}
\lefteqn{\mathbb{Y}^{\Y}(\gamma(u^L\otimes u^R); z, \bar{z})
\gamma(v^L \otimes v^R)}\nn
&&= z^{L^L(0)}
\bar{z}^{L^R(0)} \mathcal{Y}(\gamma(u^L\otimes u^R), 1)
z^{-L^L(0)}  \bar{z}^{-L^R(0)} \gamma(v^L \otimes v^R) \nn
&&= \sum_{a_{3}\in \A^{L}}\sum_{b_{3}\in \A^{R}}
\sum_{m_{a_{3}b_{3}}=1}^{h_{a_{3}b_{3}}}\gamma(
(z^{L^L(0)}\Y_{a_{1}a_{2}}^{L;m_{a_{3}b_{3}};a_{3}}(u^L, 1)
\bar{z}^{-L^L(0)}v^{L})\nn
&&\hspace{12em}\otimes (z^{L^R(0)} \Y_{b_{1}b_{2}}^{R; m_{a_{3}b_{3}}; b_{3}}
(u^R, 1) \bar{z}^{-L^R(0)}v^{R}))  \nn
&&= \sum_{a_{3}\in \A^{L}}\sum_{b_{3}\in \A^{R}}
\sum_{m_{a_{3}b_{3}}=1}^{h_{a_{3}b_{3}}}\gamma(
\Y_{a_{1}a_{2}}^{L;m_{a_{3}b_{3}};a_{3}}
(u^{L}, z)v^{L}\otimes \Y_{b_{1}b_{2}}^{R; m_{a_{3}b_{3}}; b_{3}}
(u^{R}, \bar{z})v^{R}). 
\end{eqnarray*}
The proof of (\ref{ch-ach-f}) is completely the same.
\epfv

\begin{cor}\label{r-cfa-decomp} 
Let 
$(F, m, \rho)$
be a conformal full field algebra over 
$V^L\otimes V^R$.
Then as a module for the vertex operator algebra $V^{L}\otimes V^{R}$, 
$F$ is isomorphic to (\ref{F-exp}).
Moreover, if $\gamma$ is an isomorphism from (\ref{F-exp})
to $F$, then  for any $a_{1}, a_{2}\in \A^{L}$, $b_{1}, b_{2}\in 
\A^{R}$, $1\le m_{a_{1}b_{1}}\le h_{a_{1}b_{1}}$ and $1\le m_{a_{2}b_{2}}\le 
h_{a_{2}b_{2}}$, there exist intertwining operators 
$\Y_{a_{1}a_{2}}^{L;m_{a_{3}b_{3}};a_{3}}$ and 
$\Y_{b_{1}b_{2}}^{R; m_{a_{3}b_{3}}; b_{3}}$
for $a_{3}\in \A^{L}$, $b_{3}\in \A^{R}$ and $m_{a_{3}b_{3}}=1, \dots, 
h_{a_{3}b_{3}}$ of types 
$\binom{(W^{L;a_{3}})^{(m_{a_{3}b_{3}})}}
{(W^{L;a_{1}})^{(m_{a_{1}b_{1}})}\;(W^{L;a_{2}})^{(m_{a_{2}b_{2}})}}$ and 
$\binom{(W^{R;b_{3}})^{(m_{a_{3}b_{3}})}}
{(W^{R;b_{1}})^{(m_{a_{1}b_{1}})}\;(W^{R;b_{2}})^{(m_{a_{2}b_{2}})}}$, 
respectively, such that for
$u^L\otimes u^R \in (W^{L; a_{1}})^{(m_{a_{1}b_{1}})}\otimes 
(W^{R; b_{1}})^{(m_{a_{1}b_{1}})}$ and
$v^L\otimes v^R \in (W^{L; a_{2}})^{(m_{a_{2}b_{2}})}\otimes 
(W^{R; b_{2}})^{(m_{a_{2}b_{2}})}$,  the formulas (\ref{ch-ach}) and
(\ref{ch-ach-f}) hold when
$\mathbb{Y}^{\mathcal{Y}}$ and $\mathbb{Y}^{\mathcal{Y}}_{f}$ are
replaced by $\mathbb{Y}$  and 
$\mathbb{Y}_{f}$, respectively.
\end{cor}
\pf
The first conclusion follows immediately from Corollary \ref{mod-decomp}.
Now if we consider the intertwining operator 
$\mathbb{Y}_{f}(\cdot; x, x)$, then the second conclusion 
follows immediately from Proposition \ref{int-decomp}. 
\epfv

For either the map $\mathbb{Y}_{f}^{\Y}$ in Proposition \ref{int-decomp}
or the formal full vertex operator map $\mathbb{Y}_{f}$ for a 
conformal full field algebra over 
$V^L\otimes V^R$, we can substitute $z$ and $\zeta$ for the formal variables 
$x$ and $\bar{x}$ in $\mathbb{Y}_{f}^{\Y}(\cdot; x,\bar{x})$
or $\mathbb{Y}_{f}(\cdot; x,\bar{x})$ (that is, 
substitute $e^{r\log z}$ and $e^{s\;\overline{\log \zeta}}$ for 
$x^{r}$ and $\bar{x}^{s}$, respectively, for $r, s\in \R$) to obtain 
$\mathbb{Y}_{\rm an}^{\Y}(\cdot; z, \zeta)$ (called {\it 
analytic splitting of $\Y$}) or
$\mathbb{Y}_{\rm an}(\cdot; z, \zeta)$. Then by (\ref{ch-ach-f}), 
we have:

\begin{cor}\label{fvo-decomp}
For the analytic splitting $\mathbb{Y}_{\rm an}^{\Y}$ of $\Y$
in Proposition \ref{int-decomp}, we have
\begin{eqnarray}\label{ch-ach-z-zeta}
\lefteqn{\mathbb{Y}_{\rm an}^{\Y}(\gamma(u^L\otimes u^R); z,\zeta) 
\gamma(v^L \otimes v^R)}\nn
&& =
\sum_{a_{3}\in \A^{L}}\sum_{b_{3}\in \A^{R}}
\sum_{m_{a_{3}b_{3}}=1}^{h_{a_{3}b_{3}}}\gamma(
\Y_{a_{1}a_{2}}^{L;m_{a_{3}b_{3}};a_{3}}(u^L, z)v^L \otimes 
\Y_{b_{1}b_{2}}^{R; m_{a_{3}b_{3}}; b_{3}}(u^R,\zeta)v^R)\nn
\end{eqnarray}
for
$u^L\otimes u^R \in (W^{L; a_{1}})^{(m_{a_{1}b_{1}})}\otimes 
(W^{R; b_{1}})^{(m_{a_{1}b_{1}})}$ and
$v^L\otimes v^R \in (W^{L; a_{2}})^{(m_{a_{2}b_{2}})}\otimes 
(W^{R; b_{2}})^{(m_{a_{2}b_{2}})}$. The same is also true
for the analytic full vertex operator map $\mathbb{Y}_{\rm an}$ for a 
conformal full field algebra over 
$V^L\otimes V^R$.
\end{cor}

This corollary allows us to treat the left and right variables $z$ 
and $\bar{z}$
in $\mathbb{Y}^{\Y}(\cdot; z, \bar{z})$ or $\mathbb{Y}(\cdot; z, \bar{z})$ 
independently.
In particular, we have the following 
strong versions of associativity and commutativity
for conformal full field algebra over 
$V^L\otimes V^R$:

\begin{prop}[Associativity]\label{anal-assoc}
Let $(F, m, \rho)$ 
be a conformal full field algebra over 
$V^L\otimes V^R$.
Then for
$u,v,w\in F$ and $w'\in F'$, 
\begin{eqnarray}  \label{asso-z-zeta}
\lefteqn{\langle w', \mathbb{Y}_{\rm an}(u; z_{1},\zeta_{1})
\mathbb{Y}_{\rm an}(v; z_{2}, \zeta_{2})w\rangle} \nn
&&= \langle w', 
\mathbb{Y}_{\rm an}(\mathbb{Y}_{\rm an}(u; z_{1}-z_{2}, \zeta_{1}-\zeta_{2})v;
 z_{2}, \zeta_{2})w\rangle
\end{eqnarray}
when $|z_{1}|>|z_{2}|>|z_{1}-z_{2}|>0$ and $|\zeta_{1}|>|\zeta_{2}|>|\zeta_{1}-\zeta_{2}|>0$. 
\end{prop}
\pf
Using (\ref{ch-ach-z-zeta}) and the convergence result 
proved by the first author in \cite{H5} for vertex operator algebras
satisfying the conditions assumed for $V^L$ and $V^R$ in the beginning of 
this section,  
the left-hand side of (\ref{asso-z-zeta}) converges absolutely when 
$|z_{1}|>|z_{2}|>0$ and $|\zeta_{1}|>|\zeta_{1}|>0$, and 
the right-hand side of (\ref{asso-z-zeta}) converges absolutely when
$|z_{2}|>|z_{1}-z_{2}|>0$ and $|\zeta_{2}|>|\zeta_{1}-\zeta_{2}|>0$. 
By the associativity (\ref{asso-1}), (\ref{asso-z-zeta}) 
is true when $\zeta_{1}=\bar{z}_{1}$
and $\zeta_{2}=\bar{z}_{2}$ for all $u,v,w\in F$ and $w'\in F'$. 
In particular, replacing
$u$ by $(L^L(-1))^k (L^R(-1))^l u$, $v$ by $(L^L(-1))^m L^R(-1))^n v$,
for $k,l, m, n\in \N$ and using the $L^{L}(-1)$- and $L^{R}(-1)$-derivative
properties,  we obtain 
\begin{eqnarray}  \label{asso-z-zeta-1}
\lefteqn{\frac{\partial^k}{\partial z_{1}^k} 
\frac{\partial^l}{\partial \zeta_{1}^l}
\frac{\partial^m}{\partial z_{2}^m} \frac{\partial^n}{\partial \zeta_{2}^n}
\langle w', \mathbb{Y}_{\rm an}(u; z_{1},\zeta_{1})
\mathbb{Y}_{\rm an}(v; z_{2}, \zeta_{2})w\rangle
\lbar_{\zeta_{1}=\bar{z}_{1}, \zeta_{2}=\bar{z}_{2}}
} \nn
&&= \frac{\partial^k}{\partial z_{1}^k} \frac{\partial^l}{\partial \zeta_{1}^l}
\frac{\partial^m}{\partial z_{2}^m} \frac{\partial^n}{\partial \zeta_{2}^n}
\langle w', 
\mathbb{Y}_{\rm an}(\mathbb{Y}_{\rm an}(u; z_{1}-z_{2}, \zeta_{1}-\zeta_{2})v;
 z_{2}, \zeta_{2})w\rangle\lbar_{\zeta_{1}=\bar{z}_{1}, \zeta_{2}=\bar{z}_{2}}\nn
\end{eqnarray}
for all $k,l, m, n\in \N$, when 
$|z_{1}|>|z_{2}|>0$ and
$|z_{2}|>|z_{1}-z_{2}|>0$. We know that both sides of (\ref{asso-z-zeta})
give branches of some multivalued
analytic functions in the region given by 
$|z_{1}|>|z_{2}|>0$, $|\zeta_{1}|>|\zeta_{1}|>0$,
$|z_{2}|>|z_{1}-z_{2}|>0$ and $|\zeta_{2}|>|\zeta_{1}-\zeta_{2}|>0$.
From (\ref{asso-z-zeta-1}), we know that the power series expansions of 
these branches
are equal in the neighborhood of those points satisfying 
$\zeta_{1}=\bar{z}_{1}$, $\zeta_{2}=\bar{z}_{2}$.  
Thus (\ref{asso-z-zeta-1}) holds in the region
$|z_{1}|>|z_{2}|>0$, $|\zeta_{1}|>|\zeta_{1}|>0$,
$|z_{2}|>|z_{1}-z_{2}|>0$ and $|\zeta_{2}|>|\zeta_{1}-\zeta_{2}|>0$. 
\epfv

\begin{prop}[Commutativity] \label{comm-0}
Let $(F, m, \rho)$ 
be a conformal full field algebra over 
$V^L\otimes V^R$. Then for $u,v,w\in F$ and $w'\in F'$, 
\begin{equation} \label{prod-12}
\langle w', \mathbb{Y}_{\rm an}(u; z_{1}, \zeta_{1})
\mathbb{Y}_{\rm an}(v; z_{2}, \zeta_{2} )w\rangle 
\end{equation}
and
\begin{equation} \label{prod-21}
\langle w', \mathbb{Y}_{\rm an}(v; z_{2}, \zeta_{2})
\mathbb{Y}_{\rm an}(u; z_{1}, \zeta_{1} )w\rangle
\end{equation}
are absolutely convergent when $|z_{1}|>|z_{2}|>0$, 
$|\zeta_{1}|>|\zeta_{2}|>0$ and when $|z_{2}|>|z_{1}|>0$, 
$|\zeta_{2}|>|\zeta_{1}|>0$, respectively, and 
can both be analytically extended to a same multivalued analytic function
of $(z_{1}, z_{2}; \zeta_{1}, \zeta_{2})$ for 
$(z_{1}, z_{2}; \zeta_{1}, \zeta_{2})\in M^{2}\times M^{2}$,
where $M^{2}=\{(z_{1}, z_{2})\in (\C^{\times})^{2}\;|\;z_{1}\ne z_{2}\}$. 
\end{prop}
\pf
The convergence and the existence of analytic extensions follow
immediately from Corollary \ref{fvo-decomp} and the 
convergence and the existence of analytic extensions of 
products of intertwining operators for the vertex operator algebras 
$V^{L}$ and $V^{R}$. 

By Proposition \ref{commu-general}, we know that 
these two multivalued 
analytic functions obtained by analytically extending 
(\ref{prod-12}) and (\ref{prod-21}) have equal values
at points of the form 
$(z_{1}, \bar{z}_{1}, 
\dots, z_{n}, \bar{z}_{n})$ for $(z_{1}, \dots, z_{n})\in 
\mathbb{F}_{n}(\C)$. 
Using the $L^{L}(-1)$- and $L^{R}(-1)$-conjugation properties for 
full vertex operators, we see that these two analytic functions 
are actually the same, that is, they are analytic extensions of each 
other. 
\epfv

For $(z_{1}, \dots, z_{n}), (\zeta_{1}, \dots, \zeta_{n})\in 
\mathbb{F}_{n}(\C)$, we denote the corresponding elements of 
$\mathbb{F}_{n}(\C)\times \mathbb{F}_{n}(\C)$ by 
$(z_{1}, \zeta_{1}, \dots, z_{n}, \zeta_{n})$ instead of 
$(z_{1}, \dots, z_{n}, \zeta_{1}, \dots, \zeta_{n})$.
We  have the following analyticity of the correlation functions:

\begin{prop}
Let $(F, m, \rho)$ 
be a conformal full field algebra over 
$V^L\otimes V^R$. 
For any $n\in \Z_{+}$ and $u_{1}, \dots, u_{n}$, 
there exists a multivalued analytic function
of $(z_{1}, \zeta_{1}, \dots, 
z_{n}, \zeta_{n})\in \mathbb{F}_{n}(\C)\times \mathbb{F}_{n}(\C)$ such that 
for $(z_{1}, \dots, 
z_{n})\in \mathbb{F}_{n}(\C)$, the values 
$$m_{n}(u_{1}, 
\dots, u_{n}; z_{1}, \bar{z}_{1}, \dots, 
z_{n}, \bar{z}_{n})$$
of the correlation function
is a value of this multivalued analytic function
above at the point $(z_{1}, \bar{z}_{1}, \dots, 
z_{n}, \bar{z}_{n})$. Moreover, these 
multivalued analytic functions are determined uniquely by the 
products of analytic full vertex operators in their 
regions of convergence. 
\end{prop}
\pf
The proof of this result is basically the same as the proof of 
the generalized rationality for intertwining operator algebras in 
\cite{H4}. We have proved the above strong versions of 
associativity and commutativity for analytic full vertex operators. 
Using these strong versions of associativity and commutativity,
we see that the multivalued analytic functions in various regions
obtained from 
all kinds of products and iterates of analytic full vertex operators 
are analytic extensions of each other. Thus we have such a 
global multivalued analytic function. Clearly 
these 
multivalued analytic functions are determined uniquely by the 
products of analytic full vertex operators in their 
regions of convergence.
\epfv

By the results above, we see that for a conformal
full field algebra over $V^L\otimes V^R$, the correlation function 
maps are determined uniquely by the 
products of analytic full vertex operators in their 
regions of convergence, and thus are determined uniquely by
the full vertex operator map. In view of this fact, we shall 
use also $(F, \mathbb{Y}, \rho)$ to denote 
a conformal full field algebra over 
$V^L\otimes V^R$.

We shall use 
\begin{equation}\label{anal-ext}
E(m)_{n}(u_{1}, \dots, u_{n}; z_{1}, \zeta_{1}, \dots, 
z_{n}, \zeta_{n})
\end{equation}
to denote the analytic extension obtained in the proposition above
together with the prefered values 
$$m_{n}(u_{1}, 
\dots, u_{n}; z_{1}, \bar{z}_{1}, \dots, 
z_{n}, \bar{z}_{n})$$ 
at the special points of the form $(z_{1}, \bar{z}_{1}, \dots, 
z_{n}, \bar{z}_{n})$. For $u_{1}, \dots, u_{n}\in F$ and 
 a path 
\begin{eqnarray*}
\gamma: [0, 1]&\to &\mathbb{F}_{n}(\C)\times 
\mathbb{F}_{n}(\C)\nn
t&\mapsto& (z_{1}(t), \zeta_{1}(t), \dots, z_{n}(t), \zeta_{n}(t))
\end{eqnarray*}
starting from a point of the form 
$(z_{1}, \bar{z}_{1}, \dots, 
z_{n}, \bar{z}_{n})$, we shall use 
$$E(m)_{n}(u_{1}, \dots, u_{n}; z_{1}(t), \zeta_{1}(t), \dots, 
z_{n}(t), \zeta_{n}(t))$$ 
to denote the value of (\ref{anal-ext}) at the point 
$\gamma(t)$ obtained by 
analytically extend the preferred value of (\ref{anal-ext}) 
at the starting 
point $\gamma(0)$ 
of $\gamma$ along the path $\gamma$ to the point $\gamma(t)$. 

\begin{cor}\label{anal-perm-cor}
Let $(F, \mathbb{Y}, \rho)$ 
be a conformal full field algebra over 
$V^L\otimes V^R$. Let 
\begin{eqnarray*}
\gamma: [0, 1]&\to &\mathbb{F}_{n}(\C)\times 
\mathbb{F}_{n}(\C)\nn
t&\mapsto& (z_{1}(t), \zeta_{1}(t), \dots, z_{n}(t), \zeta_{n}(t))
\end{eqnarray*}
be a path starting from  a point of the form 
$(z_{1}, \bar{z}_{1}, \dots, 
z_{n}, \bar{z}_{n})$. Then we have the following 
{\it permutation property}: For $u_{1}, \dots, u_{n}\in F$ and 
$\sigma\in S_{n}$,
\begin{eqnarray}\label{anal-perm}
\lefteqn{E(m)_{n}(u_{1}, \dots, u_{n}; z_{1}(t), \zeta_{1}(t), \dots, 
z_{n}(t), \zeta_{n}(t))}\nn
&&=E(m)_{n}(u_{\sigma(1)}, \dots, u_{\sigma(n)}; z_{\sigma(1)}(t), 
\zeta_{\sigma(1)}(t), \dots, 
z_{\sigma(n)}(t), \zeta_{\sigma(n)}(t)).
\end{eqnarray}
\end{cor}
\pf
This follows immediately from the permutation 
property for full field algebras and the uniqueness of 
analytic extensions.
\epf

\begin{cor}  \label{comm}
Let $(F, \mathbb{Y}, \rho)$ 
be a conformal full field algebra over 
$V^L\otimes V^R$. Let $r_{1}, r_{2}\in \R$ satisfying 
$r_{2}>r_{1}>0$. 
Then for $u,v,w\in F$ and $w'\in F'$, 
\begin{equation} \label{comm-cft-1-2}
\langle w', \mathbb{Y}_{\rm an}(v; r_{2}, r_{2})
\mathbb{Y}_{\rm an}(u; r_{1}, r_{1} )w\rangle, 
\end{equation}
can be obtained by analytically extending the analytic function 
(which is a branch of a multivalued function)
\begin{equation} \label{comm-cft-2-1}
\langle w', \mathbb{Y}_{\rm an}(u; z_{1}, \zeta_{1})
\mathbb{Y}_{\rm an}(v; z_{2}, \zeta_{2} )w\rangle ,
\end{equation}
defined near the point $z_{1}=\zeta_{1}=r_{2}$, $z_{2}=\zeta_{2}=r_{1}$,
in the region $|z_{1}|>|z_{2}|>0$ and 
$|\zeta_{1}|
>|\zeta_{2}|>0$, 
along the path  given by
\begin{eqnarray*}
[0, 1]&\to &M^{2}\times M^{2}\nn
t&\mapsto& ((z_{1}(t), z_{2}(t)), (\zeta_{1}(t), \zeta_{2}(t))),
\end{eqnarray*}
where
\begin{eqnarray*}
z_{1}(t)&=&\frac{r_{1}+r_{2}}{2}+e^{i\pi t}\frac{r_{2}-r_{1}}{2},\nn
z_{2}(t)&=&\frac{r_{1}+r_{2}}{2}-e^{i\pi t}\frac{r_{2}-r_{1}}{2},\nn
\zeta_{1}(t)&=&\frac{r_{1}+r_{2}}{2}+e^{-i\pi t}\frac{r_{2}-r_{1}}{2},\nn
\zeta_{2}(t)&=&\frac{r_{1}+r_{2}}{2}-e^{-i\pi t}\frac{r_{2}-r_{1}}{2},
\end{eqnarray*}
to the region $|z_{2}|>|z_{1}|>0$ and 
$|\zeta_{2}|
>|\zeta_{1}|>0$ and then evaluated at $z_{1}=\zeta_{1}=r_{1}$
and $z_{2}=\zeta_{2}=r_{2}$. 
%\begin{equation} \label{comm-path-1}
%\epsfxsize  0.8\textwidth
%\epsfysize  0.15\textwidth
%\epsfbox{comm-cft.eps}
%\end{equation}
\end{cor}
\pf
By Proposition \ref{comm-0}, we know that (\ref{comm-cft-1-2})
can indeed be obtained from (\ref{comm-cft-2-1}) by analytic 
extension. 
What we need to show now is that the analytic 
extension along the path given above gives precisely 
(\ref{comm-cft-1-2}).

Since $\mathbb{Y}(\cdot; z, \bar{z})
=\mathbb{Y}_{\rm an}(\cdot; z, \zeta)|_{\zeta=\bar{z}}$ and 
$\zeta_{1}(t)=\overline{z_{1}(t)}$ and $\zeta_{2}(t)=\overline{z_{2}(t)}$,
we see that (\ref{comm-cft-1-2}) is equal to 
$\langle w', m_{3}(v, u, w; r_{2}, r_{2}, r_{1}, r_{1}, 0, 0)\rangle$
and that 
\begin{equation}\label{comm-cft-2-1-1}
\langle w', \mathbb{Y}_{\rm an}(u; z_{1}(t), \zeta_{1}(t))
\mathbb{Y}_{\rm an}(v; z_{2}(t), \zeta_{2}(t))w\rangle
\end{equation}
when $|z_{2}(t)|>|z_{1}(t)|>0$
and 
\begin{equation}\label{comm-cft-2-1-2}
\langle w', \mathbb{Y}_{\rm an}(v; z_{2}(t), \zeta_{2}(t))
\mathbb{Y}_{\rm an}(u; z_{1}(t), \zeta_{1}(t))w\rangle
\end{equation}
when $|z_{1}(t)|>|z_{2}(t)|>0$ are
equal to 
$$\langle w', E(m)_{3}(u, v, w; z_{1}(t), \zeta_{1}(t), 
z_{2}(t), \zeta_{2}(t), 0, 0)\rangle$$
and 
$$\langle w', E(m)_{3}(v, u, w; z_{2}(t), \zeta_{2}(t), z_{1}(t), 
\zeta_{1}(t), 0, 0)\rangle,$$
respectively. 
By the permutation property of full field algebras and 
Corollary \ref{anal-perm-cor}), we see that 
(\ref{comm-cft-1-2}), (\ref{comm-cft-2-1-1}) and (\ref{comm-cft-2-1-2})
are equal to 
$$\langle w', m_{3}(u, v, w; 
r_{1}, r_{1}, r_{2}, r_{2}, 0, 0)\rangle,$$
$$\langle w', E(m)_{3}(u, v, w; z_{1}(t), \zeta_{1}(t), 
z_{2}(t), \zeta_{2}(t), 0, 0)\rangle$$ 
and 
$$\langle w', E(m)_{3}(u, v, w; z_{1}(t), \zeta_{1}(t), z_{2}(t), 
\zeta_{2}(t), 0, 0)\rangle,$$
respectively. From this fact, we see that indeed 
the analytic extension of (\ref{comm-cft-2-1}) near the point 
$z_{1}=\zeta_{1}=r_{2}$, $z_{2}=\zeta_{2}=r_{1}$,
along the path given above gives (\ref{comm-cft-1-2}).

This result can also be proved directly using the 
associativity (Proposition \ref{anal-assoc}) 
and the skew-symmetry (\ref{skew-formal}) (see \cite{K} for details). 
\epf

\begin{thm}  \label{r-conf-alg-2-prop}
A conformal full field algebra over 
$V^L\otimes V^R$ is equivalent to a module $F$ 
for the vertex operator algebra $V^L\otimes V^R$ equipped with
an intertwining operator
$\Y$ of type $\binom{F}{FF}$ and an injective 
linear map $\rho: V^L\otimes V^R \to F$, 
satisfying the following conditions:
\begin{enumerate}

\item The {\it identity property}: 
$\Y(\rho(\one^L\otimes \one^R), x)=I_{F}$.

\item The {\it creation property}: For $u\in F$, 
$\lim_{x\rightarrow 0} \Y(u, x)\rho(\one^L\otimes \one^R)=u$.

\item The {\it associativity}: The equality (\ref{asso-z-zeta}) holds when
$|z_{1}|>|z_{2}|>0$ and $|\zeta_{1}|>|\zeta_{2}|>0$.

\item The {\it single-valuedness property}: 
\begin{equation} \label{sing-val-1}
e^{2\pi i (L^L(0)-L^R(0))} = I_{F}.
\end{equation}   

\item The {\it skew symmetry}: 
\begin{equation} \label{skew-1-1}
\mathbb{Y}^{\Y}(u; 1, 1)v = e^{L^L(-1)+L^R(-1)}
\mathbb{Y}^{\Y}(v; e^{\pi i}, e^{-\pi i})u.
\end{equation}
\end{enumerate}
\end{thm}
\pf 
If $(F, \mathbb{Y}, \rho)$ is a conformal full field algebra over 
$V^L\otimes V^R$, then the results in Section 1 shows that 
$F$ is a $V^L\otimes V^R$-module, $\mathbb{Y}_{f}(\cdot; x, x)$ 
is an intertwining operator of type $\binom{F}{FF}$ and the 
five conditions are all satisfied. We now prove the converse.

Let $F$ be a module for $V^L\otimes V^R$, $\Y$ 
an intertwining operator  of type $\binom{F}{FF}$ and 
$\rho: V^L\otimes V^R \to F$ an injective linear
map, satisfying the five conditions above. 
We take the splitting $\mathbb{Y}^{\Y}$ of $\Y$ to be the 
full vertex operator map. For simplicity, we shall denote 
$\mathbb{Y}^{\Y}$ simply by $\mathbb{Y}$.
We now want to construct the 
maps $m_{n}$ for $n\in \N$ and to 
verify the convergence property. 

Using (\ref{ch-ach-z-zeta}) and the convergence property
of the intertwining operators for the vertex operator algebras
$V^{L}$ and $V^{R}$, we know that for $u_{1}, \dots, u_{n}\in F$
and $w'\in F'$,
\begin{equation}\label{prod}
\langle w', \mathbb{Y}(u_{1}; z_{1}, \zeta_{1})\cdots
\mathbb{Y}(u_{n}; z_{n}, \zeta_{n})\one\rangle
\end{equation}
is absolutely 
convergent when $|z_{1}|>\cdots>|z_{n}|>0$, 
$|\zeta_{1}|>\cdots>|\zeta_{n}|>0$,
and can be analytically 
extended to a (possibly multivalued) analytic function of 
$z_{1}, \dots, z_{n}, \zeta_{1}, \dots, \zeta_{n}$ in the region 
given by $z_{i}\ne z_{j}$, $z_{i}\ne 0$, $\zeta_{i}\ne \zeta_{j}$, 
$\zeta_{i}\ne 0$. We use 
\begin{equation}\label{correl-fn-ex}
E(m)_{n}(w', u_{1}, \dots, u_{n}; 
z_{1}, \zeta_{1}, \dots, z_{n}, \zeta_{n})
\end{equation}
to denote this function. This is a function of $z_{1}, \zeta_{1}, 
\dots, z_{n}, \zeta_{n}$
where $(z_{1}, \dots, z_{n})$, $(\zeta_{1}, \dots, \zeta_{n})\in 
\mathbb{F}_{n}(\C)$. 
So we can view this function as a function on $\mathbb{F}_{n}(\C)\times \mathbb{F}_{n}(\C)$.
In general, this function is multivalued. 
Using analytic extension,
a value of this function at a point $\mathbb{P}_{1}
\in \mathbb{F}_{n}(\C)\times \mathbb{F}_{n}(\C)$ and a path $\gamma$
in $\mathbb{F}_{n}(\C)\times \mathbb{F}_{n}(\C)$ from $\mathbb{P}_{1}$ to $\mathbb{P}_{2}\in 
\mathbb{F}_{n}(\C)\times \mathbb{F}_{n}(\C)$, determines uniquely 
a value of the function at the point 
$\mathbb{P}_{2}$. Moreover, this value depends only on the 
homotopy class of the path $\gamma$. We shall call the value of the 
function (\ref{correl-fn-ex}) at 
$\mathbb{P}_{2}$ obtained this way the {\it value 
of (\ref{correl-fn-ex}) at $\mathbb{P}_{2}$  obtained by analytically 
extending the value of (\ref{correl-fn-ex}) at $\mathbb{P}_{1}$ along 
$\gamma$}. 

We choose the correlation function
\begin{equation}\label{correl-fn}
\langle w', m_{n}(u_{1}, \dots, u_{n}; 
z_{1}, \bar{z}_{1}, \dots, z_{n}, \bar{z}_{n})\rangle
\end{equation}
as follows: For $z_{1}=n, \dots, z_{n}=1$,
we define (\ref{correl-fn}) to be (\ref{prod}) with 
$z_{1}=\zeta_{1}=n, \dots,
z_{n}=\zeta_{n}=1$. 
For general $(z_{1}, \dots, z_{n})\in \mathbb{F}_{n}(\C)$, 
we choose a path $\gamma$ from $(n, \dots, 1)$ 
to $(z_{1}, \dots, z_{n})$. Then we 
have a path $\gamma\times \overline{\gamma}$
from 
$$((n, \dots, 1), (n, \dots, 1))\in 
\mathbb{F}_{n}(\C)\times \mathbb{F}_{n}(\C)$$
to 
$$((z_{1}, \dots, z_{n}), (\bar{z}_{1}, \dots, \bar{z}_{n}))\in 
\mathbb{F}_{n}(\C)\times \mathbb{F}_{n}(\C).$$
We define 
(\ref{correl-fn}) at to be the value of (\ref{correl-fn-ex})
obtained by analytically extending the value of (\ref{correl-fn-ex}) at 
$((n, \dots, 1), (n, \dots, 1))$  along $\gamma\times \overline{\gamma}$.

The first thing we have to prove is that the correlation function 
we just defined is indeed independent of the path $\gamma$. To prove 
this fact, we need only prove that if $\gamma$ is a loop in 
$\mathbb{F}_{n}(\C)$ based at $(n, \dots, 1)$ , then the value
of (\ref{correl-fn-ex}) 
at $((n, \dots, 1), (n, \dots, 1))$ obtained by analytically 
extending the value (\ref{correl-fn}) of (\ref{correl-fn-ex}) at 
$((n, \dots, 1), (n, \dots, 1))$ 
along the loop $\gamma\times \overline{\gamma}$ is equal to 
the original value (\ref{correl-fn}) of (\ref{correl-fn-ex}) at 
$((n, \dots, 1), (n, \dots, 1))$.
In other words, we need only prove that the {\it monodromy along 
the path $\gamma\times \overline{\gamma}$ is trivial}. 
Note that  the group of the 
homotopy classes of based loops in $\mathbb{F}_n(\C)$ , 
that is, the fundamental group 
of $\mathbb{F}_{n}(\C)$, is the pure braid group of $n$ strands (see 
\cite{Bi}). This 
group is generated by the homotopy classes of the loops 
given by fixing $z_{1}, \dots, z_{j-1}$, $z_{j+1}, \dots, z_{n}$
to be $n, \dots, n-(j-2)$, $n-j, \dots, 1$, respectively, 
and moving $z_{j}$ starting from $z_{j}=n-(j-1)$ 
around $z_{i}=n-(i-1)$ once (but 
not around other points above) in the counter clockwise direction,
for $i\ne j$, $i, j=1, \dots, n$. 
Hence we need only prove that the monodromy along 
the path $\gamma\times \overline{\gamma}$ is trivial
for (the homotopy class of) such a loop $\gamma$. 

We now prove that the monodromy along 
the path $\gamma\times \bar{\gamma}$ is trivial for 
(the homotopy classes of) such a  loop $\gamma$.
Let $r$ be a positive real 
number satisfying $n-(i-1)>r>n-i$. Note that 
$r$ satisfies $r>n-(i-1)-r>0$. We know that 
\begin{eqnarray}\label{prod-1}
\lefteqn{\langle w', \mathbb{Y}(u_{1}; n, n)\cdots
\mathbb{Y}(u_{i}; n-(i-1), n-(i-1))\mathbb{Y}(u_{j}; r, r)\cdot}\nn
&&\quad\quad\cdot 
\mathbb{Y}(u_{i+1}; n-i, n-i)\cdots
\mathbb{Y}(u_{i}; n-(j-2), n-(j-2))\cdot\nn
&&\quad\quad\quad\quad\cdot
\mathbb{Y}(u_{i+1}; n-j, n-j))\cdots
\mathbb{Y}(u_{n}; 1, 1)\one\rangle
\end{eqnarray}
can be obtained by analytically extending
the value
$$%\begin{equation}\label{prod-1.1}
\langle w', \mathbb{Y}(u_{1}; n, n)\cdots
\mathbb{Y}(u_{n}; 1, 1)\one\rangle.
$$%\end{equation}
along a path from $((n, \dots, 1), (n, \dots, 1))$
to 
\begin{equation}\label{point}
((n, \dots, n-(j-2), r, n-j, \dots,
1), (n, \dots, n-(j-2), r, n-j, \dots,
1)).
\end{equation}
Such a path can always be taken to be of the form 
$\gamma_{0}\times \bar{\gamma}_{0}$ where 
$\gamma_{0}$ is a path in $\mathbb{F}_{n}(\C)$ from 
$(n, \dots, 1)$ to $(n,  \dots, n-(j-2), r,
n-j, \dots, 1)$. 
This path $\gamma_{0}$ induces an isomorphism from the fundamental group
of $\mathbb{F}_{n}(\C)$ based at $(n, \dots, 1)$ to that based at 
$(n,  \dots, n-(j-2), r, 
n-j, \dots, 1)$.
It is clear that the monodromy along a loop
based at $(n, \dots, 1)$ is trivial if and only if the monodromy along
the corresponding  loop
based at $(n, \dots, n-(j-2), r,
n-j, \dots, 1)$ is trivial. So we need only prove that 
the monodromy along a loop of the form $\gamma\times \bar{\gamma}$ 
is trivial where $\gamma$ is a loop based at 
$(n, \dots, n-(j-2), r, 
n-j, \dots, 1)$
given by fixing $z_{1}, \dots, z_{j-1}$, $z_{j+1}, \dots, z_{n}$
to be  $n, \dots, n-(j-2)$, $n-j, \dots, 1$, respectively, 
and moving $z_{j}$ starting from $z_{j}=r$ around $z_{i}=n-(i-1)$ once 
(but not around other points above) in the counter clockwise direction.
By the definition of $\gamma_{0}$,
the value of (\ref{correl-fn-ex}) at the point (\ref{point})
obtained by
analytically extending the value (\ref{correl-fn}) 
of (\ref{correl-fn-ex}) at 
the point $((n, \dots, 1), (n, \dots, 1))$
along $\gamma_{0}$ is (\ref{prod-1}).
Since we also have $r>n-(i-1)-r>0$, by associativity, 
(\ref{prod-1}) is equal to 
\begin{eqnarray*}
\lefteqn{\langle w', \mathbb{Y}(u_{1}; n, n)\cdots
\mathbb{Y}(u_{i-1}; n-(i-2), n-(i-2))\cdot}\nn
&&\quad\cdot
\mathbb{Y}(\mathbb{Y}(u_{i}; n-(i-1)-r, n-(i-1)-r)u_{j}; 
r, r)\cdot\nn
&&\quad\quad\quad\cdot
\mathbb{Y}(u_{i+1}; n-i, n-i))\cdots 
\mathbb{Y}(u_{i}; n-(j-2), n-(j-2))\cdot\nn
&&\quad\quad\quad\quad\quad \cdot
\mathbb{Y}(u_{i+1}; n-j, n-j))\cdots
\mathbb{Y}(u_{n}; r_{n}, r_{n})\one\rangle.
\end{eqnarray*}
Now let $\gamma: [0, 1]\to \mathbb{F}_{n}(\C)$ be the loop given by 
$$%\begin{equation}\label{iter-1}
t\mapsto (n, \dots, n-(i-2),
r+e^{2\pi i t}(n-(i-1)-r), n-i, \dots, 
n-(j-2), r, n-j, \dots, 1).
$$%\end{equation} 
Then the value of (\ref{correl-fn-ex}) at (\ref{point})
obtained by
analytically extending the original value (\ref{prod-1})
of (\ref{correl-fn-ex}) at 
the point (\ref{point}) 
along $\gamma$
is
\begin{eqnarray}\label{iter-2}
\lefteqn{\langle w', \mathbb{Y}(u_{1}; n, n)\cdots
\mathbb{Y}(u_{i-1}; n-(i-2), n-(i-2))\cdot}\nn
&&\quad\cdot
\mathbb{Y}(\mathbb{Y}(u_{i}; e^{2\pi i}(n-(i-1)-r),
e^{-2\pi i}(n-(i-1)-r))u_{j}; r, r)\cdot\nn
&&\quad\quad\quad\cdot \mathbb{Y}(u_{i+1}; n-i, n-i)
\cdots 
\mathbb{Y}(u_{i}; n-(j-2), n-(j-2))\cdot\nn
&&\quad\quad\quad\quad\quad\cdot
\mathbb{Y}(u_{i+1}; n-j, n-j))\cdots
\mathbb{Y}(u_{n}; r_{n}, r_{n})\one\rangle.
\end{eqnarray} 
But by the $L^{L}(0)$- and $L^{R}(0)$-conjugation properties and 
the siungle-valuedness property, we have
\begin{eqnarray}\label{single-v}
\lefteqn{\mathbb{Y}(u_{i}; e^{2\pi i}(n-(i-1)-r),
e^{-2\pi i}(n-(i-1)-r))}\nn
&&=e^{2\pi i(L^{L}(0)-L^{R}(0))}
\mathbb{Y}(e^{-2\pi i(L^{L}(0)-L^{R}(0))}u_{i}; n-(i-1)-r, 
n-(i-1)-r)\cdot\nn
&&\quad\quad\cdot e^{-2\pi i(L^{L}(0)-L^{R}(0))}\nn
&&=\mathbb{Y}(u_{i}; n-(i-1)-r, 
n-(i-1)-r).
\end{eqnarray}
Using (\ref{single-v}) and the associativity again, 
we see that (\ref{iter-2}) is equal to 
(\ref{prod-1}). Thus the analytic extension along this loop indeed
gives trivial monodromy.

Now the correlation functions and thus the maps 
$m_{n}$ for $n\in \N$ are defined. The only remaining thing to be shown is 
the convergence property.
We need to show that for any $k\in \Z_+$, $l_1, \dots, l_k\in \Z_+$,
$(z_1, \dots, z_k)\in \mathbb{F}_{n}(\C)$, 
$(z_1^{(i)}, \dots, z_{l_{i}}^{(i)}) \in 
\mathbb{F}_{l_{i}}(\C)$, $i=1, \dots, k$,
the series (\ref{series-case-2-0})
converges absolutely to 
(\ref{sum}) when 
$|z_p^{(i)}| + |z_q^{(j)}|< |z^{(0)}_i-z^{(0)}_j|$ 
for $i\ne j$, $i, j=1, \dots, k$,  
$p=1, \dots, l_i$ and $q=1, \dots, l_j$.

We use induction on $k$. 
We first prove the special case in which 
$k=2$ and $z_{2}^{(0)}=\bar{z}_{2}^{(0)}=0$. The case $k=1$ is in 
fact a special case. By the definition 
of $\mathbb{Y}$, (\ref{series-case-2-0}) becomes
\begin{eqnarray}\label{conv-pf-0}
\lefteqn{\sum_{p_1,q_1, p_{2}, q_{2}}
\mathbb{Y}(P_{p_1,q_1}m_{l_1}(u_1^{(1)}, \dots, u_{l_1}^{(1)}; 
z_1^{(1)}, \bar{z}_1^{(1)},
\dots, z_{l_1}^{(1)}, \bar{z}_{l_1}^{(1)}); 
z_1^{(0)}, \bar{z}_1^{(0)})\cdot}\nn
&&\quad\quad\quad\quad\quad\quad\cdot
P_{p_2,q_2}m_{l_2}(u_1^{(2)}, \dots, u_{l_2}^{(2)}; 
z_1^{(2)}, \bar{z}_1^{(2)},
\dots, z_{l_2}^{(2)}, \bar{z}_{l_2}^{(2)})).
\end{eqnarray}
We use induction on $l_{1}$. When $l_{1}=1$, (\ref{conv-pf-0})
becomes
\begin{eqnarray}\label{conv-pf-0.1}
\lefteqn{\sum_{p_1,q_1, p_{2}, q_{2}}
\mathbb{Y}(P_{p_1,q_1}\mathbb{Y}(u_1^{(1)}; 
z_1^{(1)}, \bar{z}_1^{(1)})\one; 
z_1^{(0)}, \bar{z}_1^{(0)})\cdot}\nn
&&\quad\quad\quad\quad\quad\quad\cdot
P_{p_2,q_2}m_{l_2}(u_1^{(2)}, \dots, u_{l_2}^{(2)}; 
z_1^{(2)}, \bar{z}_1^{(2)},
\dots, z_{l_2}^{(2)}, \bar{z}_{l_2}^{(2)})).
\end{eqnarray}
Using the construction of $\mathbb{Y}$ in terms of intertwining 
operators, the 
properties of intertwining operators and noticing that 
our condition $|z_{1}^{(1)}|+|z_{j}^{(2)}|<|z_{1}^{(0)}|$ implies
$|z_{1}^{(1)}|<|z_{1}^{(0)}|$ and 
$|z_{1}^{(1)}+z_{1}^{(0)}|>|z_{j}^{(2)}|$, we know that 
\begin{eqnarray*}\label{conv-pf-0.2}
\lefteqn{\sum_{p_1,q_1}
\mathbb{Y}(P_{p_1,q_1}\mathbb{Y}(u_1^{(1)}; 
z_1^{(1)}, \bar{z}_1^{(1)})\one; 
z_1^{(0)}, \bar{z}_1^{(0)})\cdot}\nn
&&\quad\quad\quad\quad\quad\quad\cdot
P_{p_2,q_2}m_{l_2}(u_1^{(2)}, \dots, u_{l_2}^{(2)}; 
z_1^{(2)}, \bar{z}_1^{(2)},
\dots, z_{l_2}^{(2)}, \bar{z}_{l_2}^{(2)}))
\end{eqnarray*}
is absolutely convergent to 
\begin{eqnarray*}\label{conv-pf-0.3}
\lefteqn{\mathbb{Y}(u_1^{(1)}; 
z_1^{(1)}+z_1^{(0)}, \bar{z}_1^{(1)}+\bar{z}_1^{(0)})
\mathbb{Y}(\one; 
z_1^{(0)}, \bar{z}_1^{(0)})\cdot}\nn
&&\quad\quad\quad\quad\quad\cdot
P_{p_2,q_2}m_{l_2}(u_1^{(2)}, \dots, u_{l_2}^{(2)}; 
z_1^{(2)}, \bar{z}_1^{(2)},
\dots, z_{l_2}^{(2)}, \bar{z}_{l_2}^{(2)}))\nn
&&=\mathbb{Y}(u_1^{(1)}; 
z_1^{(1)}+z_1^{(0)}, \bar{z}_1^{(1)}+\bar{z}_1^{(0)})
\cdot\nn
&&\quad\quad\quad\quad\quad\cdot
P_{p_2,q_2}m_{l_2}(u_1^{(2)}, \dots, u_{l_2}^{(2)}; 
z_1^{(2)}, \bar{z}_1^{(2)},
\dots, z_{l_2}^{(2)}, \bar{z}_{l_2}^{(2)})).
\end{eqnarray*}
Then by the construction of the correlation function maps and, in particular, 
by the fact that the correlation functions are values of 
multivalued analytic functions at certain particular points,
we know that
the right-hand side of (\ref{conv-pf-0.1}) is absolutely convergent to 
$$m_{1+l_{2}}(u_1^{(1)}, u_1^{(2)}, \dots, u_{l_2}^{(2)};
z_1^{(1)}+z_1^{(0)}, \bar{z}_1^{(1)}+\bar{z}_1^{(0)}, 
z_1^{(2)}, \bar{z}_1^{(2)},
\dots, z_{l_2}^{(2)}, \bar{z}_{l_2}^{(2)}).$$
Here we have used the fact that if a certain iterated sum 
of a series in powers of these complex variables 
is convergent to an analytic functions in the region above, then 
the multisum 
must also be absolutely convergent.

Now we assume that for $l_{1}<l$, the conclusion holds. We want to prove 
the conclusion for the case
$l_{1}=l$. We first assume that $|z_{1}^{(1)}|, \dots, |z_{l}^{(1)}|$
are all different from each other. Then in particular 
there exists $t$ such that 
$|z_{t}^{(1)}|>|z_{1}^{(1)}|, \dots, \widehat{|z_{t}^{(1)}|}, 
\dots, |z_{l}^{(1)}|$, where and also below we use $\hat{~}$ to denote 
that the item under $\hat{~}$ is missing. 
Then (\ref{conv-pf-0}) in this case is equal
to 
\begin{eqnarray}\label{conv-pf-0.4}
\lefteqn{\sum_{p_1,q_1, p_{2}, q_{2}}\sum_{r, s}
\mathbb{Y}(P_{p_1,q_1}\mathbb{Y}(u_t^{(1)}; 
z_t^{(1)}, \bar{z}_t^{(1)})\cdot}\nn
&&\cdot
P_{r, s}
m_{l-1}(u_1^{(1)}, \dots, \widehat{u_{t}^{(1)}},  \dots,
u_{l}^{(1)}; 
z_1^{(1)}, \bar{z}_1^{(1)}, \dots, \widehat{z_{t}^{(1)}},
\widehat{\bar{z}_{t}^{(1)}},
\dots, z_{l}^{(1)}, \bar{z}_{l}^{(1)}); \nn
&&\quad\quad\quad\quad\quad\quad\quad\quad\quad\quad
\quad\quad\quad\quad\quad\quad
\quad\quad\quad\quad\quad\quad\quad\quad\quad\quad
z_1^{(0)}, \bar{z}_1^{(0)})\cdot\nn
&&\quad\quad\quad\quad\quad\quad\cdot
P_{p_2,q_2}m_{l_2}(u_1^{(2)}, \dots, u_{l_2}^{(2)}; 
z_1^{(2)}, \bar{z}_1^{(2)},
\dots, z_{l_2}^{(2)}, \bar{z}_{l_2}^{(2)})).
\end{eqnarray}
Using the construction of $\mathbb{Y}$ in terms of intertwining 
operators and the 
properties of intertwining operators, and noticing that 
our condition $|z_{t}^{(1)}|+|z_{j}^{(2)}|<|z_{1}^{(0)}|$ implies
$|z_{t}^{(1)}|<|z_{1}^{(0)}|$ and 
$|z_{t}^{(1)}+z_{1}^{(0)}|>|z_{j}^{(2)}|$, we know that 
\begin{eqnarray}\label{conv-pf-0.5}
\lefteqn{\sum_{p_1,q_1}
\mathbb{Y}(P_{p_1,q_1}\mathbb{Y}(u_t^{(1)}; 
z_t^{(1)}, \bar{z}_t^{(1)})\cdot}\nn
&&\cdot
P_{r, s}
m_{l-1}(u_1^{(1)}, \dots, \widehat{u_{t}^{(1)}},  \dots,
u_{l}^{(1)}; 
z_1^{(1)}, \bar{z}_1^{(1)}, \dots, \widehat{z_{t}^{(1)}},
\widehat{\bar{z}_{t}^{(1)}},
\dots, z_{l}^{(1)}, \bar{z}_{l}^{(1)}); \nn
&&\quad\quad\quad\quad\quad\quad\quad\quad\quad\quad
\quad\quad\quad\quad\quad\quad
\quad\quad\quad\quad\quad\quad\quad\quad\quad\quad
z_1^{(0)}, \bar{z}_1^{(0)})\cdot\nn
&&\quad\quad\quad\quad\quad\quad\cdot
P_{p_2,q_2}m_{l_2}(u_1^{(2)}, \dots, u_{l_2}^{(2)}; 
z_1^{(2)}, \bar{z}_1^{(2)},
\dots, z_{l_2}^{(2)}, \bar{z}_{l_2}^{(2)}))
\end{eqnarray}
is convergent absolutely and, when $|z_t^{(1)}+z_{1}^{(0)}|>|z_{1}^{(0)}|$,
it is absolutely convergent to
\begin{eqnarray}\label{conv-pf-0.6}
\lefteqn{\mathbb{Y}(u_t^{(1)}; 
z_t^{(1)}+z_{1}^{(0)}, \bar{z}_t^{(1)}+\bar{z}_{1}^{(0)})\cdot}\nn
&&\cdot \mathbb{Y}(P_{r, s}
m_{l-1}(u_1^{(1)}, \dots, \widehat{u_{t}^{(1)}},  \dots,
u_{l}^{(1)}; 
z_1^{(1)}, \bar{z}_1^{(1)}, \dots, \widehat{z_{t}^{(1)}},
\widehat{\bar{z}_{t}^{(1)}},
\dots, z_{l}^{(1)}, \bar{z}_{l}^{(1)}); \nn
&&\quad\quad\quad\quad\quad\quad\quad\quad\quad\quad
\quad\quad\quad\quad\quad\quad
\quad\quad\quad\quad\quad\quad\quad\quad\quad\quad
z_1^{(0)}, \bar{z}_1^{(0)})\cdot\nn
&&\quad\quad\quad\quad\quad\quad\cdot
P_{p_2,q_2}m_{l_2}(u_1^{(2)}, \dots, u_{l_2}^{(2)}; 
z_1^{(2)}, \bar{z}_1^{(2)},
\dots, z_{l_2}^{(2)}, \bar{z}_{l_2}^{(2)})).
\end{eqnarray}
By the induction assumption, 
\begin{eqnarray}\label{conv-pf-0.7}
\lefteqn{\sum_{r, s, p_{2}, q_{2}}\mathbb{Y}(u_t^{(1)}; 
z_t^{(1)}+z_{1}^{(0)}, \bar{z}_t^{(1)}+\bar{z}_{1}^{(0)})\cdot}\nn
&&\cdot \mathbb{Y}(P_{r, s}
m_{l-1}(u_1^{(1)}, \dots, \widehat{u_{t}^{(1)}},  \dots,
u_{l}^{(1)}; 
z_1^{(1)}, \bar{z}_1^{(1)}, \dots, \widehat{z_{t}^{(1)}},
\widehat{\bar{z}_{t}^{(1)}},
\dots, z_{l}^{(1)}, \bar{z}_{l}^{(1)}); \nn
&&\quad\quad\quad\quad\quad\quad\quad\quad\quad\quad
\quad\quad\quad\quad\quad\quad
\quad\quad\quad\quad\quad\quad\quad\quad\quad\quad
z_1^{(0)}, \bar{z}_1^{(0)})\cdot\nn
&&\quad\quad\quad\quad\quad\quad\cdot
P_{p_2,q_2}m_{l_2}(u_1^{(2)}, \dots, u_{l_2}^{(2)}; 
z_1^{(2)}, \bar{z}_1^{(2)},
\dots, z_{l_2}^{(2)}, \bar{z}_{l_2}^{(2)}))
\end{eqnarray}
is absolutely convergent to 
\begin{eqnarray}\label{conv-pf-0.8}
\lefteqn{\mathbb{Y}(u_t^{(1)}; 
z_t^{(1)}+z_{1}^{(0)}, \bar{z}_t^{(1)}+\bar{z}_{1}^{(0)})\cdot}\nn
&&\cdot 
m_{l+l_{2}-1}(u_1^{(1)}, \dots, \widehat{u_{t}^{(1)}},  \dots,
u_{l}^{(1)}, u_1^{(2)}, \dots, u_{l_2}^{(2)}; 
z_1^{(1)}+z_1^{(0)}, \bar{z}_1^{(1)}+\bar{z}_1^{(0)}, \dots, \nn
&&\quad
\widehat{z_{t}^{(1)}+z_1^{(0)}},
\widehat{\bar{z}_{t}^{(1)}+\bar{z}_1^{(0)}},
\dots, z_{l}^{(1)}+z_1^{(0)}, \bar{z}_{l}^{(1)}+\bar{z}_1^{(0)}% \nn
%&&\quad\quad\quad\quad\quad\quad\cdot
z_1^{(2)}, \bar{z}_1^{(2)},
\dots, z_{l_2}^{(2)}, \bar{z}_{l_2}^{(2)})\nn
&&=m_{l+l_{2}}(u_1^{(1)}, \dots, 
u_{l}^{(1)}, u_1^{(2)}, \dots, u_{l_2}^{(2)}; 
z_1^{(1)}+z_1^{(0)}, \bar{z}_1^{(1)}+\bar{z}_1^{(0)}, \dots, \nn
&&\quad\quad\quad\quad\quad\quad\quad\quad\quad\quad
z_{l}^{(1)}+z_1^{(0)}, \bar{z}_{l}^{(1)}+\bar{z}_1^{(0)}% \nn
%&&\quad\quad\quad\quad\quad\quad\cdot
z_1^{(2)}, \bar{z}_1^{(2)},
\dots, z_{l_2}^{(2)}, \bar{z}_{l_2}^{(2)}).
\end{eqnarray}
We know that the right-hand side of (\ref{conv-pf-0.7}) is a value of 
the multivalued analytic function 
\begin{eqnarray*}
\lefteqn{E(m)_{l+l_{2}}(u_1^{(1)}, \dots, 
u_{l}^{(1)}, u_1^{(2)}, \dots, u_{l_2}^{(2)}; 
z_1^{(1)}+z_1^{(0)}, \zeta_1^{(1)}+\zeta_1^{(0)}, \dots,} \nn
&&\quad\quad\quad\quad\quad\quad\quad\quad\quad\quad
z_{l}^{(1)}+z_1^{(0)}, \zeta_{l}^{(1)}+\zeta_1^{(0)}% \nn
%&&\quad\quad\quad\quad\quad\quad\cdot
z_1^{(2)}, \zeta_1^{(2)},
\dots, z_{l_2}^{(2)}, \zeta_{l_2}^{(2)})
\end{eqnarray*}
at the points satisfying $\zeta_1^{(0)}=\bar{z}_1^{(0)}$,
$\zeta_{p}^{(i)}=\bar{z}_{p}^{(i)}$ for $p=1, \dots, l_{i}$,
$i=1, 2$. Since both the sum of (\ref{conv-pf-0.7}) and 
the right-hand side of (\ref{conv-pf-0.8}) are values of 
multivalued analytic functions in the same region and 
we have proved that 
their values are equal when $|z_t^{(1)}+z_{1}^{(0)}|>|z_{1}^{(0)}|$, 
(\ref{conv-pf-0.7}) must be convergent 
absolutely to
the right-hand side of (\ref{conv-pf-0.8}) even when 
$|z_t^{(1)}+z_{1}^{(0)}|>|z_{1}^{(0)}|$ is not satisfied. 
By the properties of analytic functions, we know that 
(\ref{conv-pf-0.4}) as a sum in a different order  is also 
convergent absolutely to the right-hand side of (\ref{conv-pf-0.8}).

Now we discuss the case that some of $|z_{1}^{(1)}|, \dots, |z_{l}^{(1)}|$
are equal. Let $N(z_{1}^{(1)}, \dots, z_{l}^{(1)})$ be the subset of
$\{z_{1}^{(1)}, \dots, z_{l}^{(1)}\}$ consisting of those elements
whose absolute values are equal to the absolute 
values of some other elements of $\{z_{1}^{(1)}, \dots, z_{l}^{(1)}\}$.
We use induction on the number of elements of 
$N(z_{1}^{(1)}, \dots, z_{l}^{(1)})$.
When the number is $0$, this is the case discussed above. 
Now assume that when the number is equal to $n$, the conclusion holds.
When this number is equal to $n+1$, 
let $\epsilon$ be a complex number such that the number of elements
of $N(z_{1}^{(1)}+\epsilon, \dots, z_{l}^{(1)}+\epsilon)$
is $n$ and 
$|z_{p}^{(1)}+\epsilon|+|z_{q}^{(2)}+\epsilon|<|z_1^{(0)}|$ for 
$p=1, \dots, l$ and $q=1, \dots, l_{2}$. 
Note that we can always find such an $\epsilon$ and we can 
take such an $\epsilon$ with $|\epsilon|$ to be arbitrarily small. 

By induction assumption, 
\begin{eqnarray}\label{conv-pf-0.9}
\lefteqn{\sum_{p_1,q_1, p_{2}, q_{2}}
\mathbb{Y}(P_{p_1,q_1}m_{l}(u_1^{(1)}, \dots, u_{l}^{(1)}; 
z_1^{(1)}+\epsilon, \bar{z}_1^{(1)}+\bar{\epsilon},
\dots, z_{l}^{(1)}+\epsilon, \bar{z}_{l}^{(1)}+\bar{\epsilon});} \nn
&&\quad\quad\quad\quad\quad\quad\quad\quad\quad\quad
\quad\quad\quad\quad\quad\quad
\quad\quad\quad\quad\quad\quad\quad\quad\quad\quad
z_1^{(0)}, \bar{z}_1^{(0)})\cdot\nn
&&\quad\quad\quad\cdot
P_{p_2,q_2}m_{l_2}(u_1^{(2)}, \dots, u_{l_2}^{(2)}; 
z_1^{(2)}+\epsilon, \bar{z}_1^{(2)}+\bar{\epsilon},
\dots, z_{l_2}^{(2)}+\epsilon, \bar{z}_{l_2}^{(2)}+\bar{\epsilon})\nn
&&
\end{eqnarray}
is absolutely convergent to 
\begin{eqnarray}\label{conv-pf-0.10}
\lefteqn{m_{l+l_{2}}(u_1^{(1)}, \dots, 
u_{l}^{(1)}, u_1^{(2)}, \dots, u_{l_2}^{(2)}; 
z_1^{(1)}+z_1^{(0)}+\epsilon, 
\bar{z}_1^{(1)}+\bar{z}_1^{(0)}+\bar{\epsilon}, \dots,} \nn
&&\quad\quad
z_{l}^{(1)}+z_1^{(0)}+\epsilon, 
\bar{z}_{l}^{(1)}+\bar{z}_1^{(0)}+\bar{\epsilon}% \nn
%&&\quad\quad\quad\quad\quad\quad\cdot
z_1^{(2)}+\epsilon, \bar{z}_1^{(2)}+\bar{\epsilon},
\dots, z_{l_2}^{(2)}+\epsilon, \bar{z}_{l_2}^{(2)}+\bar{\epsilon}).\nn
&&
\end{eqnarray}
We have
\begin{eqnarray}\label{conv-pf-0.115}
\lefteqn{\sum_{r_1,s_1, r_{2}, s_{2}}\sum_{p_1,q_1, p_{2}, q_{2}}
\langle e^{-\epsilon L^{L}(1)-\bar{\epsilon}L^{R}(1)}u', 
\mathbb{Y}(P_{r_{1}, s_{1}}e^{\epsilon L^{L}(-1)+\bar{\epsilon}L^{R}(-1)}
\cdot}\nn
&&\quad\quad\quad\quad\quad\quad\quad\cdot
P_{p_1,q_1}m_{l}(u_1^{(1)}, \dots, u_{l}^{(1)}; 
z_1^{(1)}, \bar{z}_1^{(1)},
\dots, z_{l}^{(1)}, \bar{z}_{l}^{(1)}); 
z_1^{(0)}, \bar{z}_1^{(0)})\cdot\nn
&&\quad\quad\;\;\cdot
P_{r_{1}, s_{1}}e^{\epsilon L^{L}(-1)+\bar{\epsilon}L^{R}(-1)}
P_{p_2,q_2}m_{l_2}(u_1^{(2)}, \dots, u_{l_2}^{(2)}; 
z_1^{(2)}, \bar{z}_1^{(2)},
\dots, z_{l_2}^{(2)}, \bar{z}_{l_2}^{(2)})\rangle\nn
&&=\sum_{r_1,s_1, r_{2}, s_{2}}
\langle e^{-\epsilon L^{L}(1)-\bar{\epsilon}L^{R}(1)}u', \nn
&&\quad\quad\quad
\mathbb{Y}(P_{r_{1}, s_{1}}m_{l}(u_1^{(1)}, \dots, u_{l}^{(1)}; 
z_1^{(1)}+\epsilon, \bar{z}_1^{(1)}+\bar{\epsilon},
\dots, z_{l}^{(1)}+\epsilon, \bar{z}_{l}^{(1)}+\bar{\epsilon}); \nn
&&\quad\quad\quad\quad\quad\quad\quad\quad\quad\quad
\quad\quad\quad\quad\quad\quad
\quad\quad\quad\quad\quad\quad\quad\quad\quad\quad
z_1^{(0)}, \bar{z}_1^{(0)})\cdot\nn
&&\quad\quad\quad\quad\cdot
P_{r_{1}, s_{1}}m_{l_2}(u_1^{(2)}, \dots, u_{l_2}^{(2)}; 
z_1^{(2)}+\epsilon, \bar{z}_1^{(2)}+\bar{\epsilon},
\dots, z_{l_2}^{(2)}+\epsilon, \bar{z}_{l_2}^{(2)}+\bar{\epsilon})\rangle.\nn
&&
\end{eqnarray}
So the right-hand side and thus also the 
left-hand side of (\ref{conv-pf-0.115}) is absolutely 
convergent to 
\begin{eqnarray}\label{conv-pf-0.12}
\lefteqn{\langle e^{-\epsilon L^{L}(1)-\bar{\epsilon}L^{R}(1)}u', 
m_{l+l_{2}}(u_1^{(1)}, \dots, 
u_{l}^{(1)}, u_1^{(2)}, \dots, u_{l_2}^{(2)}; 
z_1^{(1)}+z_1^{(0)}+\epsilon, } \nn
&&\quad\quad\quad\quad\quad\quad\quad\quad
\bar{z}_1^{(1)}+\bar{z}_1^{(0)}+\bar{\epsilon}, \dots,
z_{l}^{(1)}+z_1^{(0)}+\epsilon, 
\bar{z}_{l}^{(1)}+\bar{z}_1^{(0)}+\bar{\epsilon}, \nn
&&\quad\quad\quad\quad\quad\quad\quad\quad\quad\quad\quad\quad\quad\quad
z_1^{(2)}+\epsilon, \bar{z}_1^{(2)}+\bar{\epsilon},
\dots, z_{l_2}^{(2)}+\epsilon, \bar{z}_{l_2}^{(2)}+\bar{\epsilon})\rangle\nn
&&=\langle u', e^{-\epsilon L^{L}(-1)-\bar{\epsilon}L^{R}(-1)}
m_{l+l_{2}}(u_1^{(1)}, \dots, 
u_{l}^{(1)}, u_1^{(2)}, \dots, u_{l_2}^{(2)}; 
z_1^{(1)}+z_1^{(0)}+\epsilon,  \nn
&&\quad\quad\quad\quad\quad\quad\quad\quad
\bar{z}_1^{(1)}+\bar{z}_1^{(0)}+\bar{\epsilon}, \dots,
z_{l}^{(1)}+z_1^{(0)}+\epsilon, 
\bar{z}_{l}^{(1)}+\bar{z}_1^{(0)}+\bar{\epsilon}, \nn
&&\quad\quad\quad\quad\quad\quad\quad\quad\quad\quad\quad\quad\quad\quad
z_1^{(2)}+\epsilon, \bar{z}_1^{(2)}+\bar{\epsilon},
\dots, z_{l_2}^{(2)}+\epsilon, \bar{z}_{l_2}^{(2)}+\bar{\epsilon})\rangle\nn
&&=\langle u', 
m_{l+l_{2}}(u_1^{(1)}, \dots, 
u_{l}^{(1)}, u_1^{(2)}, \dots, u_{l_2}^{(2)}; 
z_1^{(1)}+z_1^{(0)}, 
\bar{z}_1^{(1)}+\bar{z}_1^{(0)}, \dots, \nn
&&\quad\quad\quad\quad\quad\quad\quad\quad\quad\quad
z_{l}^{(1)}+z_1^{(0)}, 
\bar{z}_{l}^{(1)}+\bar{z}_1^{(0)}% \nn
%&&\quad\quad\quad\quad\quad\quad\cdot
z_1^{(2)}, \bar{z}_1^{(2)},
\dots, z_{l_2}^{(2)}, \bar{z}_{l_2}^{(2)})\rangle.\nn
\end{eqnarray}
Since the left-hand side of (\ref{conv-pf-0.12}) is 
a value of a multivalued analytic function, any of its expansion must 
be absolutely convergent. In particular, the left-hand side of 
(\ref{conv-pf-0.115}) as an expansion of the left-hand side of 
(\ref{conv-pf-0.12}) is absolutely convergent. Thus we can exchange the 
order of the two summation signs such that the resulting 
series is still absolutely convergent to 
the left-hand side of (\ref{conv-pf-0.12}) and thus to 
the right-hand side of (\ref{conv-pf-0.12}).

But for $u'\in F'$, 
\begin{eqnarray}\label{conv-pf-0.11}
\lefteqn{\sum_{p_1,q_1, p_{2}, q_{2}}
\langle u', \mathbb{Y}(P_{p_1,q_1}m_{l}(u_1^{(1)}, \dots, u_{l}^{(1)}; 
z_1^{(1)}, \bar{z}_1^{(1)},
\dots, z_{l}^{(1)}, \bar{z}_{l}^{(1)}); } \nn
&&\quad\quad\quad\quad\quad\quad\quad\quad\quad\quad
\quad\quad\quad\quad\quad\quad
\quad\quad\quad\quad\quad\quad\quad\quad\quad\quad
z_1^{(0)}, \bar{z}_1^{(0)})\cdot\nn
&&\quad\quad\quad\quad\quad\quad\quad\quad\cdot
P_{p_2,q_2}m_{l_2}(u_1^{(2)}, \dots, u_{l_2}^{(2)}; 
z_1^{(2)}, \bar{z}_1^{(2)},
\dots, z_{l_2}^{(2)}, \bar{z}_{l_2}^{(2)})\rangle\nn
&&=\sum_{p_1,q_1, p_{2}, q_{2}}
\langle e^{-\epsilon L^{L}(1)-\bar{\epsilon}L^{R}(1)}u', 
e^{\epsilon L^{L}(-1)+\bar{\epsilon}L^{R}(-1)}\cdot\nn
&&\quad\quad\quad\quad\quad\cdot
\mathbb{Y}(P_{p_1,q_1}m_{l}(u_1^{(1)}, \dots, u_{l}^{(1)}; 
z_1^{(1)}, \bar{z}_1^{(1)},
\dots, z_{l}^{(1)}, \bar{z}_{l}^{(1)}); 
z_1^{(0)}, \bar{z}_1^{(0)})\cdot\nn
&&\quad\quad\quad\quad\quad\quad\quad\quad\quad\quad\cdot
P_{p_2,q_2}m_{l_2}(u_1^{(2)}, \dots, u_{l_2}^{(2)}; 
z_1^{(2)}, \bar{z}_1^{(2)},
\dots, z_{l_2}^{(2)}, \bar{z}_{l_2}^{(2)})\rangle\nn
&&=\sum_{p_1,q_1, p_{2}, q_{2}}
\langle e^{-\epsilon L^{L}(1)-\bar{\epsilon}L^{R}(1)}u', 
\mathbb{Y}(e^{\epsilon L^{L}(-1)+\bar{\epsilon}L^{R}(-1)}\cdot\nn
&&\quad\quad\quad\quad\quad\quad\quad\cdot
P_{p_1,q_1}m_{l}(u_1^{(1)}, \dots, u_{l}^{(1)}; 
z_1^{(1)}, \bar{z}_1^{(1)},
\dots, z_{l}^{(1)}, \bar{z}_{l}^{(1)}); 
z_1^{(0)}, \bar{z}_1^{(0)})\cdot\nn
&&\quad\quad\quad\quad\quad\cdot
e^{\epsilon L^{L}(-1)+\bar{\epsilon}L^{R}(-1)}
P_{p_2,q_2}m_{l_2}(u_1^{(2)}, \dots, u_{l_2}^{(2)}; 
z_1^{(2)}, \bar{z}_1^{(2)},
\dots, z_{l_2}^{(2)}, \bar{z}_{l_2}^{(2)})\rangle\nn
&&=\sum_{p_1,q_1, p_{2}, q_{2}}\sum_{r_1,s_1, r_{2}, s_{2}}
\langle e^{-\epsilon L^{L}(1)-\bar{\epsilon}L^{R}(1)}u', 
\mathbb{Y}(P_{r_{1}, s_{1}}e^{\epsilon L^{L}(-1)+\bar{\epsilon}L^{R}(-1)}
\cdot\nn
&&\quad\quad\quad\quad\quad\quad\quad\cdot
P_{p_1,q_1}m_{l}(u_1^{(1)}, \dots, u_{l}^{(1)}; 
z_1^{(1)}, \bar{z}_1^{(1)},
\dots, z_{l}^{(1)}, \bar{z}_{l}^{(1)}); 
z_1^{(0)}, \bar{z}_1^{(0)})\cdot\nn
&&\quad\quad\;\;\cdot
P_{r_{1}, s_{1}}e^{\epsilon L^{L}(-1)+\bar{\epsilon}L^{R}(-1)}
P_{p_2,q_2}m_{l_2}(u_1^{(2)}, \dots, u_{l_2}^{(2)}; 
z_1^{(2)}, \bar{z}_1^{(2)},
\dots, z_{l_2}^{(2)}, \bar{z}_{l_2}^{(2)})\rangle.\nn
&&
\end{eqnarray}
We have shown that the right-hand side of (\ref{conv-pf-0.11})
is absolute convergent to the right-hand side of (\ref{conv-pf-0.12}).
Thus the left-hand side of (\ref{conv-pf-0.11})
is also absolute convergent to the right-hand side 
of (\ref{conv-pf-0.12}). So we have proved the convergence 
property when the number of elements of 
$N(z_{1}^{(1)}, \dots, z_{l}^{(1)})$ is $n+1$. 
Thus the convergence property is proved when 
some of $|z_{1}^{(1)}|, \dots, |z_{l}^{(1)}|$
are equal. 

By induction principle, we 
have proved the convergence property in this special case.

We now assume that when $k<K$, 
(\ref{series-case-2-0})
converges absolutely to 
(\ref{sum}) when $z_{p}^{(0)}\ne 
z_{q}^{(0)}$ for $p, q=1, \dots, K$, $z_{p}^{(i)}\ne 
z_{q}^{(i)}$ for  $p, q=1, \dots, l_{i}$ and $i=1, \dots, K$
$1\leq p, q\leq l_i$, and
$|z_p^{(i)}| + |z_q^{(j)}|< |z^{(0)}_i-z^{(0)}_j|$ 
for $p=1, \dots, l_i$, $q=1, \dots, l_j$,
$i, j=1, \dots, K$, $i\ne j$. Now we consider the 
case $k=K$. We first consider the case that $z_p^{(i)}\in \R_{+}\cup \{0\}$ 
for $p=1, \dots, l_{i}$ and $i=0, \dots, K$ and 
$z_{1}^{(0)}>\cdots >z_{K}^{(0)}$. 
By the definition of the correlation function maps, 
we know that (\ref{series-case-2-0}) in this case
is equal to 
\begin{eqnarray}   \label{conv-pf-1}
\lefteqn{\sum_{p_1,q_1, \dots, p_K, q_K}\sum_{r, s}
\mathbb{Y}(P_{p_1,q_1}m_{l_1}(u_1^{(1)}, \dots, u_{l_1}^{(1)}; 
z_1^{(1)}, \bar{z}_1^{(1)},
\dots, z_{l_1}^{(1)}, \bar{z}_{l_1}^{(1)}); 
z_1^{(0)}, \bar{z}_1^{(0)})\cdot} \nn
&&\cdot
P_{r, s}m_{K-1}(P_{p_{2}, q_{2}}m_{l_{2}}(u_1^{(2)},\dots, u_{l_2}^{(2)};
z_1^{(2)}, 
\bar{z}_1^{(2)}, 
\dots, z_{l_2}^{(2)}, \bar{z}_{l_2}^{(2)}), \dots, \nn
&&\quad\quad\quad\quad\quad\quad\quad\quad
P_{p_K,q_K}m_{l_K}(u_1^{(K)},\dots, u_{l_K}^{(K)}; z_1^{(K)}, 
\bar{z}_1^{(K)}, 
\dots, z_{l_K}^{(K)}, \bar{z}_{l_K}^{(K)}); \nn
&&\quad\quad\quad\quad\quad\quad\quad\quad\quad\quad\quad\quad 
\quad\quad\quad\quad\quad\quad \quad\quad\quad 
z_2^{(0)}, \bar{z}_2^{(0)},
\dots, 
z_K^{(0)}, \bar{z}_K^{(0)}).\nn
&&
\end{eqnarray}
Using 
the induction assumption, we have
\begin{eqnarray}   \label{conv-pf-2}
\lefteqn{\sum_{r, s, p_1,q_1}\sum_{p_{2}, q_{2}, \dots, p_K, q_K}
\mathbb{Y}(P_{p_1,q_1}m_{l_1}(u_1^{(1)}, \dots, u_{l_1}^{(1)}; 
z_1^{(1)}, \bar{z}_1^{(1)},
\dots, z_{l_1}^{(1)}, \bar{z}_{l_1}^{(1)}); 
z_1^{(0)}, \bar{z}_1^{(0)})\cdot} \nn
&&\cdot
P_{r, s}m_{K-1}(P_{p_{2}, q_{2}}m_{l_{2}}(u_1^{(2)},\dots, u_{l_2}^{(2)};
z_1^{(2)}, 
\bar{z}_1^{(2)}, 
\dots, z_{l_2}^{(2)}, \bar{z}_{l_2}^{(2)}), \dots, \nn
&&\quad\quad\quad\quad\quad\quad\quad\quad
P_{p_K,q_K}m_{l_K}(u_1^{(K)},\dots, u_{l_K}^{(K)}; z_1^{(K)}, 
\bar{z}_1^{(K)}, 
\dots, z_{l_K}^{(K)}, \bar{z}_{l_K}^{(K)}); \nn
&&\quad\quad\quad\quad\quad\quad\quad\quad\quad\quad\quad\quad 
\quad\quad\quad\quad\quad\quad \quad\quad\quad 
z_2^{(0)}, \bar{z}_2^{(0)},
\dots, 
z_K^{(0)}, \bar{z}_K^{(0)})\nn
&&=\sum_{r, s, p_1,q_1}
\mathbb{Y}(P_{p_1,q_1}m_{l_1}(u_1^{(1)}, \dots, u_{l_1}^{(1)}; 
z_1^{(1)}, \bar{z}_1^{(1)},
\dots, z_{l_1}^{(1)}, \bar{z}_{l_1}^{(1)}); 
z_1^{(0)}, \bar{z}_1^{(0)})\cdot \nn
&&\quad\quad\cdot
P_{r, s}m_{l_{2}+\cdots+l_{K}}(u_1^{(2)},\dots, u_{l_2}^{(2)}, 
\dots, u_1^{(K)},\dots, u_{l_K}^{(K)};
z_1^{(2)}+z_2^{(0)}, \nn
&&\quad \quad \quad\quad\quad\quad\quad\quad\quad
\bar{z}_1^{(2)}+\bar{z}_2^{(0)}, 
\dots, 
z_{l_2}^{(2)}+z_2^{(0)}, \bar{z}_{l_2}^{(2)}+\bar{z}_2^{(0)}, \dots, \nn
&&\quad \quad \quad\quad\quad\quad\quad
z_1^{(K)}+z_K^{(0)}, 
\bar{z}_1^{(K)}+\bar{z}_K^{(0)}, 
\dots, z_{l_K}^{(K)}+z_K^{(0)}, \bar{z}_{l_K}^{(K)}+\bar{z}_K^{(0)}).\nn
&&
\end{eqnarray}
Since $z_p^{(1)} + z_q^{(i)}< z^{(0)}_1-z^{(0)}_i$
for $p=1, \dots, l_{1}$, $q=1, \dots, l_{i}$ and $i=2, \dots, K$, 
we have $z_p^{(1)} + (z_q^{(i)}+z^{(0)}_i)< z^{(0)}_1-0$. Thus
by the special case we proved above, 
the right-hand side of (\ref{conv-pf-2}) is 
absolutely convergent to
\begin{eqnarray}\label{conv-pf-2.5}
\lefteqn{m_{l_{1}+\cdots +l_{K}}(u_1^{(1)},\dots, u_{l_1}^{(1)}, 
\dots, u_1^{(K)},\dots, u_{l_K}^{(K)};
z_1^{(1)}+z_1^{(0)}, }\nn
&&\quad \quad \quad\quad\quad\quad
\bar{z}_1^{(1)}+\bar{z}_1^{(0)}, 
\dots, 
z_{l_1}^{(1)}+z_1^{(0)}, \bar{z}_{l_1}^{(1)}+\bar{z}_1^{(0)}, \dots, 
z_1^{(K)}+z_K^{(0)}, \nn
&&\quad \quad \quad\quad\quad\quad\quad\quad
\bar{z}_1^{(K)}+\bar{z}_K^{(0)}, 
\dots, z_{l_K}^{(K)}+z_K^{(0)}, \bar{z}_{l_K}^{(K)}+\bar{z}_K^{(0)}).
\end{eqnarray}

Note that  (\ref{conv-pf-2.5}) 
is a value of the $\overline{F}$-valued multivalued analytic function 
\begin{eqnarray}   \label{conv-pf-3}
\lefteqn{E(m)_{l_{1}+\cdots +l_{K}}(u_1^{(1)},\dots, u_{l_1}^{(1)}, 
\dots, u_1^{(K)},\dots, u_{l_K}^{(K)};
z_1^{(1)}+z_1^{(0)},} \nn
&&\quad \quad \quad\quad\quad\quad
\zeta_1^{(1)}+\zeta_1^{(0)}, 
\dots, 
z_{l_1}^{(1)}+z_1^{(0)}, \zeta_{l_1}^{(1)}+\zeta_1^{(0)}, \dots, 
z_1^{(K)}+z_K^{(0)}, \nn
&&\quad \quad \quad\quad\quad\quad\quad\quad
\zeta_1^{(K)}+\zeta_K^{(0)}, 
\dots, z_{l_K}^{(K)}+z_K^{(0)}, \zeta_{l_K}^{(K)}+\zeta_K^{(0)}).
\end{eqnarray}
at the point $z_{i}^{(j)}=z_{i}^{(j)}$, $\zeta_{i}^{(j)}=\bar{z}_{i}^{(j)}$.
Thus its expansions, no matter in which ways, must be convergent absolutely. 
In particular, (\ref{conv-pf-1}) as one expansion of 
(\ref{conv-pf-2.5})
must be convergent absolutely to 
(\ref{conv-pf-2.5}), proving the convergence 
in this special case of the case $k=K$. 

We know that for a series in powers of several variables, if 
it is absolutely convergent when these variables are equal to 
some real numbers, then it is also convergent when the variables 
are equal to complex numbers whose absolute values are 
equal to these real numbers. Using this property, 
we see that
\begin{eqnarray}   \label{conv-pf-4}
\lefteqn{\sum_{p_1,q_1, \dots, p_K, q_K}\sum_{r, s}
\mathbb{Y}(P_{p_1,q_1}E(m)_{l_1}(u_1^{(1)}, \dots, u_{l_1}^{(1)}; 
z_1^{(1)}, \zeta_1^{(1)},
\dots, z_{l_1}^{(1)}, \zeta_{l_1}^{(1)}); 
z_1^{(0)}, \zeta_1^{(0)})\cdot} \nn
&&\cdot
P_{r, s}E(m)_{K-1}(P_{p_{2}, q_{2}}E(m)_{l_{2}}(u_1^{(2)},\dots, u_{l_2}^{(2)};
z_1^{(2)}, 
\zeta_1^{(2)}, 
\dots, z_{l_2}^{(2)}, \zeta_{l_2}^{(2)}), \dots, \nn
&&\quad\quad\quad\quad\quad\quad\quad\quad
P_{p_K,q_K}E(m)_{l_K}(u_1^{(K)},\dots, u_{l_K}^{(K)}; z_1^{(K)}, 
\zeta_1^{(K)}, 
\dots, z_{l_K}^{(K)}, \zeta_{l_K}^{(K)}); \nn
&&\quad\quad\quad\quad\quad\quad\quad\quad\quad\quad\quad\quad 
\quad\quad\quad\quad\quad\quad \quad\quad\quad 
z_2^{(0)}, \zeta_2^{(0)},
\dots, 
z_K^{(0)}, \zeta_K^{(0)}).\nn
&&
\end{eqnarray}
is convergent absolutely to a branch of (\ref{conv-pf-3}) when 
$z_{p}^{(0)}\ne 
z_{q}^{(0)}$, $\zeta_{p}^{(0)}\ne 
\zeta_{q}^{(0)}$ for $p, q=1, \dots, K$, $z_{p}^{(i)}\ne 
z_{q}^{(i)}$, $\zeta_{p}^{(i)}\ne 
\zeta_{q}^{(i)}$ for  $p, q=1, \dots, l_{i}$, $i=1, \dots, K$
$|z_p^{(i)}| + |z_q^{(j)}|< ||z^{(0)}_i|-|z^{(0)}_j||$,
$|\zeta_p^{(i)}| + |\zeta_q^{(j)}|< ||\zeta^{(0)}_i|-|\zeta^{(0)}_j||$,
for $p=1, \dots, l_i$,  $q=1, \dots, l_j$, $i, j=1, \dots, K$, 
$i\ne j$, and $|z_{1}^{(0)}|>\cdots >|z_{K}^{(0)}|$,
$|\zeta_{1}^{(0)}|>\cdots >|\zeta_{K}^{(0)}|$. 
Using the permutation property
for the correlation functions, we obtain that 
(\ref{conv-pf-4}) is 
convergent absolutely to a branch of (\ref{conv-pf-3}) when 
$|z_{\sigma(1)}^{(0)}|> \cdots> |z_{\sigma(K)}^{(0)}|$ and 
$|\zeta_{\sigma(1)}^{(0)}|> \cdots> |\zeta_{\sigma(K)}^{(0)}|$
for some $\sigma\in S_{K}$.

Now for fixed $z_{1}^{(0)}, \dots, z_{K}^{(0)}$,
$\zeta_{1}^{(0)}, \dots, \zeta_{K}^{(0)}$ satisfying 
$|z_{\sigma(1)}^{(0)}|> \cdots> |z_{\sigma(K)}^{(0)}|$ and 
$|\zeta_{\sigma(1)}^{(0)}|> \cdots> |\zeta_{\sigma(K)}^{(0)}|$
for some $\sigma\in S_{K}$, any branch of 
(\ref{conv-pf-3})
can be expanded as a series in powers of
$z_{p}^{(i)}$ and $\zeta_{p}^{(i)}$, $p=1, \dots, l_{i}$, 
$i=1, \dots, K$ in the region
\begin{eqnarray}\label{conv-pf-5}
\lefteqn{\{(z_{1}^{(1)}, \dots, z_{l_{K}}^{(K)}, 
\zeta_{1}^{(1)}, \dots, \zeta_{l_{K}}^{(K)})\;|\;z_{p}^{(i)}\ne 
z_{q}^{(i)}, \zeta_{p}^{(i)}\ne 
\zeta_{q}^{(i)}}\nn
&&\quad\quad\quad\quad\quad\quad\quad\quad\quad\quad
\quad\quad\quad\quad\quad \;\mbox{\rm for}\;  p, q=1, \dots, l_{i}, 
i=1, \dots, K,\nn
&&\quad |z_p^{(i)}| + |z_q^{(j)}|< |z^{(0)}_i-z^{(0)}_j|,
|\zeta_p^{(i)}| + |\zeta_q^{(j)}|< |\zeta^{(0)}_i-\zeta^{(0)}_j|,\nn
&&\quad\quad\quad\quad\quad\quad\quad
\;\mbox{\rm for}\;
p=1, \dots, l_i, q=1, \dots, l_j, i, j=1, \dots, K, 
i\ne j\}.\nn
&&
\end{eqnarray}
But in the region 
\begin{eqnarray}\label{conv-pf-6}
\lefteqn{\{(z_{1}^{(1)}, \dots, z_{l_{K}}^{(K)}, 
\zeta_{1}^{(1)}, \dots, \zeta_{l_{K}}^{(K)})\;|\;z_{p}^{(i)}\ne 
z_{q}^{(i)}, \zeta_{p}^{(i)}\ne 
\zeta_{q}^{(i)}}\nn
&&\quad\quad\quad\quad\quad\quad\quad\quad\quad\quad
\quad\quad\quad\quad\quad \;\mbox{\rm for}\;  p, q=1, \dots, l_{i}, 
i=1, \dots, K,\nn
&&\quad |z_p^{(i)}| + |z_q^{(j)}|< ||z^{(0)}_i|-|z^{(0)}_j||,
|\zeta_p^{(i)}| + |\zeta_q^{(j)}|< ||\zeta^{(0)}_i|-|\zeta^{(0)}_j||,\nn
&&\quad\quad\quad\quad\quad\quad\quad
\;\mbox{\rm for}\;
p=1, \dots, l_i, q=1, \dots, l_j, i, j=1, \dots, K, 
i\ne j\},\nn
&&
\end{eqnarray}
we have proved that one branch of (\ref{conv-pf-3}) can be 
expanded as the series (\ref{conv-pf-4}), which can be further
expanded as a series in powers of 
$z_{p}^{(i)}$ and $\zeta_{p}^{(i)}$, $p=1, \dots, l_{i}$, 
$i=1, \dots, K$ in this region. Since the region (\ref{conv-pf-6})
is contained in the region (\ref{conv-pf-5}) 
and the coefficients of the expansion
can be determined completely using the values of the branch 
in the region (\ref{conv-pf-6}), we see that the restriction to 
the region (\ref{conv-pf-6}) of the expansion 
in the region (\ref{conv-pf-5}) is the same as the expansion 
in the region (\ref{conv-pf-5}). Thus, the series (\ref{conv-pf-4})
is convergent absolutely to a branch of (\ref{conv-pf-3}) 
in the region (\ref{conv-pf-5}). 

In the region (\ref{conv-pf-5}), when $\zeta_{p}^{(0)}=z_{p}^{(0)}\in \R$
for $p=1, \dots, K$, $\zeta_{p}^{(i)}=z_{p}^{(i)}\in \R$ for 
$p=1, \dots, l_{i}$, $i=1, \dots, K$, we have proved that 
(\ref{series-case-2-0}) in this case is convergent absolutely to 
the right-hand side of (\ref{conv-pf-2}). Thus in the region 
(\ref{conv-pf-5}),
(\ref{series-case-2-0}) with $k=K$ is convergent absolutely to 
the right-hand side of (\ref{conv-pf-2}),  the value of a branch 
of (\ref{conv-pf-3}).

Finally we consider the case that some of $|z_{1}^{(0)}|, \dots,
|z_{K}^{(0)}|$ are equal. Recall 
the subset $N(z_{1}^{(0)}, \dots, z_{K}^{(0)})$ of
$\{z_{1}^{(0)}, \dots, z_{K}^{(0)}\}$ consisting of those elements
whose absolute values are equal to the absolute 
values of some other elements of $\{z_{1}^{(0)}, \dots, z_{K}^{(0)}\}$.
We use induction on the number of elements of 
$N(z_{1}^{(0)}, \dots, z_{K}^{(0)})$.
When the number is $0$, this is the case discussed above. 
Now assume that when the number is equal to $n$, the conclusion holds.
When this number is equal to $n+1$, 
let $\epsilon$ be a complex number such that the number of elements
of $N(z_{1}^{(0)}+\epsilon, \dots, z_{K}^{(0)}+\epsilon)$
is $n$ and that the other conditions are still 
satisfied. 
Note that we can always find such an $\epsilon$ and we can 
take such an $\epsilon$ with $|\epsilon|$ to be arbitrary small. 
By induction assumption, 
\begin{eqnarray*}  
\lefteqn{\sum_{p_1,q_1, \dots, p_K, q_K} 
m_K(P_{p_1,q_1}m_{l_1}(u_1^{(1)},\dots, u_{l_1}^{(1)}; 
z_1^{(1)}, \bar{z}_1^{(1)}, 
\dots, z_{l_1}^{(1)}, \bar{z}_{l_1}^{(1)}), 
\dots,} \nn
&&\quad\quad\quad\quad\quad\quad
P_{p_K,q_K}m_{l_K}(u_1^{(K)},\dots, u_{l_K}^{(K)}; z_1^{(K)}, 
\bar{z}_1^{(K)}, 
\dots, z_{l_K}^{(K)}, \bar{z}_{l_K}^{(K)}); \nn
&&\quad\quad\quad\quad\quad\quad\quad\quad\quad\quad\quad\quad\quad
\quad
z_1^{(0)}+\epsilon, \bar{z}_1^{(0)}+\bar{\epsilon}, \dots, 
z_K^{(0)}+\epsilon, \bar{z}_K^{(0)}+\bar{\epsilon})
\end{eqnarray*}
is absolutely convergent to 
\begin{eqnarray*}  \label{conv-pf-8}
\lefteqn{
m_{l_{1}+\cdots +l_{K}}(u_1^{(1)},\dots, u_{l_1}^{(1)}, 
\dots, u_1^{(K)},\dots, u_{l_K}^{(K)};
z_1^{(1)}+z_1^{(0)}+\epsilon, }\nn
&&\quad 
\bar{z}_1^{(1)}+\bar{z}_1^{(0)}+\bar{\epsilon}, 
\dots, 
z_{l_1}^{(1)}+z_1^{(0)}+\epsilon, 
\bar{z}_{l_1}^{(1)}+\bar{z}_1^{(0)}+\bar{\epsilon}, \dots, 
z_1^{(K)}+z_K^{(0)}+\epsilon, \nn
&&\quad \quad \quad\quad\quad\quad\quad\quad\quad 
\bar{z}_1^{(K)}+\bar{z}_K^{(0)}+\bar{\epsilon}, 
\dots, z_{l_K}^{(K)}+z_K^{(0)}+\epsilon, 
\bar{z}_{l_K}^{(K)}+\bar{z}_K^{(0)}+\bar{\epsilon}).
\end{eqnarray*}
Thus for $u'\in F$, 
\begin{eqnarray*}  \label{conv-pf-9}
\lefteqn{\sum_{p_1,q_1, \dots, p_K, q_K} 
\langle u', m_K(P_{p_1,q_1}m_{l_1}(u_1^{(1)},\dots, u_{l_1}^{(1)}; 
z_1^{(1)}, \bar{z}_1^{(1)}, 
\dots, z_{l_1}^{(1)}, \bar{z}_{l_1}^{(1)}), 
\dots,} \nn
&&\quad  \quad\quad\quad\quad\quad\quad
P_{p_K,q_K}m_{l_K}(u_1^{(K)},\dots, u_{l_K}^{(K)}; z_1^{(K)}, 
\bar{z}_1^{(K)}, 
\dots, z_{l_K}^{(K)}, \bar{z}_{l_K}^{(K)}); \nn
&&\quad\quad\quad\quad\quad\quad\quad\quad\quad\quad\quad\quad\quad
\quad\quad\quad\quad\quad\quad\quad\quad
z_1^{(0)}, \bar{z}_1^{(0)}, \dots, 
z_K^{(0)}, \bar{z}_K^{(0)})\rangle\nn
&&=\sum_{p_1,q_1, \dots, p_K, q_K} 
\langle e^{\epsilon L^{L}(1)+\bar{\epsilon}L^{R}(1)}u', 
e^{-\epsilon L^{L}(-1)-\bar{\epsilon}L^{R}(-1)}\cdot\nn
&&\quad\quad\quad\quad\cdot
m_K(P_{p_1,q_1}m_{l_1}(u_1^{(1)},\dots, u_{l_1}^{(1)}; 
z_1^{(1)}, \bar{z}_1^{(1)}, 
\dots, z_{l_1}^{(1)}, \bar{z}_{l_1}^{(1)}), 
\dots, \nn
&&\quad \quad\quad\quad\quad\quad\quad
P_{p_K,q_K}m_{l_K}(u_1^{(K)},\dots, u_{l_K}^{(K)}; z_1^{(K)}, 
\bar{z}_1^{(K)}, 
\dots, z_{l_K}^{(K)}, \bar{z}_{l_K}^{(K)}); \nn
&&\quad\quad\quad\quad\quad\quad\quad\quad\quad\quad\quad\quad\quad
\quad\quad\quad\quad\quad\quad\quad\quad
z_1^{(0)}, \bar{z}_1^{(0)}, \dots, 
z_K^{(0)}, \bar{z}_K^{(0)})\rangle\nn
&&=\sum_{p_1,q_1, \dots, p_K, q_K} 
\langle e^{\epsilon L^{L}(1)+\bar{\epsilon}L^{R}(1)}u', \nn
&&\quad\quad\quad\quad
m_K(P_{p_1,q_1}m_{l_1}(u_1^{(1)},\dots, u_{l_1}^{(1)}; 
z_1^{(1)}, \bar{z}_1^{(1)}, 
\dots, z_{l_1}^{(1)}, \bar{z}_{l_1}^{(1)}), 
\dots, \nn
&&\quad \quad\quad\quad\quad\quad\quad
P_{p_K,q_K}m_{l_K}(u_1^{(K)},\dots, u_{l_K}^{(K)}; z_1^{(K)}, 
\bar{z}_1^{(K)}, 
\dots, z_{l_K}^{(K)}, \bar{z}_{l_K}^{(K)}); \nn
&&\quad\quad\quad\quad\quad\quad\quad\quad\quad\quad\quad\quad\quad
z_1^{(0)}+\epsilon, \bar{z}_1^{(0)}+\bar{\epsilon},, \dots, 
z_K^{(0)}+\epsilon,, \bar{z}_K^{(0)}+\bar{\epsilon})\rangle.\nn
&&
\end{eqnarray*}
is absolutely convergent to 
\begin{eqnarray*}  \label{conv-pf-10}
\lefteqn{\langle e^{\epsilon L^{L}(1)+\bar{\epsilon}L^{R}(1)}u', 
m_{l_{1}+\cdots +l_{K}}(u_1^{(1)},\dots, u_{l_1}^{(1)}, 
\dots, u_1^{(K)},\dots, u_{l_K}^{(K)};
z_1^{(1)}+z_1^{(0)}+\epsilon, }\nn
&&\quad 
\bar{z}_1^{(1)}+\bar{z}_1^{(0)}+\bar{\epsilon}, 
\dots, 
z_{l_1}^{(1)}+z_1^{(0)}+\epsilon, 
\bar{z}_{l_1}^{(1)}+\bar{z}_1^{(0)}+\bar{\epsilon}, \dots, 
z_1^{(K)}+z_K^{(0)}+\epsilon, \nn
&&\quad \quad \quad\quad\quad\quad\quad\quad\quad
\bar{z}_1^{(K)}+\bar{z}_K^{(0)}+\bar{\epsilon}, 
\dots, z_{l_K}^{(K)}+z_K^{(0)}+\epsilon, 
\bar{z}_{l_K}^{(K)}+\bar{z}_K^{(0)}+\bar{\epsilon})\rangle\nn
&&=\langle u', m_{l_{1}+\cdots +l_{K}}(u_1^{(1)},\dots, u_{l_1}^{(1)}, 
\dots, u_1^{(K)},\dots, u_{l_K}^{(K)};
z_1^{(1)}+z_1^{(0)}, \nn
&&\quad \quad \quad\quad\quad\quad
\bar{z}_1^{(1)}+\bar{z}_1^{(0)}, 
\dots, 
z_{l_1}^{(1)}+z_1^{(0)}, \bar{z}_{l_1}^{(1)}+\bar{z}_1^{(0)}, \dots, 
z_1^{(K)}+z_K^{(0)}, \nn
&&\quad \quad \quad\quad\quad\quad\quad\quad\quad\quad\quad\quad\quad
\bar{z}_1^{(K)}+\bar{z}_K^{(0)}, 
\dots, z_{l_K}^{(K)}+z_K^{(0)}, \bar{z}_{l_K}^{(K)}+\bar{z}_K^{(0)})
\rangle.
\end{eqnarray*}
Since $u'$ is arbitrary, we have proved that 
(\ref{series-case-2-0}) with $k=K$ is convergent absolutely to 
the right-hand side of (\ref{conv-pf-2}) in the case that 
the number of elements of $N(z_{1}^{(0)}, \dots,
z_{K}^{(0)})$ is $n+1$. Thus we have proved this conclusion 
in the case that some of $|z_{1}^{(0)}|, \dots,
|z_{K}^{(0)}|$ are equal. 

By the principle of induction, the convergence property is proved.
\epfv

%Notice that the first three conditions in above proposition exactly amounts
%to an open-string vertex operator algebra structure on $F$. 
%Hence we have the following corollary.
%\begin{cor} \label{ffa-op-voa-rema} 
%An open-string vertex operator
%algebra $F$ including $(V^L\otimes V^R)_{voa}$ in its meromorphic center 
%and satisfying the single-valuedness (\ref{sing-val-1}) and skew symmetry
%(\ref{skew-1-1}) is equivalent to a $V^L\otimes V^R$-rational 
%full field algebra structure on $F$. 
%\end{cor}

\begin{rema} {\rm
Although the definition of full field algebra in Definition \ref{ffa} is 
very general, it is not easy to verify all the axioms directly. 
Theorem \ref{r-conf-alg-2-prop}  gives
an equivalent definition of conformal full field algebra over 
$V^L\otimes V^R$ and the axioms in this definition are much easier 
to verify than those in Definitions \ref{ffa}, \ref{RR-grading-ffa}
and \ref{g-r-RR-g-ffa}. In our construction of full field algebras
in the next section, we shall use this definition to verify the 
structure we construct is indeed a full field algebra. }
\end{rema}

\renewcommand{\theequation}{\thesection.\arabic{equation}}
\renewcommand{\thethm}{\thesection.\arabic{thm}}
\setcounter{equation}{0}
\setcounter{thm}{0}

\section{A construction of full field algebras with nondegenerate invariant
bilinear forms}

Let $V$ be a simple vertex operator algebra and $C_{2}(V)$ the subspace of 
$V$ spanned by $u_{-2}v$ for $u, v\in V$. In this section,
we assume that $V$ satisfies the 
following conditions:

\begin{enumerate}

\item $V_{(n)}=0$ for $n<0$, $V_{(0)}=\mathbb{C}\mathbf{1}$
and $W_{(0)}=0$ for any irreducible $V$-module $W$ which is not equivalent
to $V$. %$V'$ as a $V$-module is isomorphic to $V$.

\item Every $\mathbb{N}$-gradable weak $V$-module is completely 
reducible.

\item $V$ is $C_{2}$-cofinite, that is, $\dim V/C_{2}(V)<\infty$. 

\end{enumerate}
(%Note that if 
%$V$ satisfies (i)' $V_{(n)}=0$ for $n<0$, $V_{(0)}=\mathbb{C}\mathbf{1}$,
%$W_{(0)}=0$ for any irreducible $V$-module $W$ which is not equivalent
%to $V$, then (i) is satisfied. 
Note %also 
that by results of Li \cite{L} and  Abe, Buhl and Dong \cite{ABD}, 
Conditions 2 and 3 can be 
replaced by a single condition that every weak $V$-module is completely 
reducible.)

Since $V$ satisfies the conditions above, all the results 
in \cite{H7} can be used. We shall use all the 
notations, conventions and choices used in this paper. 
In particular, we use the following notations and choices:
$\mathcal{A}$ is the (finite) set of
equivalence classes of irreducible $V$-modules; $e$ is the
equivalence class containing $V$;
$^{\prime}:
\mathcal{A}\to \mathcal{A}$
is the map induced from the functor given by taking contragredient 
modules; for $a\in \mathcal{A}$, 
$W^{a}$ is a representative of $a$;
$(\cdot, \cdot)$ is the nondegenerate bilinear form 
on $V$ normalized by $(\mathbf{1}, \mathbf{1})=1$; 
for $a_{1}, a_{2}, a_{3}\in \mathcal{A}$, 
$\mathcal{V}_{a_{1}a_{2}}^{a_{3}}$ 
are the spaces of intertwining operators of 
type $\binom{W^{a_{3}}}{W^{a_{1}}W^{a_{2}}}$;
$\sigma_{12}$ and $\sigma_{23}$ are actions of $(12)$ and 
$(23)$ on $\mathcal{V}$ and they
generate an action of $S_{3}$ on 
$$\mathcal{V}=\coprod_{a_{1}, a_{2}, a_{3}\in 
\mathcal{A}}\mathcal{V}_{a_{1}a_{2}}^{a_{3}};$$ 
for any bases 
$\Y_{a_{1}a_{2}; i}^{a_{3}; (p)}$, $i=1, \dots, 
N_{a_{1}a_{2}}^{a_{3}}$, $p=1, 2, 3, 4, 5, 6, \dots$ and $a_{1}, a_{2}, 
a_{3}\in \mathcal{A}$, of $\mathcal{V}_{a_{1}a_{2}}^{a_{3}}$,
$$F(\Y_{a_{1}a_{5}; i}^{a_{4}; (1)}\otimes \Y_{a_{2}a_{3}; j}^{a_{5}; (2)}; 
\Y_{a_{6}a_{3}; l}^{a_{4}; (3)}
\otimes \Y_{a_{1}a_{2}; k}^{a_{6}; (4)}) \in \C$$
are matrix elements of the fusing isomorphism;
for $a\in \A$, $\Y_{ea; 1}^{a}$, $\Y_{ae; 1}^{a}$ and 
$\Y_{aa'; 1}^{e}$ are bases of $\V_{ea}^{a}$, $\V_{ae}^{a}$ and 
$\V_{aa'}^{e}$ chosen in \cite{H7};
for $a\in \A$, there exists $h_{a}\in \Q$ such that
$W^{a}=\coprod_{n\in h_{a}+\N}W_{(n)}$.

For $a_{1}, a_{2}, a_{3}\in \A$, we now want to 
introduce a pairing between 
$\mathcal{V}_{a_{1}a_{2}}^{a_{3}}$ and $\mathcal{V}_{a'_{1}a'_{2}}^{a'_{3}}$.

For $a\in \A$, $w_{a}\in W^{a}$ and $w'_{a}\in (W^{a})'$,
we shall use $\tilde{w}_{a}$ and $\tilde{w}'_{a}$ to denote
$e^{-L(1)}w_{a}$ and $e^{-L(1)}w'_{a}$, respectively. Then we have 
$$\langle e^{L(1)}\tilde{w}'_{a}, e^{L(1)}\tilde{w}_{a}\rangle
=\langle w'_{a}, w_{a}\rangle.$$ 
We have:

\begin{lemma}\label{res}
For $a\in \A$, $w_{a}\in W^{a}$ and $w'_{a}\in (W^{a})'$,
$$\res_{z=0}z^{-1}\Y_{a'a;1}^{e}(z^{L(0)}e^{\pi i(L(0)-h_{a})}
\tilde{w}'_{a}, z)
z^{L(0)}\tilde{w}_{a}=\langle w'_{a}, w_{a}\rangle\one.$$
\end{lemma}
\pf
Since $V_{(0)}=\C\one$ and $\Y_{a'a;1}^{e}=\sigma_{23}(\Y_{a'e; 1}^{a'})$, 
\begin{eqnarray*}
\lefteqn{\res_{z=0}z^{-1}\Y_{a'a;1}^{e}(z^{L(0)}
e^{\pi i(L(0)-h_{a})}\tilde{w}'_{a}, z)
z^{L(0)}\tilde{w}_{a}}\nn
&&=(\mathbf{1}, 
\Y_{a'a;1}^{e}(1^{L(0)}e^{\pi i(L(0)-h_{a})}\tilde{w}'_{a}, 1)
\tilde{w}_{a})\one\nn
&&=(\mathbf{1}, 
\sigma_{23}(\Y_{a'e;1}^{a'})(e^{\pi i(L(0)-h_{a})}\tilde{w}'_{a_{i}}, 1)
\tilde{w}_{a})\one\nn
&&=\langle e^{\pi ih_{a}}\Y_{a'e;1}^{a'}(e^{L(1)}e^{-\pi iL(0)}
e^{\pi i(L(0)-h_{a})}
\tilde{w}'_{a}, 1)
\mathbf{1}, 
\tilde{w}_{a}\rangle \one\nn
&&=\langle \Y_{a'e;1}^{a'}(e^{L(1)}\tilde{w}'_{a}, 1)
\one, \tilde{w}_{a}\rangle\one\nn
&&=\langle e^{L(-1)}e^{L(1)}\tilde{w}'_{a}, \tilde{w}_{a}\rangle\one\nn
&&=\langle e^{L(1)}\tilde{w}'_{a}, e^{L(1)}\tilde{w}_{a}\rangle\one\nn
&&=\langle w'_{a}, w_{a}\rangle \one.
\end{eqnarray*}
\epfv

For a single-valued
branch $f_{1}(z_{1}, z_{2})$ of a multivalued
analytic function in a region $A$, we use $E(f_{1}(z_{1}, z_{2}))$ 
to denote the multivalued analytic extension together with the 
preferred branch $f_{1}(z_{1}, z_{2})$. 
Let $w_{1}=w_{1}(z_{1}, z_{2})$ and $w_{2}=w_{2}(z_{1}, z_{2})$
be a change of variables and $f_{2}(z_{1}, z_{2})$ 
a branch 
of $E(f_{1}(z_{1}, z_{2}))$ in a region $B$ containing
$w_{1}(z_{1}, z_{2})=0$
and $w_{2}(z_{1}, z_{2})=0$ such that $A\cap B\ne \emptyset$ and 
$f_{1}(z_{1}, z_{2})=f(z_{1}, z_{2})$ for $(z_{1}, z_{2})\in A\cap B$.
Then we use 
$$\res_{w_{1}=0\;|\;w_{2}}E(f_{1}(z_{1}, z_{2}))$$
to denote the coefficient of $w_{1}^{-1}$ in the expansion of 
$f_{2}(z_{1}, z_{2})$ as a series in powers of $w_{1}$  whose coefficients
are analytic functions of $w_{2}$. By definition, we have 
\begin{equation}\label{residue}
\res_{w_{1}=0\;|\;C_{1}w_{2}+C_{2}}E(f_{1}(z_{1}, z_{2}))
=\res_{w_{1}=0\;|\;w_{2}}E(f_{1}(z_{1}, z_{2}))
\end{equation}
for any $C_{1}\in \C^{\times}, 
C_{2}\in \C$ independent of $z_{1}$ and $z_{2}$.
We have:

\begin{prop}\label{pairing}
For $a_{1}, a_{2}, a_{3}\in \A$, $w_{a_{1}}\in W^{a_{1}}$, 
$w_{a_{2}}\in W^{a_{2}}$, 
$w_{a'_{1}}\in (W^{a_{1}})'$, $w'_{a_{2}}\in (W^{a_{1}})'$, 
$\Y_{1}\in 
\V_{a_{1}a_{2}}^{a_{3}}$ and $\Y_{2}\in \V_{a'_{1}a'_{2}}^{a'_{3}}$, 
there exists a constant $\langle \Y_{1},
\Y_{2}\rangle_{\V_{a_{1}a_{2}}^{a_{3}}}\in \C$  such that 
\begin{eqnarray}\label{inner}
\lefteqn{\res_{1-z_{1}-z_{2}=0\;|\;z_{2}}(1-z_{1}-z_{2})^{-1}
E(\langle e^{L(1)}\Y_{2}((1-z_{1}-z_{2})^{L(0)}
\tilde{w}'_{a_{1}}, z_{1})
\tilde{w}'_{a_{2}}, }\nn
&&\quad\quad\quad\quad\quad\quad\quad\quad\quad\quad
e^{L(1)}\Y_{1}((1-z_{1}-z_{2})^{L(0)}\tilde{w}_{a_{1}}, z_{2})
\tilde{w}_{a_{2}}\rangle)\nn
&&=\langle w'_{a_{1}}, w_{a_{1}}\rangle\langle w'_{a_{2}}, w_{a_{2}}\rangle
\langle \Y_{1},
\Y_{2}\rangle_{\V_{a_{1}a_{2}}^{a_{3}}}.
\end{eqnarray}
Explicitly, for any bases $\{\Y_{a_{1}a_{2}; i}^{a_{3}; (1)}\;|\;
i=1, \dots, N_{a_{1}a_{2}}^{a_{3}}\}$ and
$\{\Y_{a'_{1}a'_{2}; i}^{a'_{3}; (2)}\;|\;
i=1, \dots, N_{a'_{1}a'_{2}}^{a'_{3}}\}$ of 
$\V_{a_{1}a_{2}}^{a_{3}}$ and 
$\V_{a'_{1}a'_{2}}^{a'_{3}}$, respectively, and for  $m, n, k, l\in \Z_{+}$, 
$i=1, \dots, N_{a_{1}a_{2}}^{a_{3}}$ and 
$j=1, \dots, N_{a'_{1}a'_{2}}^{a'_{3}}$,
we have
\begin{eqnarray}\label{inner-fusing}
\langle \Y_{a_{1}a_{2}; j}^{a_{3}; (1)}, \Y_{a'_{1}a'_{2}; i}^{a'_{3}; (2)}
\rangle_{\V_{a_{1}a_{2}}^{a_{3}}}&=&
F(\sigma_{23}(\Y_{a'_{1}a'_{2}; i}^{a'_{3}; (2)})
\otimes \Y_{a_{1}a_{2}; j}^{a_{3}; (1)}; \Y_{ea_{2}; 1}^{a_{2}}\otimes
\Y_{a'_{1}a_{1}; 1}^{e})\nn
&=&F(\sigma_{23}(\Y_{a_{1}a_{2}; j}^{a_{3}; (1)})
\otimes \Y_{a'_{1}a'_{2}; i}^{a'_{3}; (2)}; \Y_{ea'_{2}; 1}^{a'_{2}}\otimes
\Y_{a_{1}a'_{1}; 1}^{e}).
\end{eqnarray}
\end{prop}
\pf
We prove (\ref{inner}) in the case $\Y_{1}=\Y_{a_{1}a_{2}; j}^{a_{3}; (1)}$
and $\Y_{2}=\Y_{a'_{1}a'_{2}; i}^{a'_{3}; (2)}$ for 
$i=1, \dots, N_{a'_{1}a'_{2}}^{a'_{3}}$ and 
$j=1, \dots, N_{a_{1}a_{2}}^{a_{3}}$,
respectively, or equivalently, we prove (\ref{inner-fusing}).
The general case follows 
immediately from the bilinearity in $\Y_{1}$ 
and $\Y_{2}$ of the right-hand side of (\ref{inner}). 

For $a_{1}, a_{2}\in \A$, $a_{1}, a_{2}\ne e$, let 
$\{\Y_{a_{1}a_{2}; i}^{a_{2}}\;|\;i=1, \dots, N_{a_{1}a_{2}}^{a_{2}}\}$
and $\{\Y_{a'_{1}a_{1}; i}^{a_{2}}\;|\;i=1, \dots, N_{a'_{1}a_{1}}^{a_{2}}\}$
be an arbitrary basis of $\V_{a_{1}a_{2}}^{a_{2}}$ and 
$\V_{a'_{1}a_{1}}^{a_{2}}$, respectively. 

For $w_{a_{1}}\in W^{a_{1}}$, 
$w_{a_{2}}\in W^{a_{2}}$, 
$w_{a'_{1}}\in (W^{a_{1}})'$, $w'_{a_{2}}\in (W^{a_{1}})'$, we have
\begin{eqnarray*}
\lefteqn{\res_{1-z_{1}-z_{2}=0\;|\;z_{2}}(1-z_{1}-z_{2})^{-1}
E(\langle e^{L(1)}
\Y_{a'_{1}a'_{2}; i}^{a'_{3}; (2)}
((1-z_{1}-z_{2})^{L(0)}\tilde{w}'_{a_{1}}, z_{1})
\tilde{w}'_{a_{2}}, }\nn
&&\quad\quad\quad\quad\quad\quad\quad\quad\quad\quad\quad\quad
e^{L(1)}\Y_{a_{1}a_{2}; j}^{a_{3}; (1)}
((1-z_{1}-z_{2})^{L(0)}\tilde{w}_{a_{1}}, z_{2})
\tilde{w}_{a_{2}}\rangle)\nn
&&=\res_{1-z_{1}-z_{2}=0\;|\;z_{2}}(1-z_{1}-z_{2})^{-1}\cdot\nn
&&\quad\quad\cdot
E(\langle \Y_{a'_{1}a'_{2}; i}^{a'_{3}; (2)}
(e^{(1-z_{1})L(1)}(1-z_{1})^{-2L(0)}\cdot\nn
&&\quad\quad\quad\quad\quad\quad\quad\quad\cdot
(1-z_{1}-z_{2})^{L(0)}\tilde{w}'_{a_{1}}, (1-z_{1})^{-1})
e^{L(-1)}e^{L(1)}\tilde{w}'_{a_{2}}, \nn
&&\quad\quad\quad\quad\quad\quad\quad\quad\quad \quad\quad\quad\quad\quad
\Y_{a_{1}a_{2}; j}^{a_{3}; (1)}((1-z_{1}-z_{2})^{L(0)}
\tilde{w}_{a_{1}}, z_{2})
\tilde{w}_{a_{2}}\rangle)\nn
&&=\res_{1-z_{1}-z_{2}=0\;|\;z_{2}}(1-z_{1}-z_{2})^{-1}\nn
&&\quad E(\langle 
e^{-L(1)}e^{L(1)}\tilde{w}'_{a_{2}}, 
\sigma_{23}(\Y_{a'_{1}a'_{2}; i}^{a'_{3}; (2)})((1-z_{1}-z_{2})^{L(0)}
e^{\pi i(L(0)-h_{a'_{1}})}\tilde{w}'_{a_{1}}, 1-z_{1})\cdot\nn
&&\quad\quad\quad\quad\quad\quad\quad\quad\quad\quad\quad\quad\quad\quad\cdot
\Y_{a_{1}a_{2}; j}^{a_{3}; (1)}((1-z_{1}-z_{2})^{L(0)}
\tilde{w}_{a_{1}}, z_{2})
\tilde{w}_{a_{2}}\rangle)\nn
&&=\sum_{a_{4}\in \A}\sum_{p=1}^{N_{a_{4}a_{2}}^{a_{2}}}
\sum_{q=1}^{N_{a'_{1}a_{1}}^{a_{4}}}
F(\sigma_{23}(\Y_{a'_{1}a'_{2}; i}^{a'_{3}; (2)})
\otimes \Y_{a_{1}a_{2}; j}^{a_{3}; (1)}; \Y_{a_{4}a_{2}; p}^{a_{2}}\otimes
\Y_{a'_{1}a_{1}; q}^{a_{4}})\cdot \nn
&&\quad\cdot\res_{1-z_{1}-z_{2}=0\;|\;z_{2}}(1-z_{1}-z_{2})^{-1}\cdot\nn
&&\quad\quad\quad\quad\quad\quad\cdot E(\langle 
e^{L(-1)}e^{L(1)}\tilde{w}'_{a_{2}}, 
\Y_{a_{4}a_{2}; p}^{a_{2}}(\Y_{a'_{1}a_{1}; q}^{a_{4}}
((1-z_{1}-z_{2})^{L(0)}\cdot\nn
&&\quad\quad\quad\quad\quad\quad\quad\quad\quad\quad\quad\quad
\cdot e^{\pi i(L(0)-h_{a'_{1}})}
\tilde{w}'_{a_{1}}, 1-z_{1}-z_{2})\cdot\nn
&&\quad\quad\quad\quad\quad\quad\quad\quad\quad\quad\quad\quad
\quad\quad\quad\quad
\cdot (1-z_{1}-z_{2})^{L(0)}
\tilde{w}_{a_{1}}, z_{2})
\tilde{w}_{a_{2}}\rangle)\nn
&&=F(\sigma_{23}(\Y_{a'_{1}a'_{2}; i}^{a'_{3}; (2)})
\otimes \Y_{a_{1}a_{2}; j}^{a_{3}; (1)}; \Y_{ea_{2}; p}^{a_{2}}\otimes
\Y_{a'_{1}a_{1}; 1}^{e})\cdot \nn
&&\quad\cdot\res_{1-z_{1}-z_{2}=0\;|\;z_{2}}(1-z_{1}-z_{2})^{-1}\cdot\nn
&&\quad\quad\quad\quad\quad\quad\cdot E(\langle 
e^{L(-1)}e^{L(1)}\tilde{w}'_{a_{2}}, 
\Y_{ea_{2}; 1}^{a_{2}}(\Y_{a'_{1}a_{1}; 1}^{e}
((1-z_{1}-z_{2})^{L(0)}\cdot\nn
&&\quad\quad\quad\quad\quad\quad\quad\quad\quad\quad\quad\quad\cdot
e^{\pi i(L(0)-h_{a'_{1}})}
\tilde{w}'_{a_{1}}, 1-z_{1}-z_{2})\cdot\nn
&&\quad\quad\quad\quad\quad\quad\quad\quad\quad\quad\quad\quad
\quad\quad\quad\quad
\cdot
(1-z_{1}-z_{2})^{L(0)}\tilde{w}_{a_{1}}, z_{2})
\tilde{w}_{a_{2}}\rangle)\nn
&&=F(\sigma_{23}(\Y_{a'_{1}a'_{2}; i}^{a'_{3}; (2)})
\otimes \Y_{a_{1}a_{2}; j}^{a_{3}; (1)}; \Y_{ea_{2}; p}^{a_{2}}\otimes
\Y_{a'_{1}a_{1}; 1}^{e})\cdot \nn
&&\quad\quad\quad\quad\quad\quad\quad\quad\quad\quad\quad\cdot\langle 
e^{L(-1)}e^{L(1)}\tilde{w}'_{a_{2}}, 
\Y_{ea_{2}; 1}^{a_{2}}(\langle w'_{a_{1}}, w_{a_{1}}\rangle \one, z_{2})
\tilde{w}_{a_{2}}\rangle\nn
&&=\langle w'_{a_{1}}, w_{a_{1}}\rangle
F(\sigma_{23}(\Y_{a'_{1}a'_{2}; i}^{a'_{3}; (2)})
\otimes \Y_{a_{1}a_{2}; j}^{a_{3}; (1)}; \Y_{ea_{2}; p}^{a_{2}}\otimes
\Y_{a'_{1}a_{1}; 1}^{e})\langle 
e^{L(1)}\tilde{w}'_{a_{2}}, 
e^{L(1)}\tilde{w}_{a_{2}}\rangle\nn
&&=\langle w'_{a_{1}}, w_{a_{1}}\rangle
\langle w'_{a_{2}}, w_{a_{2}}\rangle
F(\sigma_{23}(\Y_{a'_{1}a'_{2}; i}^{a'_{3}; (2)})
\otimes \Y_{a_{1}a_{2}; j}^{a_{3}; (1)}; \Y_{ea_{2}; p}^{a_{2}}\otimes
\Y_{a'_{1}a_{1}; 1}^{e}),
\end{eqnarray*}
where we have used the fact that $W_{(0)}^{a_{4}}=0$ for $a_{4}\ne e$.
This proves (\ref{inner}) and also the first equality in (\ref{inner-fusing}).
The second equality in (\ref{inner-fusing}) can be proved similarly or 
can be simply obtained using the first equality in (\ref{inner-fusing}) 
and symmetry. 
\epfv

Clearly, $\langle \Y_{1},
\Y_{2}\rangle_{\V_{a_{1}a_{2}}^{a_{3}}}$ is bilinear in $\Y_{1}$ and $\Y_{2}$.
Thus we have a pairing $\langle \cdot, \cdot\rangle_{\V_{a_{1}a_{2}}^{a_{3}}}: 
\V_{a_{1}a_{2}}^{a_{3}}\otimes \V_{a'_{1}a'_{2}}^{a'_{3}}
\to \C$.

We need the following lemma:

\begin{lemma}\label{zero}
For $a_{1}, a_{2}, a_{3}, a_{4}, a_{5}\in \A$, $\Y_{1}\in 
\V_{a'_{1}a'_{2}}^{a'_{3}}$,
$\Y_{2}\in \V_{a_{4}a_{5}}^{a_{3}}$, $w'_{a_{1}}\in (W^{a_{1}})'$,
$w'_{a_{2}}\in (W^{a_{2}})'$, $w_{a_{4}}\in W^{a_{1}}$, $w_{a_{4}}
\in W^{a_{5}}$, if $a_{1}\ne a_{4}$ or $a_{2}\ne a_{5}$, then
\begin{eqnarray*}
\lefteqn{\res_{1-z_{1}-z_{2}=0\;|\;z_{2}}(1-z_{1}-z_{2})^{-1}E(\langle e^{L(1)}
\Y_{1}((1-z_{1}-z_{2})^{L(0)}\tilde{w}'_{a_{1}}, z_{1})
\tilde{w}'_{a_{2}}, }\nn
&&\quad\quad\quad\quad\quad\quad\quad\quad\quad\quad
e^{L(1)}\Y_{2}((1-z_{1}-z_{2})^{L(0)}\tilde{w}_{a_{4}}, z_{2})
\tilde{w}_{a_{5}}\rangle)=0.
\end{eqnarray*}
\end{lemma}
\pf
Using the $L(1)$- and $L(-1)$-conjugation formulas
for intertwining operators, the definition of $\sigma_{23}$ and 
the associativity of intertwining operators, we know that there 
exist a $V$-module $W$ and 
intertwining operators $\Y_{3}$ and $\Y_{4}$ of types
$\binom{W^{a_{2}}}{WW^{a_{5}}}$ and 
$\binom{W}{W^{a'_{1}}W^{a_{4}}}$, respectively, such that
\begin{eqnarray}\label{zero-1}
\lefteqn{\res_{1-z_{1}-z_{2}=0\;|\;z_{2}}(1-z_{1}-z_{2})^{-1}E(\langle e^{L(1)}
\Y_{1}((1-z_{1}-z_{2})^{L(0)}\tilde{w}'_{a_{1}}, z_{1})
\tilde{w}'_{a_{2}}, }\nn
&&\quad\quad\quad\quad\quad\quad\quad\quad\quad\quad\quad\quad\quad\quad
e^{L(1)}\Y_{2}((1-z_{1}-z_{2})^{L(0)}\tilde{w}_{a_{4}}, z_{2})
\tilde{w}_{a_{5}}\rangle)\nn
&&=\res_{1-z_{1}-z_{2}=0\;|\;z_{2}}(1-z_{1}-z_{2})^{-1}\cdot\nn
&&\quad\quad\cdot E(\langle 
e^{L(-1)}e^{L(1)}\tilde{w}'_{a_{2}}, 
\sigma_{23}(\Y_{1})((1-z_{1}-z_{2})^{L(0)}\tilde{w}'_{a_{1}}, 1-z_{1})
\cdot\nn
&&\quad\quad\quad\quad\quad\quad\quad\quad\quad\quad\quad\quad\quad\quad
\cdot
\Y_{2}((1-z_{1}-z_{2})^{L(0)}\tilde{w}_{a_{4}}, z_{2})
\tilde{w}_{a_{5}}\rangle)\nn
&&=\res_{1-z_{1}-z_{2}=0\;|\;z_{2}}(1-z_{1}-z_{2})^{-1}\cdot\nn
&&\quad\quad\cdot E(\langle 
e^{L(-1)}e^{L(1)}\tilde{w}'_{a_{2}}, 
\Y_{3}(\Y_{4}((1-z_{1}-z_{2})^{L(0)}\tilde{w}'_{a_{1}}, 1-z_{1}-z_{2})
\cdot\nn
&&\quad\quad\quad\quad\quad\quad\quad\quad\quad\quad\quad\quad\quad\quad
\cdot
(1-z_{1}-z_{2})^{L(0)}\tilde{w}_{a_{4}}, z_{2})
\tilde{w}_{a_{5}}\rangle).
\end{eqnarray}
If $a_{1}\ne a_{4}$, $W^{a_{4}}$ is not equivalent to $W^{a_{1}}$.
Thus $\V_{a'_{1}a_{4}}^{e}=0$. So it is possible to find such
a $V$-module $W$ which
does not contain a summand equivalent to $V$. By the assumption on $V$,
we have $W_{(0)}=0$. So the right-hand side of (\ref{zero-1}) is $0$, 
proving the lemma in this case. If $a_{1}= a_{4}$, 
$\V_{a'_{1}a_{1}}^{e}$ is one-dimensional. We can choose 
$W$ to contain one and only one copy of $V$. If $a_{2}\ne a_{5}$, any
intertwining operator of type $\binom{a_{2}}{ea_{5}}$ (that is, type
$\binom{W^{a_{2}}}{VW^{a_{5}}}$) must be $0$.
So $\Y_{3}(\one, z_{2})=0$. 
Since $W_{(0)}=\C\one$,  there exists $\lambda\in \C$ such that 
the right-hand side of (\ref{zero-1})
is equal to 
$$\lambda\langle 
e^{L(-1)}e^{L(1)}\tilde{w}'_{a_{2}}, 
\Y_{3}(\one, z_{2})
\tilde{w}_{a_{2}}\rangle=0,$$
proving the lemma in the case $a_{2}\ne a_{5}$.
\epfv

As in \cite{H7}, we now choose a canonical 
basis of $\V_{a_{1}a_{2}}^{a_{3}}$ for $a_{1}, a_{2}, a_{3}\in
\A$ when one of $a_{1}, a_{2}, a_{3}$ is $e$: For $a\in \A$, we choose $\Y_{ea; 1}^{a}$ to 
be the vertex operator 
$Y_{W^{a}}$ defining the module structure on $W^{a}$ and we choose 
$\Y_{ae; 1}^{a}$ to be the intertwining operator defined using 
the action of $\sigma_{12}$, or equivalently the skew-symmetry 
in this case,
\begin{eqnarray*}
\Y_{ae; 1}^{a}(w_{a}, x)u&=&\sigma_{12}(\Y_{ea; 1}^{a})(w_{a}, x)u\nn
&=&e^{xL(-1)}\Y_{ea; 1}^{a}(u, -x)w_{a}\nn
&=&e^{xL(-1)}Y_{W^{a}}(u, -x)w_{a}
\end{eqnarray*}
for $u\in V$ and $w_{a}\in W^{a}$. 
Since $V'$ as a $V$-module is isomorphic to $V$, we have 
$e'=e$. From \cite{FHL}, we know that there is a nondegenerate
invariant 
bilinear form $(\cdot, \cdot)$ on $V$ such that $(\mathbf{1}, 
\mathbf{1})=1$. 
We choose $\Y_{aa'; 1}^{e}=\Y_{aa'; 1}^{e'}$
to be the intertwining operator defined using the action of 
$\sigma_{23}$ by
$$\Y_{aa'; 1}^{e'}=\sigma_{23}(\Y_{ae; 1}^{a}),$$
that is,
$$(u, \Y_{aa'; 1}^{e}(w_{a}, x)w_{a'})
=e^{\pi i h_{a}}\langle \Y_{ae; 1}^{a}(e^{xL(1)}(e^{-\pi i}x^{-2})^{L(0)}w_{a}, x^{-1})u, 
w_{a'}\rangle$$
for $u\in V$, $w_{a}\in W^{a}$ and $w_{a'}\in W^{a'}$. Since the actions of
$\sigma_{12}$
and $\sigma_{23}$ generate the action of $S_{3}$ on $\mathcal{V}$, we have
$$\Y_{a'a; 1}^{e}=\sigma_{12}(\Y_{aa'; 1}^{e})$$
for any $a\in \mathcal{A}$.

\begin{thm}
The pairing $\langle \cdot, \cdot\rangle_{\V_{a_{1}a_{2}}^{a_{3}}}: 
\V_{a_{1}a_{2}}^{a_{3}}\otimes \V_{a'_{1}a'_{2}}^{a'_{3}}
\to \C$ is nondegenerate. In particular, 
$N_{a'_{1}a'_{2}}^{a'_{3}}
=N_{a_{1}a_{2}}^{a_{3}}.$
\end{thm}
\pf 
For $a_{1}, a_{2}, a_{3}\in
\A$ such that one of $a_{1}, a_{2}, a_{3}$ is $e$, we have a canonical 
basis $\mathcal{Y}_{a_{1}a_{2}; 1}^{a_{3}}$ given above.  
For $a_{1}, a_{2}, a_{3}\ne e$, let
$\mathcal{Y}_{a_{1}a_{2}; i}^{a_{3}}$, $i=1, \dots, 
N_{a_{1}a_{2}}^{a_{3}}$,  be an arbitrary basis
of $\mathcal{V}_{a_{1}a_{2}}^{a_{3}}$.

For $a_{1}, a_{2}, a_{3}\in
\A$, let 
\begin{eqnarray*}
\Y_{a_{1}a_{2}; j}^{a_{3}; (1)}
&=&\sigma_{123}(\Y_{a_{2}a_{3}'; j}^{a_{1}'}),
\\
\Y_{a'_{1}a'_{2}; i}^{a'_{3}; (2)}
&=&\sigma_{23}(\Y_{a'_{1}a_{3}; i}^{a_{2}}).
\end{eqnarray*}
Then the first equality of (\ref{inner-fusing}) gives
\begin{eqnarray}\label{inner-fusing-1}
\langle \sigma_{23}(\Y_{a'_{1}a_{3}; i}^{a_{2}}),
\sigma_{123}(\Y_{a_{2}a_{3}'; j}^{a_{1}'})
\rangle_{\V_{a_{1}a_{2}}^{a_{3}}}
&=&\langle \Y_{a'_{1}a'_{2}; i}^{a'_{3}; (2)},
\Y_{a_{1}a_{2}; j}^{a_{3}; (1)}
\rangle_{\V_{a_{1}a_{2}}^{a_{3}}}\nn
&=&F(\Y_{a'_{1}a_{3}; i}^{a_{2}}
\otimes \sigma_{123}(\Y_{a_{2}a_{3}'; j}^{a_{1}'}); 
\Y_{ea_{2}; 1}^{a_{2}}\otimes
\Y_{a'_{1}a_{1}; 1}^{e}).\nn
&&
\end{eqnarray}
In \cite{H7}, the first author proved the following 
formula ((4.9) in \cite{H7}):
\begin{eqnarray}\label{formula1}
\lefteqn{\sum_{k=1}^{N_{a'_{1}a_{3}}^{a_{2}}}
F(\Y_{a_{2}e; 1}^{a_{2}}\otimes \Y_{a'_{3}a_{3}; 1}^{e}; 
\Y_{a'_{1}a_{3}; k}^{a_{2}}\otimes \Y_{a_{2}a'_{3}; j}^{a'_{1}})\cdot}\nn
&&\quad\quad\quad\quad\quad\cdot 
F(\Y_{a_{1}'a_{3}; k}^{a_{2}}\otimes 
\sigma_{123}(\Y_{a_{2}a'_{3}; i}^{a'_{1}});
\Y_{ea_{2}; 1}^{a_{2}}\otimes \Y_{a'_{1}a_{1}; 1}^{e})\nn
&&=\delta_{ij}
F(\Y_{a_{2}e; 1}^{a_{2}} \otimes \Y_{a'_{2}a_{2}; 1}^{e};
\Y_{ea_{2}; 1}^{a_{2}}\otimes \Y_{a_{2}a'_{2}; 1}^{e}).
\end{eqnarray}
In the same paper \cite{H7}, the first author also proved that 
$$F(\Y_{a_{2}e; 1}^{a_{2}} \otimes \Y_{a'_{2}a_{2}; 1}^{e};
\Y_{ea_{2}; 1}^{a_{2}}\otimes \Y_{a_{2}a'_{2}; 1}^{e})\ne 0.$$
Thus from (\ref{inner-fusing-1}) and (\ref{formula1}), 
we see that the matrix
\begin{equation}\label{bl-form-mat}
(\alpha_{ij})=(\langle \sigma_{23}(\Y_{a'_{1}a_{3}; i}^{a_{2}}),
\sigma_{123}(\Y_{a_{2}a_{3}'; j}^{a_{1}'})
\rangle_{\V_{a_{1}a_{2}}^{a_{3}}})
\end{equation}
is left invertible. Note that when $a_{1}, a_{2}, a_{3}\ne e$, 
$\Y_{a'_{1}a_{3}; i}^{a_{2}}$ and $\Y_{a_{2}a_{3}'; j}^{a_{1}'}$
in (\ref{bl-form-mat})
are arbitrary bases of $\mathcal{V}_{a'_{1}a_{3}}^{a_{2}}$
and $\mathcal{V}_{a_{2}a_{3}'}^{a_{1}'}$, respectively.

We now show that (\ref{bl-form-mat}) is also right invertible. 
By definition, the bilinear form $\langle\cdot, 
\cdot\rangle_{\V_{a_{1}a_{2}}^{a_{3}}}$ is symmetric in the sense that
$$\langle \Y_{1}, \Y_{2}\rangle_{\V_{a_{1}a_{2}}^{a_{3}}}
=\langle \Y_{2}, \Y_{1}\rangle_{\V_{a'_{1}a'_{2}}^{a'_{3}}}$$
for $a_{1}, a_{2}, a_{3}\in
\A$. So 
\begin{eqnarray}\label{nondeg-1}
\lefteqn{\langle \sigma_{23}(\Y_{a'_{1}a_{3}; i}^{a_{2}}),
\sigma_{123}(\Y_{a_{2}a_{3}'; j}^{a_{1}'})
\rangle_{\V_{a_{1}a_{2}}^{a_{3}}}}\nn
&&=\langle \sigma_{123}(\Y_{a_{2}a_{3}'; j}^{a_{1}'}),
\sigma_{23}(\Y_{a'_{1}a_{3}; i}^{a_{2}})
\rangle_{\V_{a'_{1}a'_{2}}^{a'_{3}}}\nn
&&=\langle \sigma_{23}(\sigma_{13}(\Y_{a_{2}a_{3}'; j}^{a_{1}'})),
\sigma_{123}(\sigma_{13}(\Y_{a'_{1}a_{3}; i}^{a_{2}}))
\rangle_{\V_{a'_{1}a'_{2}}^{a'_{3}}}.
\end{eqnarray}
Note that for $a_{1}, a_{2}, a_{3}\in \A$, 
$\sigma_{13}(\Y_{a_{3}'a_{2}; i}^{a_{1}'})$ is a basis of 
$\mathcal{V}_{a_{1}a_{2}}^{a_{3}}$ such that when one of the elements
$a_{1}, a_{2}, a_{3}\in \A$ is $e$, these basis elements are equal to 
the special ones we chosen above. Thus by the result we obtained above,
the matrix 
$$(\beta_{ij})=(\langle \sigma_{23}(\sigma_{13}(\Y_{a_{2}a_{3}'; i}^{a_{1}'})),
\sigma_{123}(\sigma_{13}(\Y_{a'_{1}a_{3}; j}^{a_{2}}))
\rangle_{\V_{a'_{1}a'_{2}}^{a'_{3}}})$$ 
must be left invertible. So the transpose of $(\beta_{ij})$, that is, the matrix
$$(\gamma_{kl})
=(\langle \sigma_{23}(\sigma_{13}(\Y_{a_{2}a_{3}'; l}^{a_{1}'})),
\sigma_{123}(\sigma_{13}(\Y_{a'_{1}a_{3}; k}^{a_{2}}))
\rangle_{\V_{a'_{1}a'_{2}}^{a'_{3}}}),$$ 
is right invertible. By (\ref{nondeg-1}), we see that 
(\ref{bl-form-mat})
is also right invertible. 

Now we have shown that the matrix (\ref{bl-form-mat})
is in fact invertible. This is equivalent to the 
nondegeneracy of the bilinear 
form. It also implies $N_{a'_{1}a'_{2}}^{a'_{3}}
=N_{a_{1}a_{2}}^{a_{3}}.$
\epfv

For $a\in \A$, let 
$$F_{a}=F(\Y_{ae; 1}^{a} \otimes \Y_{a'a; 1}^{e};
\Y_{ea; 1}^{a}\otimes \Y_{aa'; 1}^{e})\ne 0.$$
Then by (3.12) in \cite{H7}, $F_{a'}=F_{a}$ for $a\in \A$.

\begin{lemma}
If for $a\in \A$,
$\Y_{ea; 1}^{a}$, $\Y_{ae; 1}^{a}$, $\Y_{aa'; 1}^{e}$
are the canonical bases of $\V_{ea}^{a}$, $\V_{ae}^{a}$, $\V_{aa'}^{e}$,
respectively, chosen in \cite{H7} and above, then their dual bases in 
$\V_{ea'}^{a'}$, $\V_{a'e}^{a'}$, $\V_{a'a}^{e}$
with respect to the pairing $\langle\cdot, 
\cdot\rangle_{\V_{ea}^{a}}$, $\langle\cdot, 
\cdot\rangle_{\V_{ae}^{a}}$, $\langle\cdot, 
\cdot\rangle_{\V_{aa'}^{e}}$, respectively, are equal to 
$\Y_{ea'; 1}^{a'}$, $\Y_{a'e; 1}^{a}$, $\frac{\Y_{a'a; 1}^{e}}{F_{a}}$.
\end{lemma}
\pf
This result follows immediately from the 
definition of the canonical bases in \cite{H7}.
\epfv

We have:

\begin{prop}\label{fusing}
For $a_{1}, a_{2}, a_{3}\in \A$, let
$\{\Y_{a_{1}a_{2}; i}^{a_{3}}\;|\;i=1, \dots, N_{a_{1}a_{2}}^{a_{3}}\}$
be bases
of $\V_{a_{1}a_{2}}^{a_{3}}$ and let 
$\{\Y_{a'_{1}a'_{2}; i}^{\prime; a'_{3}}\;|\;i=1, \dots, 
N_{a'_{1}a'_{2}}^{a'_{3}}\}$ be the dual bases of 
$\{\Y_{a_{1}a_{2}; i}^{a_{3}}\;|\;i=1, \dots, N_{a_{1}a_{2}}^{a_{3}}\}$
with respect to the pairing $\langle\cdot, 
\cdot\rangle_{\V_{a_{1}a_{2}}^{a_{3}}}$. Assume that 
for $a\in \A$,
$\Y_{ea; 1}^{a}$, $\Y_{ae; 1}^{a}$, $\Y_{aa'; 1}^{e}$
are the canonical bases of $\V_{ea}^{a}$, $\V_{ae}^{a}$, $\V_{aa'}^{e}$,
respectively,
we have chosen.
Then for $a_{1}, a_{2}, a_{3}, a_{4}\in \A$, 
\begin{eqnarray*}
\lefteqn{\sum_{a_{5}\in \A}\sum_{p=1}^{N_{a_{1}a_{5}}^{a_{4}}}
\sum_{q=1}^{N_{a_{2}a_{3}}^{a_{5}}}
F(\Y_{a_{1}a_{5}; p}^{a_{4}}\otimes \Y_{a_{2}a_{3}; q}^{a_{5}};
\Y_{a_{6}a_{3}; m}^{a_{4}}\otimes \Y_{a_{1}a_{2}; k}^{a_{6}})\cdot}\nn
&&\quad\quad\quad\quad\quad\quad\cdot
F(\Y_{a'_{1}a'_{5}; p}^{\prime; a'_{4}}\otimes \Y_{a'_{2}a'_{3}; q}^{\prime;
a_{5}};
\Y_{a'_{7}a'_{3}; n}^{\prime; a'_{4}}\otimes \Y_{a'_{1}a'_{2}; l}^{\prime; 
a'_{7}})\nn
&&=\delta_{a_{6}a_{7}}\delta_{mn}\delta_{kl}.
\end{eqnarray*}
\end{prop}
\pf 
For $a_{1}, a_{2}, a_{3}, a_{4}\in \A$, $w_{a_{i}}\in W^{a_{i}}$
and $w'_{a_{i}}\in (W^{a_{i}})'$  satisfying 
$\langle w'_{a_{i}}, w_{a_{i}}\rangle=1$
for $i=1, 2, 3$, using (\ref{inner}) and Lemma \ref{zero},
we have
\begin{eqnarray}\label{fusing-1}
\lefteqn{\res_{1-z_{1}-z_{3}=0\;|\; z_{3}}\res_{1-z_{2}-z_{4}=0\;|\;z_{4}}
(1-z_{1}-z_{3})^{-1}(1-z_{2}-z_{4})^{-1}}\nn
&&\quad\quad 
E(\langle e^{L(1)}\Y_{a'_{1}a'_{5}; k}^{\prime; a'_{4}}
((1-z_{1}-z_{3})^{L(0)}\tilde{w}'_{a_{1}}, z_{1})\cdot\nn
&&\quad\quad\quad\quad\quad\quad \quad\quad\quad \quad\quad\quad \cdot
\Y_{a'_{2}a'_{3};l}^{\prime; a'_{5}}
((1-z_{2}-z_{4})^{L(0)}\tilde{w}'_{a_{2}; q}, z_{2})\tilde{w}'_{a_{3}},\nn
&&\quad\quad\quad
e^{L(1)}\Y_{a_{1}a_{6}; m}^{a'_{4}}((1-z_{1}-z_{3})^{L(0)}\tilde{w}_{a_{1}}, z_{3})
\cdot\nn
&&\quad\quad\quad\quad\quad\quad \quad\quad\quad \quad\quad\quad \cdot
\Y_{a_{2}a_{3};n}^{a_{6}}
((1-z_{2}-z_{4})^{L(0)}\tilde{w}_{a_{2}; q}, z_{4})\tilde{w}_{a_{3}}\rangle)\nn
&&=\delta_{a_{5}a_{6}}\delta_{km}\res_{1-z_{2}-z_{4}=0\;|\;z_{4}}
(1-z_{2}-z_{4})^{-1}\nn
&&\quad\quad
E(\langle\Y_{a'_{2}a'_{3};l}^{\prime; a'_{5}}
((1-z_{2}-z_{4})^{L(0)}\tilde{w}'_{a_{2}}, z_{2})
\tilde{w}'_{a_{3}},\nn
&&\quad\quad\quad\quad\quad\quad \quad\quad\quad \quad\quad\quad 
\Y_{a_{2}a_{3};n}^{a_{5}}
((1-z_{2}-z_{4})^{L(0)}\tilde{w}_{a_{2}}, z_{4})\tilde{w}_{a_{3}}\rangle)\nn
&&=\delta_{a_{5}a_{6}}\delta_{km}\delta_{ln}.
\end{eqnarray}

On the other hand, by 
the associativity of intertwining operators and
Lemma \ref{zero}, we have
\begin{eqnarray}\label{fusing-2}
\lefteqn{\res_{1-z_{1}-z_{3}=0\;|\;z_{3}}\res_{1-z_{2}-z_{4}=0\;|\;z_{4}}
(1-z_{1}-z_{3})^{-1}(1-z_{2}-z_{4})^{-1}}\nn
&&\quad\quad E(\langle e^{L(1)}\Y_{a'_{1}a'_{5}; k}^{\prime; a'_{4}}
((1-z_{1}-z_{3})^{L(0)}\tilde{w}'_{a_{1}}, z_{1})\cdot\nn
&&\quad\quad\quad\quad\quad\quad \quad\quad\quad \quad\quad\quad \cdot
\Y_{a'_{2}a'_{3};l}^{\prime; a'_{5}}((1-z_{2}-z_{4})^{L(0)}
\tilde{w}'_{a_{2}}, z_{2})\tilde{w}'_{a_{3}},
\nn
&&\quad\quad\quad
e^{L(1)}\Y_{a_{1}a_{6}; m}^{a'_{4}}
((1-z_{1}-z_{3})^{L(0)}\tilde{w}_{a_{1}}, z_{3})\cdot\nn
&&\quad\quad\quad\quad\quad\quad \quad\quad\quad \quad\quad\quad \cdot
\Y_{a_{2}a_{3};n}^{a_{6}}((1-z_{2}-z_{4})^{L(0)}
\tilde{w}_{a_{2}}, z_{4})\tilde{w}_{a_{3}}\rangle)\nn
&&=\sum_{a_{7}, a_{8}\in \A}\sum_{i, j, s, t}
F(\Y_{a'_{1}a'_{5}; k}^{\prime; a'_{4}}\otimes 
\Y_{a'_{2}a'_{3};l}^{\prime; a'_{5}};
\Y_{a'_{7}a'_{3};i}^{\prime; a'_{4}}\otimes 
\Y_{a'_{1}a'_{2};j}^{\prime; a'_{7}})\cdot \nn
&&\quad\quad\quad\quad\quad\quad\quad
\cdot
F(\Y_{a_{1}a_{6}; m}^{a_{4}}\otimes \Y_{a_{2}a_{3};n}^{a_{6}};
\Y_{a_{8}a_{3};s}^{a_{4}}\otimes \Y_{a_{1}a_{2};t}^{a_{8}})\cdot\nn
&&\quad\quad\cdot\res_{1-z_{1}-z_{3}=0\;|\;z_{3}}\res_{1-z_{2}-z_{4}=0\;|\;
z_{4}}
(1-z_{1}-z_{3})^{-1}(1-z_{2}-z_{4})^{-1} \nn
&&\quad\quad\quad
E(\langle e^{L(1)}\Y_{a'_{7}a'_{3}; i}^{\prime; a'_{4}}(
\Y_{a'_{1}a'_{2};j}^{\prime; a'_{7}}((1-z_{1}-z_{3})^{L(0)}
\tilde{w}'_{a_{1}}, z_{1}-z_{2})\cdot\nn
&&\quad\quad\quad\quad\quad\quad \quad\quad\quad \quad\quad\quad
\quad\quad \cdot
(1-z_{2}-z_{4})^{L(0)}\tilde{w}'_{a_{2}}, z_{2})\tilde{w}'_{a_{3}}, \nn
&&\quad\quad\quad\quad
e^{L(1)}\Y_{a_{8}a_{3}; s}^{a_{4}}
(\Y_{a_{1}a_{2};t}^{a_{8}}((1-z_{1}-z_{3})^{L(0)}\tilde{w}_{a_{1}}, 
z_{3}-z_{4})\cdot\nn
&&\quad\quad\quad\quad\quad\quad \quad\quad\quad \quad\quad\quad
\quad\quad \cdot
(1-z_{2}-z_{4})^{L(0)}\tilde{w}_{a_{2}}, z_{4})\tilde{w}_{a_{3}}\rangle)\nn
&&=\sum_{a_{7}, a_{8}\in \A}\sum_{i, j, s, t}
F(\Y_{a'_{1}a'_{5}; k}^{\prime; a'_{4}}\otimes 
\Y_{a'_{2}a'_{3};l}^{\prime; a'_{5}};
\Y_{a'_{7}a'_{3};i}^{\prime; a'_{4}}\otimes 
\Y_{a'_{1}a'_{2};j}^{\prime; a'_{7}})\cdot \nn
&&\quad\quad\quad\quad\quad\quad\quad
\cdot
F(\Y_{a_{1}a_{6}; m}^{a_{4}}\otimes \Y_{a_{2}a_{3};n}^{a_{6}};
\Y_{a_{8}a_{3};s}^{a_{4}}\otimes \Y_{a_{1}a_{2};t}^{a_{8}})\cdot\nn
&&\quad\quad\cdot\res_{1-z_{2}-z_{4}=0\;|\;z_{4}}
\res_{1-z_{1}-z_{3}=0\;|\;z_{3}}
(1-z_{1}-z_{3})^{-1}(1-z_{2}-z_{4})^{-1} \nn
&&\quad\quad\quad
E\Biggl(\Biggl\langle e^{L(1)}\Y_{a'_{7}a'_{3}; i}^{\prime; a'_{4}}\Biggl(
(1-z_{2}-z_{4})^{L(0)}
\Y_{a'_{1}a'_{2};j}^{\prime; a'_{7}}
\Biggl(\Biggl(\frac{1-z_{1}-z_{3}}{1-z_{2}-z_{4}}\Biggr)^{L(0)}
\cdot\nn
&&\quad\quad\quad\quad\quad\quad\quad\quad\quad\quad\quad\quad\quad \cdot
\tilde{w}'_{a_{1}}, \frac{z_{1}-z_{2}}{1-z_{2}-z_{4}}\Biggr)
\tilde{w}'_{a_{2}}, z_{2}\Biggr)\tilde{w}'_{a_{3}}, \nn
&&\quad\quad\quad\quad
e^{L(1)}\Y_{a_{8}a_{3}; s}^{a_{4}}
\Biggl((1-z_{2}-z_{4})^{L(0)}
\Y_{a_{1}a_{2};t}^{a_{8}}
\Biggl(\Biggl(\frac{1-z_{1}-z_{3}}{1-z_{2}-z_{4}}\Biggr)^{L(0)}\cdot\nn
&&\quad\quad\quad\quad\quad\quad\quad\quad\quad\quad\quad\quad\quad \cdot
\tilde{w}_{a_{1}}, 
\frac{z_{3}-z_{4}}{1-z_{2}-z_{4}}\Biggr)
\tilde{w}_{a_{2}}, z_{4}\Biggr)\tilde{w}_{a_{3}}\Biggr\rangle\Biggr).
\end{eqnarray}

We now change the variables $z_{1}$ and $z_{3}$ to 
$$z_{5}=\frac{z_{1}-z_{2}}{1-z_{2}-z_{4}}$$
and 
$$z_{6}=\frac{z_{3}-z_{4}}{1-z_{2}-z_{4}}.$$
Then 
\begin{eqnarray*}
1-z_{5}-z_{6}&=&\frac{1-z_{1}-z_{3}}{1-z_{2}-z_{4}},\\
z_{3}&=&(1-z_{2}-z_{4})z_{6}+z_{4}.
\end{eqnarray*}
For  any branch $f(z_{1}, z_{2}, z_{3}, z_{4})$
of a multivalued analytic function 
of $z_{1}$, $z_{2}$, $z_{3}$ and $z_{4}$
on a suitable region $A$ such that it is equal to the restriction 
to $A\cap B$ of
a branch of the same analytic function on a region $B$ containing
the point $1-z_{1}-z_{3}=0$,
by definition, we have
\begin{eqnarray}\label{residue2}
\lefteqn{\res_{1-z_{1}-z_{3}=0\;|\;z_{3}}E(f(z_{1}, z_{2}, z_{3}, z_{4}))}\nn
&&=\res_{1-z_{5}-z_{6}=0\;|\;(1-z_{2}-z_{4})z_{6}+z_{4}}
E(f(z_{1}, z_{2}, z_{3}, z_{4}))\frac{1-z_{1}-z_{3}}{1-z_{5}-z_{6}}.
\end{eqnarray}
By (\ref{residue}), we have 
\begin{eqnarray}\label{residue3}
\lefteqn{\res_{1-z_{5}-z_{6}=0\;|\;(1-z_{2}-z_{4})z_{6}+z_{4}}
E(f(z_{1}, z_{2}, z_{3}, z_{4}))\frac{1-z_{1}-z_{3}}{1-z_{5}-z_{6}}}\nn
&&=\res_{1-z_{5}-z_{6}=0\;|\;z_{6}}
E(f(z_{1}, z_{2}, z_{3}, z_{4}))\frac{1-z_{1}-z_{3}}{1-z_{5}-z_{6}}.
\end{eqnarray}
From (\ref{residue2}) and (\ref{residue3}), we obtain
\begin{eqnarray}\label{residue4}
\lefteqn{\res_{1-z_{1}-z_{3}=0\;|\;z_{3}}E(f(z_{1}, z_{2}, z_{3}, z_{4}))}\nn
&&=\res_{1-z_{5}-z_{6}=0\;|\;z_{6}}
E(f(z_{1}, z_{2}, z_{3}, z_{4}))\frac{1-z_{1}-z_{3}}{1-z_{5}-z_{6}}.
\end{eqnarray}

Using (\ref{residue4}), the definition of the pairings
$\langle\cdot, 
\cdot\rangle_{\V_{a_{1}a_{2}}^{a_{3}}}$,
Lemma \ref{zero}, and the fact that 
for $a_{1}, a_{2}, a_{3}\in \A$, 
$\{\Y_{a'_{1}a'_{2}; i}^{\prime; a'_{3}}\;|\;i=1, \dots, 
N_{a'_{1}a'_{2}}^{a'_{3}}\}$ are the dual bases of 
$\{\Y_{a_{1}a_{2}; i}^{a_{3}}\;|\;i=1, \dots, N_{a_{1}a_{2}}^{a_{3}}\}$
with respect to the pairing $\langle\cdot, 
\cdot\rangle_{\V_{a_{1}a_{2}}^{a_{3}}}$, we see that the right-hand side of 
(\ref{fusing-2}) is equal to
\begin{eqnarray}\label{fusing-3}
\lefteqn{\sum_{a_{7}, a_{8}\in \A}\sum_{i, j, s, t}
F(\Y_{a'_{1}a'_{5}; k}^{\prime; a'_{4}}\otimes 
\Y_{a'_{2}a'_{3};l}^{\prime; a'_{5}};
\Y_{a'_{7}a'_{3};i}^{\prime; a'_{4}}\otimes 
\Y_{a'_{1}a'_{2};j}^{\prime; a'_{7}})\cdot} \nn
&&\quad\quad\quad\quad\quad\quad\quad
\cdot
F(\Y_{a_{1}a_{6}; m}^{a_{4}}\otimes \Y_{a_{2}a_{3};n}^{a_{6}};
\Y_{a_{8}a_{3};s}^{a_{4}}\otimes \Y_{a_{1}a_{2};t}^{a_{8}})\cdot\nn
&&\quad\quad\cdot\res_{1-z_{2}-z_{4}=0\;|\;z_{4}}
\res_{1-z_{5}-z_{6}=0\;|\;z_{6}}
(1-z_{5}-z_{6})^{-1}(1-z_{2}-z_{4})^{-1} \nn
&&\quad\quad\quad
E(\langle e^{L(1)}\Y_{a'_{7}a'_{3}; i}^{\prime; a'_{4}}(
(1-z_{2}-z_{4})^{L(0)}
\cdot\nn
&&\quad\quad\quad\quad\quad\quad\quad\quad\quad\quad \cdot
\Y_{a'_{1}a'_{2};j}^{\prime; a'_{7}}
((1-z_{5}-z_{6})^{L(0)}
\tilde{w}'_{a_{1}}, z_{5})
\tilde{w}'_{a_{2}}, z_{2})\tilde{w}'_{a_{3}}, \nn
&&\quad\quad\quad\quad
e^{L(1)}\Y_{a_{8}a_{3}; s}^{a_{4}}
((1-z_{2}-z_{4})^{L(0)}
\cdot\nn
&&\quad\quad\quad\quad\quad\quad\quad\quad\quad\quad \cdot
\Y_{a_{1}a_{2};t}^{a_{8}}
((1-z_{5}-z_{6})^{L(0)}
\tilde{w}_{a_{1}}, 
z_{6})
\tilde{w}_{a_{2}}, z_{4})\tilde{w}_{a_{3}}\rangle)\nn
&&=\sum_{a_{7}, a_{8}\in \A}\sum_{i, j, s, t}
F(\Y_{a'_{1}a'_{5}; k}^{\prime; a'_{4}}\otimes 
\Y_{a'_{2}a'_{3};l}^{\prime; a'_{5}};
\Y_{a'_{7}a'_{3};i}^{\prime; a'_{4}}\otimes 
\Y_{a'_{1}a'_{2};j}^{\prime; a'_{7}})\cdot \nn
&&\quad\quad\quad\quad\quad\quad\quad
\cdot
F(\Y_{a_{1}a_{6}; m}^{a_{4}}\otimes \Y_{a_{2}a_{3};n}^{a_{6}};
\Y_{a_{7}a_{3};s}^{a_{4}}\otimes \Y_{a_{1}a_{2};t}^{a_{7}})\cdot\nn
&&\quad\quad\cdot
\res_{1-z_{5}-z_{6}=0\;|\;z_{6}}
\res_{1-z_{2}-z_{4}=0\;|\;z_{4}}
(1-z_{5}-z_{6})^{-1}(1-z_{2}-z_{4})^{-1} \nn
&&\quad\quad\quad
E(\langle e^{L(1)}\Y_{a'_{7}a'_{3}; i}^{\prime; a'_{4}}(
(1-z_{2}-z_{4})^{L(0)}
\cdot\nn
&&\quad\quad\quad\quad\quad\quad\quad\quad\quad\quad \cdot
\Y_{a'_{1}a'_{2};j}^{\prime; a'_{7}}
((1-z_{5}-z_{6})^{L(0)}
\tilde{w}'_{a_{1}}, z_{5})
\tilde{w}'_{a_{2}}, z_{2})\tilde{w}'_{a_{3}}, \nn
&&\quad\quad\quad\quad
e^{L(1)}\Y_{a_{8}a_{3}; s}^{a_{4}}
((1-z_{2}-z_{4})^{L(0)}
\cdot\nn
&&\quad\quad\quad\quad\quad\quad\quad\quad\quad\quad \cdot
\Y_{a_{1}a_{2};t}^{a_{7}}
((1-z_{5}-z_{6})^{L(0)}
\tilde{w}_{a_{1}}, 
z_{6})
\tilde{w}_{a_{2}}, z_{4})\tilde{w}_{a_{3}}\rangle)\nn
&&=\sum_{a_{7}\in \A}\sum_{i, j, s, t}
F(\Y_{a'_{1}a'_{5}; k}^{\prime; a'_{4}}\otimes 
\Y_{a'_{2}a'_{3};l}^{\prime; a'_{5}};
\Y_{a'_{7}a'_{3};i}^{\prime; a'_{4}}\otimes 
\Y_{a'_{1}a'_{2};j}^{\prime; a'_{7}})\cdot \nn
&&\quad\quad\quad\quad\quad\quad\quad
\cdot
F(\Y_{a_{1}a_{6}; m}^{a_{4}}\otimes \Y_{a_{2}a_{3};n}^{a_{6}};
\Y_{a_{7}a_{3};s}^{a_{4}}\otimes \Y_{a_{1}a_{2};t}^{a_{7}})\cdot\nn
&&\quad\quad\cdot\res_{1-z_{5}-z_{6}=0\;|\;z_{6}}
(1-z_{5}-z_{6})^{-1} \delta_{is}\nn
&&\quad\quad\quad
E(\langle e^{L(1)}\Y_{a'_{1}a'_{2};j}^{\prime; a'_{7}}
((1-z_{5}-z_{6})^{L(0)}
\tilde{w}'_{a_{1}}, z_{5})
\tilde{w}'_{a_{2}}, \nn
&&\quad\quad\quad\quad\quad\quad\quad\quad\quad\quad\quad\quad
e^{L(1)}\Y_{a_{1}a_{2};t}^{a_{8}}
((1-z_{5}-z_{6})^{L(0)}
\tilde{w}_{a_{1}}, 
z_{6})
\tilde{w}_{a_{2}}\rangle)\nn
&&=\sum_{a_{7}\in \A}\sum_{i, j, s, t}
F(\Y_{a'_{1}a'_{5}; k}^{\prime; a'_{4}}\otimes 
\Y_{a'_{2}a'_{3};l}^{\prime; a'_{5}};
\Y_{a'_{7}a'_{3};i}^{\prime; a'_{4}}\otimes 
\Y_{a'_{1}a'_{2};j}^{\prime; a'_{7}})\cdot \nn
&&\quad\quad\quad\quad\quad\quad\quad
\cdot
F(\Y_{a_{1}a_{6}; m}^{a_{4}}\otimes \Y_{a_{2}a_{3};n}^{a_{6}};
\Y_{a_{8}a_{3};s}^{a_{4}}\otimes \Y_{a_{1}a_{2};t}^{a_{8}})\delta_{is}\delta_{jt}\nn
&&=\sum_{a_{7}\in \A}\sum_{i, j}
F(\Y_{a'_{1}a'_{5}; k}^{\prime; a'_{4}}\otimes 
\Y_{a'_{2}a'_{3};l}^{\prime; a'_{5}};
\Y_{a'_{7}a'_{3};i}^{\prime; a'_{4}}\otimes 
\Y_{a'_{1}a'_{2};j}^{\prime; a'_{7}})\cdot \nn
&&\quad\quad\quad\quad\quad\quad\quad
\cdot
F(\Y_{a_{1}a_{6}; m}^{a_{4}}\otimes \Y_{a_{2}a_{3};n}^{a_{6}};
\Y_{a_{7}a_{3};i}^{a_{4}}\otimes \Y_{a_{1}a_{2};j}^{a_{7}}).
\end{eqnarray}
From (\ref{fusing-1})--(\ref{fusing-3}), we see that 
the right inverse of the matrix with entries 
$$F(\Y_{a'_{1}a'_{5}; k}^{\prime; a'_{4}}\otimes 
\Y_{a'_{2}a'_{3};l}^{\prime; a'_{5}};
\Y_{a'_{7}a'_{3};i}^{\prime; a'_{4}}\otimes 
\Y_{a'_{1}a'_{2};j}^{\prime; a'_{7}})$$
is the transpose of the matrix with entries
$$F(\Y_{a_{1}a_{6}; m}^{a_{4}}\otimes \Y_{a_{2}a_{3};n}^{a_{6}};
\Y_{a_{7}a_{3};i}^{a_{4}}\otimes \Y_{a_{1}a_{2};j}^{a_{7}}).$$
Since for square matrices, right inverses are also left inverses,
the proposition is proved.
\epfv

For $a\in \A$, we use $\sqrt{F_{a}}$ to denote the 
square root $\sqrt{|F_{a}|}e^{\frac{i\arg F_{a}}{2}}$
of $F_{a}$. 
For $a_{1}, a_{2}, a_{3}\in \A$,
consider the modified pairings 
$$\frac{\sqrt{F_{a_{3}}}}{\sqrt{F_{a_{1}}}\sqrt{F_{a_{2}}}}
\langle \cdot, \cdot\rangle_{\V_{a_{1}a_{2}}^{a_{3}}}.$$
These pairings give a nondegenerate bilinear form 
$(\cdot, \cdot)_{\V}$ on $\V$. 
For any $\sigma\in S_{3}$, $\{\sigma(\Y_{a_{1}a_{2}; i}^{a_{3}})\;
|\;i=\dots, N_{a_{1}a_{2}}^{a_{3}}\}$ is a basis of 
$\sigma(\V_{a_{1}a_{2}}^{a_{3}})$.

We have:

\begin{prop}\label{skew}
The nondegenerate bilinear form $(\cdot, \cdot)_{\V}$
is invariant with respect to the action of $S_{3}$ on $\V$, that is,
for $a_{1}, a_{2}, a_{3}\in \A$, $\sigma\in S_{3}$, $\Y_{1}\in 
\V_{a_{1}a_{2}}^{a_{3}}$
and $\Y_{2}\in \V_{a'_{1}a'_{2}}^{a'_{3}}$, 
$$(\sigma(\Y_{1}), \sigma(\Y_{2}))_{\V}= 
(\Y_{1}, \Y_{2})_{\V}.$$
Equivalently, for $a_{1}, a_{2}, a_{3}\in \A$,
$$\left\{\left.\frac{\sqrt{F_{a_{1}}}\sqrt{F_{a_{2}}}
\sqrt{F_{a_{\sigma^{-1}(3)}}}}
{\sqrt{F_{a_{\sigma^{-1}(1)}}}\sqrt{F_{a_{\sigma^{-1}(2)}}}\sqrt{F_{a_{3}}}}
\sigma(\Y_{a'_{1}a'_{2}; i}^{\prime; a'_{3}})\;\right|\;i=\dots, 
N_{a_{1}a_{2}}^{a_{3}}\right\}$$
is the dual basis of $\{\sigma(\Y_{a_{1}a_{2}; i}^{a_{3}})\;
|\;i=\dots, N_{a_{1}a_{2}}^{a_{3}}\}$.
\end{prop}
\pf
The equivalence of the first conclusion and 
the second conclusion is clear. 

We first prove the result for 
$\sigma=\sigma_{12}$. In this case, we need to show that 
$\{\sigma_{12}(\Y_{a'_{1}a'_{2}; j}^{a'_{3}})\;|\;
i=\dots, N_{a_{1}a_{2}}^{a_{3}}\}$ is the dual basis 
of $\{\sigma_{12}(\Y_{a_{1}a_{2}; i}^{a_{3}})\;
|\;i=\dots, N_{a_{1}a_{2}}^{a_{3}}\}$.
For $i, j=1, \dots, N_{a_{1}a_{2}}^{a_{3}}$, by (\ref{inner-fusing}),
we have
\begin{equation}\label{skew-1}
\langle \sigma_{12}(\Y_{a'_{1}a'_{2}; i}^{\prime; a'_{3}}),
\sigma_{12}(\Y_{a_{1}a_{2}; j}^{a_{3}})\rangle=
F(\sigma_{23}(\sigma_{12}(\Y_{a'_{1}a'_{2}; i}^{\prime; a'_{3}}))
\otimes \sigma_{12}(\Y_{a_{1}a_{2}; j}^{a_{3}}); \Y_{ea_{1}; 1}^{a_{1}}\otimes
\Y_{a'_{2}a_{2}; 1}^{e}).
\end{equation}
By Proposition 3.4 in \cite{H7}, the right-hand side of (\ref{skew})
is equal to
\begin{eqnarray}\label{skew-2}
&F(\sigma_{132}(\sigma_{12}(\Y_{a_{1}a_{2}; j}^{a_{3}}))
\otimes \sigma_{123}(\sigma_{23}(\sigma_{12}
(\Y_{a'_{1}a'_{2}; i}^{\prime; a'_{3}}))); 
\sigma_{123}(\Y_{a'_{2}a_{2}; 1}^{e})\otimes
\sigma_{132}(\Y_{ea_{1}; 1}^{a_{1}}))&\nn
&=F(\sigma_{23}(\Y_{a_{1}a_{2}; j}^{a_{3}})
\otimes \Y_{a'_{1}a'_{2}; i}^{\prime; a'_{3}}; 
\Y_{ea'_{2}; 1}^{a'_{2}}\otimes
\Y_{a_{1}a'_{1}; 1}^{e})&\nn
&\;\;=F(\sigma_{23}(\Y_{a''_{1}a''_{2}; j}^{a''_{3}})
\otimes \Y_{a'_{1}a'_{2}; i}^{\prime; a'_{3}}; 
\Y_{ea'_{2}; 1}^{a'_{2}}\otimes
\Y_{a''_{1}a'_{1}; 1}^{e}).&
\end{eqnarray}
By (\ref{inner-fusing}) again, the right-hand side of (\ref{skew-2}) is 
equal to 
$$\langle \Y_{a''_{1}a''_{2}; j}^{a''_{3}}, 
\Y_{a'_{1}a'_{2}; i}^{\prime; a'_{3}}\rangle
=\delta_{ij},$$
proving the case of $\sigma=\sigma_{12}$. 

Next we prove the result for $\sigma=\sigma_{23}$.
We need to find the relation between
the matrices
$$\langle \sigma_{23}(\Y_{a'_{1}a'_{2}; i}^{a'_{3}}),
\sigma_{23}(\Y_{a_{1}a_{2}'; j}^{a_{3}})
\rangle_{\V_{a_{1}a'_{3}}^{a'_{2}}}$$
and 
$$\langle \Y_{a'_{1}a'_{2}; i}^{a'_{3}},
\Y_{a_{1}a_{2}'; j}^{a_{3}}\rangle_{\V_{a_{1}a_{2}}^{a_{3}}}.$$
By definition, we need to find the relation between the 
matrices 
\begin{equation}\label{form-matrix}
F(\Y_{a'_{1}a'_{2}; i}^{a'_{3}}
\otimes \sigma_{23}(\Y_{a_{1}a_{2}; j}^{a_{3}}); 
\Y_{ea'_{3}; 1}^{a'_{3}}\otimes
\Y_{a'_{1}a_{1}; 1}^{e})
\end{equation}
and
\begin{equation}\label{form-matrix-1}
F(\sigma_{23}(\Y_{a'_{1}a'_{2}; i}^{a'_{3}; (2)})
\otimes \Y_{a_{1}a_{2}; j}^{a_{3}; (1)}; \Y_{ea_{2}; 1}^{a_{2}}\otimes
\Y_{a'_{1}a_{1}; 1}^{e}).
\end{equation}
From (4.9) in \cite{H7} (or (\ref{formula1})), we see that 
the inverse of the matrix (\ref{form-matrix}) is 
\begin{equation}\label{form-matrix-inv}
\frac{F(\Y_{a'_{3}e; 1}^{a'_{3}}\otimes \Y_{a_{2}a'_{2}; 1}^{e}; 
\Y_{a'_{1}a'_{2}; k}^{a'_{3}}\otimes 
\sigma_{132}(\sigma_{23}(\Y_{a_{1}a_{2}; j}^{a_{3}})))}
{F_{a'_{3}}}.
\end{equation}
By Proposition 3.4 in \cite{H7} and the fact 
$F_{a'}=F_{a}$ for $a\in \A$, (\ref{form-matrix-inv})
is equal to 
\begin{eqnarray*}
\lefteqn{\frac{F(\Y_{a'_{2}e; 1}^{a'_{2}}\otimes \Y_{a_{3}a'_{3}; 1}^{e}; 
\sigma_{23}(\Y_{a_{1}a_{2}; j}^{a_{3}}) \otimes 
\sigma_{132}(\Y_{a'_{1}a'_{2}; k}^{a'_{3}}))}
{F_{a_{3}}}}\nn
&&=\frac{F_{a_{2}}}{F_{a_{3}}}
\frac{F(\Y_{a'_{2}e; 1}^{a'_{2}}\otimes \Y_{a_{3}a'_{3}; 1}^{e}; 
\sigma_{23}(\Y_{a_{1}a_{2}; j}^{a_{3}}) \otimes 
\sigma_{132}(\Y_{a'_{1}a'_{2}; k}^{a'_{3}}))}
{F_{a'_{2}}},
\end{eqnarray*}
which by (4.9) in \cite{H7} (or (\ref{formula1})) again is equal 
to $\frac{F_{a_{2}}}{F_{a_{3}}}$ times the inverse of 
(\ref{form-matrix-1}). So the inverse of the matrix (\ref{form-matrix})
is equal to $\frac{F_{a_{2}}}{F_{a_{3}}}$ times the inverse of 
(\ref{form-matrix-1}). Thus the matrix (\ref{form-matrix})
is equal to $\frac{F_{a_{3}}}{F_{a_{2}}}$ times the matrix
(\ref{form-matrix-1}), or equivalently, 
$$(\sigma_{23}(\Y_{a_{1}a_{2}; j}^{a_{3}}), 
\sigma_{23}(\Y_{a'_{1}a'_{2}; k}^{a'_{3}}))_{\V}
=(\Y_{a_{1}a_{2}; j}^{a_{3}}, 
\Y_{a'_{1}a'_{2}; k}^{a'_{3}})_{\V}.$$

Since $S_{3}$ is 
generated by $\sigma_{12}$
and $\sigma_{23}$, the conclusion of the proposition follows.
\epfv

We are ready to construct a full field algebra using 
the bases of intertwining operators we have chosen.
Let 
$$F=\oplus_{a\in \A} W^{a}\otimes W^{a'}.$$
For $w_{a_{1}}\in W^{a_{1}}$, $w_{a_{2}}\in W^{a_{2}}$, 
$w_{a'_{1}}\in W^{a'_{1}}$ and $w_{a'_{2}}\in W^{a'_{2}}$, 
we define 
\begin{eqnarray*}
\lefteqn{\mathbb{Y}((w_{a_{1}}\otimes w_{a'_{1}}), z, \zeta)
(w_{a_{2}}\otimes w_{a'_{2}})}\nn
&&=\sum_{a_{3}\in \A}\sum_{p=1}^{N_{a_{1}a_{2}}^{a_{3}}}
\Y_{a_{1}a_{2}; p}^{a_{3}}(w_{a_{1}}, z) w_{a_{2}}
\otimes \Y_{a'_{1}a'_{2}; p}^{\prime; a'_{3}}
(w_{a'_{1}}, \zeta)w_{a'_{2}}.
\end{eqnarray*}

\begin{thm}
The quadruple $(F, \mathbb{Y}, \one\otimes \one, \omega\otimes \one,
\one \otimes \omega)$ is a conformal full field algebra over 
$V\otimes V$. 
\end{thm}
\pf 
The identity property, the creation property and the single-valuedness
property are clear. We prove the associativity and the skew-symmetry here.

We prove associativity first.
For $a_{1}, a_{2}\in \A$, $w_{a_{1}}\in W^{a_{1}}$, $w_{a_{2}}\in W^{a_{2}}$, 
$w'_{a_{1}}\in (W^{a_{1}})'$, $w'_{a_{2}}\in (W^{a_{2}})'$, 
using the associativity of intertwining operators and Proposition 
\ref{fusing},
we have
\begin{eqnarray*}
\lefteqn{\BY((w_{a_{1}}\otimes w'_{a_{1}}), z_{1}, \zeta_{1})
\BY((w_{a_{2}}\otimes w'_{a_{2}}), z_{2}, \zeta_{2})}\nn
&&=\sum_{a_{3}, a_{4}, a_{5}\in\A}\sum_{p=1}^{N_{a_{1}a_{5}}^{a_{4}}}
\sum_{q=1}^{N_{a_{2}a_{3}}^{a_{5}}}(\Y_{a_{1}a_{5}; p}^{a_{4}}
(w_{a_{1}}, z_{1})
\Y_{a_{2}a_{3}; q}^{a_{5}}(w_{a_{2}}, z_{2}))\nn
&&\quad\quad\quad\quad\quad\quad\quad\quad\quad\quad\quad\quad
\quad\quad\otimes
(\Y_{a'_{1}a'_{5}; p}^{\prime; a'_{4}}(w'_{a_{1}}, \zeta_{1})
\Y_{a'_{2}a'_{3}; q}^{\prime; a'_{5}}(w'_{a_{2}}, \zeta_{2}))\nn
&&=\sum_{a_{3}, a_{4}, a_{5}, a_{6},a_{7}\in\A}
\sum_{p=1}^{N_{a_{1}a_{5}}^{a_{4}}}
\sum_{q=1}^{N_{a_{2}a_{3}}^{a_{5}}}\sum_{m=1}^{N_{a_{6}a_{3}}^{a_{5}}}
\sum_{n=1}^{N_{a_{1}a_{2}}^{a_{6}}}\sum_{k=1}^{N_{a'_{7}a'_{3}}^{a'_{5}}}
\sum_{l=1}^{N_{a'_{1}a'_{2}}^{a'_{7}}}\nn
&&\quad\quad\quad\quad\quad \cdot
F(\Y_{a_{1}a_{5}; p}^{a_{4}}\otimes \Y_{a_{2}a_{3}; q}^{a_{5}};
\Y_{a_{6}a_{3}; m}^{a_{4}}\otimes \Y_{a_{1}a_{2}; n}^{a_{6}})\cdot\nn
&&\quad\quad\quad\quad\quad\cdot
F(\Y_{a'_{1}a'_{5}; p}^{\prime; a'_{4}}\otimes 
\Y_{a'_{2}a'_{3}; q}^{\prime; a'_{5}};
\Y_{a'_{7}a'_{3}; k}^{a'_{4}}\otimes \Y_{a'_{1}a'_{2}; l}^{a'_{7}})\cdot\nn
&&\quad\quad\quad\quad\quad\quad
\cdot ((\Y_{a_{6}a_{3}; m}^{a_{4}}
(\Y_{a_{1}a_{2}; n}^{a_{6}}(w_{a_{1}}, z_{1}-z_{2})
w_{a_{2}}, z_{2}))\nn
&&\quad\quad\quad\quad\quad\quad\quad\quad\quad\quad\quad\quad
\otimes
(\Y_{a'_{7}a'_{3}; k}^{\prime; a'_{4}}
(\Y_{a'_{1}a'_{2}; l}^{\prime; a'_{7}}(w'_{a_{1}}, 
\zeta_{1}-\zeta_{2})w'_{a_{2}}, \zeta_{2})))\nn
&&=\sum_{a_{3}, a_{4}, a_{6}, a_{7}\in\A}\sum_{m=1}^{N_{a_{6}a_{3}}^{a_{5}}}
\sum_{n=1}^{N_{a_{1}a_{2}}^{a_{6}}}\sum_{k=1}^{N_{a'_{7}a'_{3}}^{a'_{5}}}
\sum_{l=1}^{N_{a'_{1}a'_{2}}^{a'_{7}}}
\delta_{a_{6}a_{7}}\delta_{mk}\delta_{nl}\cdot\nn
&&\quad\quad\quad\quad\quad\quad\quad\quad\quad
\cdot 
(\Y_{a_{6}a_{3}; m}^{a_{4}}(\Y_{a_{1}a_{2}; n}^{a_{6}}(w_{a_{1}}, z_{1}-z_{2})
w_{a_{2}}, z_{2}))\nn
&&\quad\quad\quad\quad\quad\quad\quad\quad\quad\quad\quad\quad
\quad\quad\otimes
(\Y_{a'_{7}a'_{3}; k}^{\prime; a'_{4}}
(\Y_{a'_{1}a'_{2}; l}^{\prime; a'_{7}}(w'_{a_{1}}, 
\zeta_{1}-\zeta_{2})w'_{a_{2}}, \zeta_{2}))\nn
&&=\sum_{a_{3}, a_{4}, a_{6}\in\A}\sum_{m=1}^{N_{a_{6}a_{3}}^{a_{5}}}
\sum_{n=1}^{N_{a_{1}a_{2}}^{a_{6}}}
(\Y_{a_{6}a_{3}; m}^{a_{4}}(\Y_{a_{1}a_{2}; n}^{a_{6}}(w_{a_{1}}, z_{1}-z_{2})
w_{a_{2}}, z_{2}))\nn
&&\quad\quad\quad\quad\quad\quad\quad\quad\quad\quad\quad\quad
\quad\quad\otimes
(\Y_{a'_{6}a'_{3}; m}^{\prime; a'_{4}}
(\Y_{a'_{1}a'_{2}; n}^{\prime; a'_{6}}(w'_{a_{1}}, 
\zeta_{1}-\zeta_{2})w'_{a_{2}}, \zeta_{2}))\nn
&&=\BY(\BY((w_{a_{1}}\otimes w'_{a_{1}}), z_{1}-z_{2}, \zeta_{1}-\zeta_{2})
(w_{a_{2}}\otimes w'_{a_{2}}), z_{2}, \zeta_{2}).
\end{eqnarray*}

We now prove the skew-symmetry. By Proposition \ref{skew},
$\{\sigma_{12}(\Y_{a'_{1}a'_{2}; i}^{\prime; a'_{3}})\;|\; i=1, \dots,
N_{a_{1}a_{2}}^{a_{3}}\}$ is the  dual basis of 
$\{\sigma_{12}(\Y_{a_{1}a_{2}; i}^{a_{3}})\;|\; i=1, \dots,
N_{a_{1}a_{2}}^{a_{3}}\}$. Thus
for $a_{1}, a_{2}\in \A$, $w_{a_{1}}\in W^{a_{1}}$, $w_{a_{2}}\in W^{a_{2}}$, 
$w'_{a_{1}}\in (W^{a_{1}})'$, $w'_{a_{2}}\in (W^{a_{2}})'$,  we have
\begin{eqnarray*}
\lefteqn{\BY((w_{a_{2}}\otimes w'_{a_{2}}), z, \zeta)
(w_{a_{1}}\otimes w'_{a_{1}})}\nn
&&=\sum_{a_{3}\in \A}\sum_{p=1}^{N_{a_{2}a_{1}}^{a_{3}}}
\sigma_{12}(\Y_{a_{1}a_{2}; p}^{a_{3}})(w_{a_{2}}, z) w_{a_{1}}
\otimes \sigma_{12}(\Y_{a'_{1}a'_{2}; p}^{\prime; a'_{3}})
(w_{a'_{2}}, \zeta)w_{a'_{1}}\nn
&&=\sum_{a_{3}\in \A}\sum_{p=1}^{N_{a_{1}a_{2}}^{a_{3}}}
e^{-\pi i\Delta(\Y_{a_{1}a_{2}; p}^{a_{3}})}
e^{zL(-1)}\Y_{a_{1}a_{2}; p}^{a_{3}}(w_{a_{1}}, e^{\pi i}z) w_{a_{2}}\nn
&&\quad\quad\quad\quad\quad\quad\quad\quad\quad\quad\quad\quad
\otimes e^{\pi i\Delta(\Y_{a'_{1}a'_{2}; p}^{a'_{3}})}
e^{\zeta L(-1)}\Y_{a'_{1}a'_{2}; p}^{\prime; a'_{3}})
(w_{a'_{1}}, e^{-\pi i}\zeta)w_{a'_{2}}\nn
&&=(e^{zL(-1)}\otimes e^{\zeta L(-1)})\cdot\nn
&&\quad\quad\quad\quad\quad\sum_{a_{3}\in \A}
\sum_{p=1}^{N_{a_{1}a_{2}}^{a_{3}}}
\Y_{a_{1}a_{2}; p}^{a_{3}}(w_{a_{1}}, e^{\pi i}z) w_{a_{2}}
\otimes \Y_{a'_{1}a'_{2}; p}^{\prime; a'_{3}}
(w_{a'_{1}}, e^{-\pi i}\zeta)w_{a'_{2}}\nn
&&=(e^{zL(-1)}\otimes e^{\zeta L(-1)})\BY((w_{a_{1}}\otimes w'_{a_{1}}),
e^{\pi i}z, e^{-\pi i}\zeta)(w_{a_{2}}\otimes w'_{a_{2}}).
\end{eqnarray*}
\epfv

\begin{defn}
{\rm A nondegenerate bilinear form 
$(\cdot, \cdot)$ on a conformal full field algebra 
$$(F, m, \one, \omega^{L}, \omega^{R})$$
is said to be {\it invariant}
if for $u, v, w\in F$, 
\begin{eqnarray*}
\lefteqn{(\mathbb{Y}(u; z, \bar{z})v, w)}\nn
&&=(v, \mathbb{Y}(e^{zL^{L}(1)+\bar{z}L^{R}(1)}
e^{\pi iL^{L}(0)-\pi iL^{R}(0)}z^{-2L^{L}(0)}\bar{z}^{-2L^{R}(0)}
u; z^{-1}, \bar{z}^{-1})w).
\end{eqnarray*}}
\end{defn}

The conformal full field algebra $F$ we constructed above
has a natural nondegenerate 
bilinear form $(\cdot, \cdot)_{F}: F\otimes F\to \C$ given by
$$((w_{a_{1}}\otimes w'_{a_{1}}), (w_{a_{2}}\otimes w'_{a_{2}}))_{F}
=\left\{\begin{array}{ll}0&a_{1}\ne a'_{2}\\
F_{a_{1}}\langle w_{a_{1}}, w_{a_{2}}\rangle\langle w'_{a_{1}}, w'_{a_{2}}
\rangle&a_{1}=a'_{2}
\end{array}\right.$$
for $a_{1}, a_{2}\in \A$, $w_{a_{1}}\in W^{a_{1}}$, $w_{a_{2}}\in W^{a_{2}}$, 
$w'_{a_{1}}\in (W^{a_{1}})'$, $w'_{a_{2}}\in (W^{a_{2}})'$.
We have:

\begin{thm}
The nondegenerate bilinear form $(\cdot, \cdot)_{F}$ is invariant.
\end{thm}
\pf
For $a_{1}, a_{2}, a_{3}\in \A$, $w_{a_{1}}\in W^{a_{1}}$, 
$w_{a_{2}}\in W^{a_{2}}$, 
$w_{a_{3}}\in W^{a_{3}}$, 
$w'_{a_{1}}\in (W^{a_{1}})'$, $w'_{a_{2}}\in (W^{a_{2}})'$, 
$w'_{a_{3}}\in (W^{a_{3}})'$,
using Proposition \ref{skew} for the case $\sigma=\sigma_{23}$, we have
\begin{eqnarray*}
\lefteqn{(\Y((w_{a_{1}}\otimes w'_{a_{1}}), z, \zeta)(w_{a_{2}}
\otimes w'_{a_{2}}),
(w_{a_{3}}\otimes w'_{a_{3}}))_{F}}\nn
&&=\sum_{a_{4}\in \A}\sum_{p=1}^{N_{a_{1}a_{2}}^{a_{4}}}
((\Y_{a_{1}a_{2}; p}^{a_{4}}(w_{a_{1}}, z) w_{a_{2}}
\otimes \Y_{a'_{1}a'_{2}; p}^{\prime; a'_{4}}
(w'_{a_{1}}, \zeta)w'_{a_{2}}), (w_{a_{3}}\otimes w'_{a_{3}}))_{F}\nn
&&=\sum_{p=1}^{N_{a_{1}a_{2}}^{a_{3}}}
F_{a_{3}}
\langle \Y_{a_{1}a_{2}; p}^{a_{3}}(w_{a_{1}}, z) w_{a_{2}}, w_{a_{3}}\rangle
\langle \Y_{a'_{1}a'_{2}; p}^{\prime; a'_{3}}
(w'_{a_{1}}, \zeta)w'_{a_{2}}, w'_{a_{3}}\rangle\nn
&&=\sum_{p=1}^{N_{a_{1}a_{2}}^{a_{3}}}
F_{a_{3}}\langle \sigma_{23}(\sigma_{23}(\Y_{a_{1}a_{2}; p}^{a_{3}})
(w_{a_{1}}, z) w_{a_{2}}, w_{a_{3}}\rangle
\langle \sigma_{23}(\sigma_{23}(\Y_{a'_{1}a'_{2}; p}^{\prime; a'_{3}})
(w'_{a_{1}}, \zeta)w'_{a_{2}}, w'_{a_{3}}\rangle\nn
&&=\sum_{p=1}^{N_{a_{1}a_{2}}^{a_{3}}}
F_{a_{3}}\langle  w_{a_{2}}, e^{\pi ih_{a_{1}}}\sigma_{23}(\Y_{a_{1}a_{2}; p}^{a_{3}})
(e^{zL(1)}e^{-\pi i L(0)}z^{-2L(0)}w_{a_{1}}, z^{-1})w_{a_{3}}\rangle\cdot\nn
&&\quad\quad\quad\quad\quad\quad\cdot
\langle w'_{a_{2}}, e^{-\pi ih_{a_{1}}}
\sigma_{23}(\Y_{a'_{1}a'_{2}; p}^{\prime; a'_{3}})
(e^{\zeta L(1)}e^{\pi i L(0)}\zeta^{-2L(0)}w'_{a_{1}}, \zeta^{-1})
w'_{a_{3}}\rangle\nn
&&=\sum_{p=1}^{N_{a_{1}a_{2}}^{a_{3}}}
F_{a_{2}}\langle  w_{a_{2}}, \sigma_{23}(\Y_{a_{1}a_{2}; p}^{a_{3}})
(e^{zL(1)}e^{-\pi i L(0)}z^{-2L(0)}w_{a_{1}}, z^{-1})w_{a_{3}}\rangle\cdot\nn
&&\quad\quad\quad\quad\cdot
\left\langle w'_{a_{2}}, \left(\frac{F_{a_{3}}}{F_{a_{2}}}
\sigma_{23}(\Y_{a'_{1}a'_{2}; p}^{\prime; a'_{3}})\right)
(e^{\zeta L(1)}e^{\pi i L(0)}\zeta^{-2L(0)}w'_{a_{1}}, \zeta^{-1})
w'_{a_{3}}\right\rangle\nn
&&=\sum_{p=1}^{N_{a_{1}a_{2}}^{a_{3}}}
((w_{a_{2}}\otimes w'_{a_{2}}), (\sigma_{23}(\Y_{a_{1}a_{2}; p}^{a_{3}})
(e^{zL(1)}e^{-\pi i L(0)}z^{-2L(0)}w_{a_{1}}, z^{-1})w_{a_{3}}\nn
&&\quad\quad\quad\quad\quad\quad\quad\otimes 
\left(\frac{F_{a_{3}}}{F_{a_{2}}}
\sigma_{23}(\Y_{a'_{1}a'_{2}; p}^{\prime; a'_{3}})\right)
(e^{\zeta L(1)}e^{\pi i L(0)}\zeta^{-2L(0)}w'_{a_{1}}, \zeta^{-1})
w'_{a_{3}}))_{F}\nn
&&=\sum_{a_{4}\in \A}\sum_{p=1}^{N_{a_{1}a_{4}}^{a_{3}}}
((w_{a_{2}}\otimes w'_{a_{2}}), (\sigma_{23}(\Y_{a_{1}a_{4}; p}^{a_{3}})
(e^{zL(1)}e^{-\pi i L(0)}z^{-2L(0)}w_{a_{1}}, z^{-1})w_{a_{3}}\nn
&&\quad\quad\quad\quad\quad\quad\quad\otimes 
\left(\frac{F_{a_{3}}}{F_{a_{2}}}
\sigma_{23}(\Y_{a'_{1}a'_{4}; p}^{\prime; a'_{3}})\right)
(e^{\zeta L(1)}e^{\pi i L(0)}\zeta^{-2L(0)}w'_{a_{1}}, \zeta^{-1})
w'_{a_{3}}))_{F}\nn
&&=((w_{a_{2}}\otimes w'_{a_{2}}), \Y((e^{zL(1)}
e^{-\pi i L(0)}z^{-2L(0)}w_{a_{1}}
\nn
&&\quad\quad\quad\quad\quad\quad\quad
\otimes (e^{\zeta L(1)}e^{\pi i L(0)}\zeta^{-2L(0)}w'_{a_{1}}), z^{-1}, 
\zeta^{-1})
(w_{a_{3}}\otimes w'_{a_{3}}))_{F},
\end{eqnarray*}
proving the invariance.
\epfv

\noindent {\small \sc Department of Mathematics, Rutgers University,
110 Frelinghuysen Rd., Piscataway, NJ 08854-8019}

\noindent {\em E-mail address}: yzhuang@math.rutgers.edu

\vspace{1em}

\noindent {\small \sc Department of Mathematics, Rutgers University,
110 Frelinghuysen Rd., Piscataway, NJ 08854-8019} 

\noindent and 

\noindent {\small \sc Max Planck Institute for Mathematics
in the Sciences, Inselstrasse 22, D-04103, Leipzig, Germany} 
(current address)

\noindent {\em E-mail address}: lkong@mis.mpg.de

\end{document}